\documentclass[letterpaper,11pt]{article}

\usepackage{ucs}
\usepackage[utf8x]{inputenc}
\usepackage{graphicx}
\usepackage{amsfonts}
\usepackage{dsfont}
\usepackage{amssymb}
\usepackage{amsmath}
\usepackage{amsthm}
\usepackage{enumerate}
\usepackage{stmaryrd}
\usepackage{fullpage}
\usepackage{ifthen}
\usepackage{subfigure}
\usepackage{epic}
\usepackage{authblk}
\usepackage{textcomp}

\usepackage[hypertexnames=false,colorlinks=true,linkcolor=blue,citecolor=blue]{hyperref}
\usepackage[numbers,comma,square,sort&compress]{natbib}
\usepackage[letterpaper,text={7in,9in},centering]{geometry}

\usepackage{color}
\usepackage{titlesec}
\graphicspath{{eps/}{pdf/}}

\newcommand{\bqq}{\begin{equation}}
\newcommand{\eqq}{\end{equation}}
\newcommand{\bqs}{\begin{equation*}}
\newcommand{\eqs}{\end{equation*}}

\newcommand{\C}{\mathbb{C}}

\newcommand{\R}{\mathbb{R}} 
\newcommand{\Z}{\mathbb{Z}}

\newcommand{\A}{\mathcal{A}}
\newcommand{\B}{\mathcal{B}}
\newcommand{\cC}{\mathcal{C}}

\newcommand{\F}{\mathcal{F}}
\newcommand{\G}{\mathcal{G}}

\newcommand{\I}{\mathcal{I}}
\newcommand{\J}{\mathcal{J}}
\newcommand{\K}{\mathcal{K}}
\newcommand{\cL}{\mathcal{L}}
\newcommand{\M}{\mathcal{M}}
\newcommand{\cN}{\mathcal{N}}
\newcommand{\cO}{\mathcal{O}}

\newcommand{\cR}{\mathcal{R}}
\newcommand{\cS}{\mathcal{S}}
\newcommand{\T}{\mathcal{T}}
\newcommand{\U}{\mathcal{U}}
\newcommand{\V}{\mathcal{V}}

\newcommand{\X}{\mathcal{X}}
\newcommand{\Y}{\mathcal{Y}}
\newcommand{\cZ}{\mathcal{Z}}

\newcommand{\bA}{\mathbf{A}}
\newcommand{\bC}{\mathbf{C}}
\newcommand{\bD}{\mathbf{D}}
\newcommand{\bF}{\mathbf{F}}
\newcommand{\bG}{\mathbf{G}}
\newcommand{\bL}{\mathbf{L}}
\newcommand{\bQ}{\mathbf{Q}}
\newcommand{\bU}{\mathbf{U}}
\newcommand{\bW}{\mathbf{W}}

\newcommand{\bu}{\mathbf{u}}
\newcommand{\bv}{\mathbf{v}}
\newcommand{\bw}{\mathbf{w}}
\newcommand{\bh}{\mathbf{h}}
\newcommand{\bg}{\mathbf{g}}

\numberwithin{equation}{section}

\newtheorem{lem}{Lemma}[section]
\newtheorem{thm}{Theorem}
\newtheorem{prop}[lem]{Proposition}
\newtheorem{cor}[lem]{Corollary}
\newtheorem{rmk}[lem]{Remark}

\newenvironment{Hypothesis}[1]%
  {\begin{trivlist}\item[]{\bf Hypothesis #1 }\em}{\end{trivlist}}

\makeatletter\@addtoreset{figure}{section}\makeatother

\makeatletter \@addtoreset{equation}{section} \makeatother

\newenvironment{Proof}[1][.]%
 {\begin{trivlist}\item[]\textbf{Proof#1 }}%
 {\hspace*{\fill}$\rule{0.3\baselineskip}{0.35\baselineskip}$\end{trivlist}}
\graphicspath{{Figures/}}

\title{Existence of pulses in excitable media with nonlocal coupling \footnote{GF was partially supported by the National Science Foundation through grant NSF-DMS-1311414. AS was partially supported by the National Science Foundation through grants NSF- DMS-0806614 and NSF-DMS-1311740.}}
\author[1]{Gr\'egory Faye}
\author[2]{Arnd Scheel}
\affil[1,2]{\small University of Minnesota,
School of Mathematics,
206 Church Street S.E.,
Minneapolis, MN 55455, USA}

\begin{document}
\maketitle

\begin{abstract}
We prove the existence of fast traveling pulse solutions in excitable media with non-local coupling. Existence results had been known, until now, in the case of local, diffusive coupling  and in the case of a discrete medium, with finite-range, non-local coupling. Our approach replaces methods from geometric singular perturbation theory, that had been crucial in previous existence proofs, by a PDE oriented approach, relying on exponential weights, Fredholm theory, and commutator estimates. 
\end{abstract}


\section{Introduction}

Excitable media play a central role in our understanding of complex systems. Chemical reactions \cite{bz,coox}, calcium waves \cite{calcium}, and neural field models \cite{bressloff:14,bressloff:12} are among the examples that motivate our present study. A prototypical model of excitable kinetics are the FitzHugh-Nagumo kinetics, derived first as a simplification of the Hodgkin-Huxley model for the propagation of electric signals through nerve fibers \cite{hodgkin-huxley:52},
\begin{subequations}
\label{eq:fhnk}
\begin{align}
\frac{\mathrm{d} u}{\mathrm{d} t} &=f(u)-v,\label{eq:fhnk1} \\
\frac{\mathrm{d} v}{\mathrm{d} t} &=\epsilon(u-\gamma v)\label{eq:fhnk2},
\end{align}
\end{subequations}
where, for instance, $f(u)=u(1-u)(u-a)$. For $0<a<1/2$ and $\gamma >0$, not too large,  all trajectories in this system converge to the trivial equilibrium $u=v=0$. The system is however excitable in the sense that finite-size perturbations of $u$, past the excitability threshold $a$, away from the stable equilibrium $u=v=0$, can induce a long transient, where $f(u)\sim v$, $u>1/2$.  During these transients, which last for times $\cO(1/\epsilon)$, $u$ is said to be in the excited state; eventually, $u$ returns to values $f(u)\sim v$, $u<1/2$, the quiescent state. 

Interest in these systems stems from the fact that, although kinetics are very simple and ubiquitous in nature, with convergence of all trajectories to a simple stable equilibrium,  spatial coupling can induce quite complex dynamics. The simplest example is the propagation of a stable excitation pulse, more complicated examples include two-dimensional spiral waves and spatio-temporal chaos. Intuitively, a local excitation  can trigger excitations of neighbors before decaying back to the quiescent state in a spatially coupled system. After initial transients, one then observes a spatially propagating region where $u$ belongs to the excited state. 

Rigorous approaches to the existence of such excitation pulses have been based on singular perturbation methods.  Consider, for example, 
\begin{subequations}
\label{eq:fhn}
\begin{align}
\partial_t \bu(x,t)&=\partial_{xx}\bu(x,t)+f(\bu(x,t))-\bv(x,t),\label{eq:fhn1} \\
\partial_t \bv(x,t)&=\epsilon(\bu(x,t)-\gamma \bv(x,t)).\label{eq:fhn2}
\end{align}
\end{subequations}
with $x\in\R$. One looks for solutions of the form 
\begin{equation}\label{eq:tw}
(\bu,\bv)(x,t)=(\bu,\bv)(x-ct),
\end{equation}
and finds first-order ordinary differential equation for $\bu,\bu_x,\bv$, in which one looks for a homoclinic solution to the origin. The small parameter $\epsilon$ introduces a singularly perturbed structure into the problem which allows one to find such a homoclinic orbit by tracking stable and unstable manifolds along fast intersections and slow, normally hyberbolic manifolds \cite{fhnex1,fhnex2,jones-etal:96}. This approach has been successfully applied in many other contexts with slow-fast like structures, with higher- or even infinite-dimensional slow-fast ODEs; see for instance \cite{doelmankaper,hupkes-sandstede:10,hupkes-sandstede:13}.

Our interest is in media with \emph{infinite-range} coupling. We will focus on \emph{linear} coupling through convolutions, although we believe that the existence result extends to a variety of other problems.  To fix ideas, we consider 
\begin{subequations}
\label{eq:nl}
\begin{align}
\partial_t \bu(x,t)&=-\bu(x,t)+\int_{\R}\K(x-y)\bu(y,t)dy+f(\bu(x,t))-\bv(x,t),\label{eq:nl1} \\
\partial_t \bv(x,t)&=\epsilon(\bu(x,t)-\gamma \bv(x,t)).\label{eq:nl2}
\end{align}
\end{subequations}
Our assumptions on the non-local coupling term $-\bu+\K \ast \bu$ in \eqref{eq:nl1} roughly require exponential localization and exponential stability of the excited and quiescent branch; see below for details. Our assumptions on $f$ and $\gamma$ encode excitability. In addition, we only require the existence of non-degenerate back and front solutions for the $u$-equation with frozen $\bv\equiv const$. Existence of such \emph{scalar} front solutions has been shown in many circumstances, for instance when $\K$ is positive. Non-degeneracy requires that the zero-eigenvalue of front and back, induced by translation, is algebraically simple. Again, such degeneracy is a consequence of monotonicity properties in many particular cases. Our main result states the existence of a traveling-wave solution (\ref{eq:tw}) for equations (\ref{eq:nl}). 

Traveling pulse solutions are stationary profiles $(\bu(\xi),\bv(\xi))$ of \eqref{eq:nl} in a comoving frame $\xi=x-ct$ that are localized so that $(\bu(\xi),\bv(\xi))\rightarrow 0$ as $\xi\rightarrow\pm\infty$. They satisfy the equations
\begin{subequations}
\label{eq:TP}
\begin{align}
-c \frac{d}{d\xi}\bu(\xi)&=-\bu(\xi)+\K\ast\bu\left(\xi\right)+f(\bu(\xi))- \bv(\xi), \\
-c \frac{d}{d\xi}\bv(\xi)&=\epsilon(\bu(\xi)-\gamma \bv(\xi)),
\end{align}
\end{subequations}
for some positive wave speed $c>0$. Due to the convolution term $\K\ast \bu$, the derivative of the state variables $\bu,\bv$ at a point $\xi$ in \eqref{eq:TP} depends on both advanced and retarded terms. Such systems are usually referred to as functional differential equations of mixed type. Considered as evolution equations in the time-like variable $\xi$, such equations present two major challenges:
\begin{itemize}
\item[(i)]the initial-value problem is ill-posed due to the presence of both advanced and retarded terms;
\item[(ii)] even for functional differential equations with only retarded terms, the infinite time horizon caused by the infinite range of the convolution kernel introduces technical difficulties.
\end{itemize}
The first difficulty has been overcome in various contexts, using exponential dichotomies as a major technical tool, instead of more geometric methods such as graph transforms; \cite{pss,hss,mvl}. In particular, existence and stability of both fronts and pulses have been established for such forward-backward systems with finite-range coupling; see for instance  \cite{mallet-paret:99a,mallet-paret:99b,hupkes-sandstede:09,hupkes-sandstede:10}. The second difficulty has not been addressed in the context of mixed-type equations. While for one-sided, retarded, say, coupling, several approaches are known that guarantee local well-posedness on suitable function spaces \cite{infinitedelay,infinitedelay2}, it is not clear how the constructions in \cite{pss,hss,mvl} would extend. 

Our approach avoids such complications, relying on more direct functional analytic tools instead of dynamical systems methods. We will give a precise statement of our result in the next section and conclude this introduction with a comparison of our results with results elsewhere in the literature.

Our result was primarily motivated by neural field equations. In fact, the existence problem for pulses in nonlocal excitable media was first addressed in the context of neural field equations with linear adaptation 
\cite{pinto-ermentrout:01,bressloff:12,bressloff:14}. Neural field equations  are nonlocal integro-differential equations of the form
\begin{subequations}
\label{eq:nfe}
\begin{align}
\partial_t \bu(x,t)&=-\bu(x,t)+\int_{\R}\K(x-y)S(\bu(y,t))dy-\bv(x,t),\label{eq:nfe1} \\
\partial_t \bv(x,t)&=\epsilon(\bu(x,t)-\gamma \bv(x,t)),\label{eq:nfe2}
\end{align}
\end{subequations}
where $\bu(x, t)$ represents the local activity of a population of neurons at position $x\in\R$ in the cortex, and the neural field $\bv(x,t)$ represents a form of negative feedback mechanism. The nonlinearity $S$ is the firing rate function and is often assumed to be of sigmoidal shape. Note that the main difference between systems \eqref{eq:nl} and \eqref{eq:nfe} is whether the nonlinearity acts inside or outside the convolution, a difference that does not affect the techniques we employ here. We note that in this context, kernels $\K$ are usually assumed to be positive, symmetric, and localized  \cite{ermentrout-mcleod:93,bates-etal:97,pinto-ermentrout:01}, matching the constraints that we will impose below. 

We conclude this introduction by mentioning two results on existence of pulses in nonlocal excitable media in the literature. Pinto \& Ermentrout \cite{pinto-ermentrout:01} use a formal singular singular limit to construct a leading order traveling pulse solution. They noticed that in a suitable spatial scaling, the convolutions converge to point evaluations, which allow one to construct a leading-order approximation of the profile in excited and recovery phases. The authors do not attempt to estimate or control errors of this leading-order approximation. Our paper can be viewed as doing just that, introducing a number of technical tools on the way. On the other hand, Faye \cite{faye:13}  exploited a special form of the kernel $\K$, which allows one to reduce the nonlocal problem to an equivalent local differential system. One can then rely again on geometric singular perturbation theory. The approach is, however, intrinsically limited to special, ``exponential type'' kernels that can be interpreted has 
Green's functions to linear differential equations.

\paragraph{Outline.} The remainder of this paper is organized as follows. We give a precise statement of our assumptions and state our main Theorem \ref{thm:existence} in Section \ref{main}. We also give a short sketch of proof, in particular relating techniques used here to the geometric methods used elsewhere. In Section \ref{slow}, we construct quiescent and excited pieces of the excitation pulse. We use those together with fast front and back solutions from the scalar problem in a leading-order Ansatz in Section \ref{ansatz}. Section \ref{proof} then puts all pieces together and concludes the proof of our main Theorem \ref{thm:existence}.

\section{Existence of excitation pulses --- main result}\label{main}
We formulate our main hypotheses, Section \ref{s:21}, state our main result, Section \ref{s:22}, and give an outline of the proof, Section \ref{subsec:sketch}. 
\subsection{Notation and hypotheses}\label{s:21}
We are interested in the existence of solutions $(\bu(\xi),\bv(\xi))$ of the system
\begin{subequations}
\label{eq:TPbis}
\begin{align}
-c \frac{d}{d\xi}\bu(\xi)&=-\bu(\xi)+\K\ast\bu\left(\xi\right)+f(\bu(\xi))- \bv(\xi), \\
-c \frac{d}{d\xi}\bv(\xi)&=\epsilon(\bu(\xi)-\gamma \bv(\xi)),
\end{align}
\end{subequations}
which are spatially localized, 
\bqs
\underset{\xi\rightarrow \pm\infty}{\lim}(\bu(\xi),\bv(\xi))=(0,0).
\eqs
Here, $c>0$ is the wave speed that needs to be determined as part of the problem and $0<\epsilon\ll 1$ is a small but fixed parameter.

\begin{figure}[t]
\centering
\includegraphics[width=0.48\textwidth]{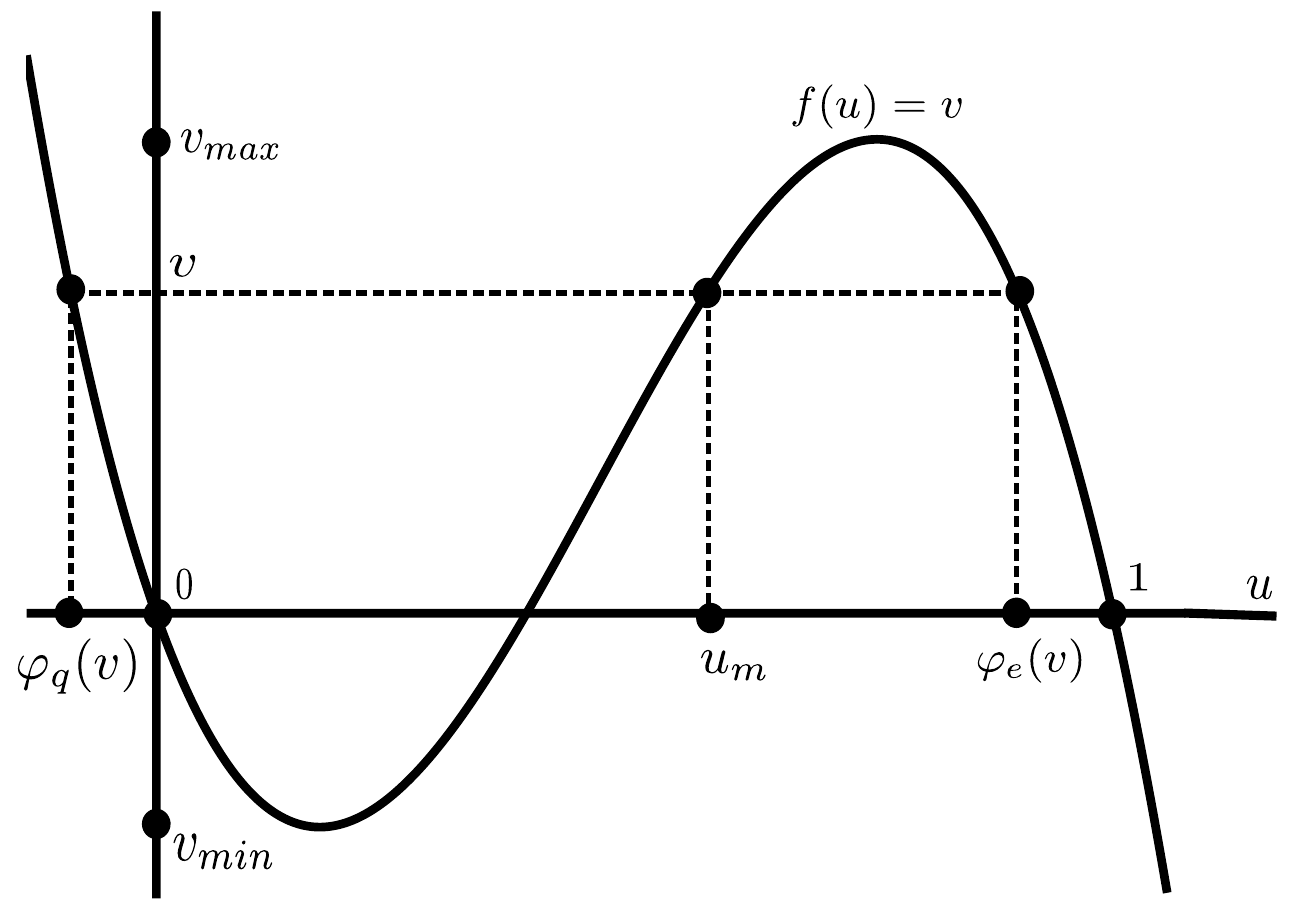}\hfill
\includegraphics[width=0.48\textwidth]{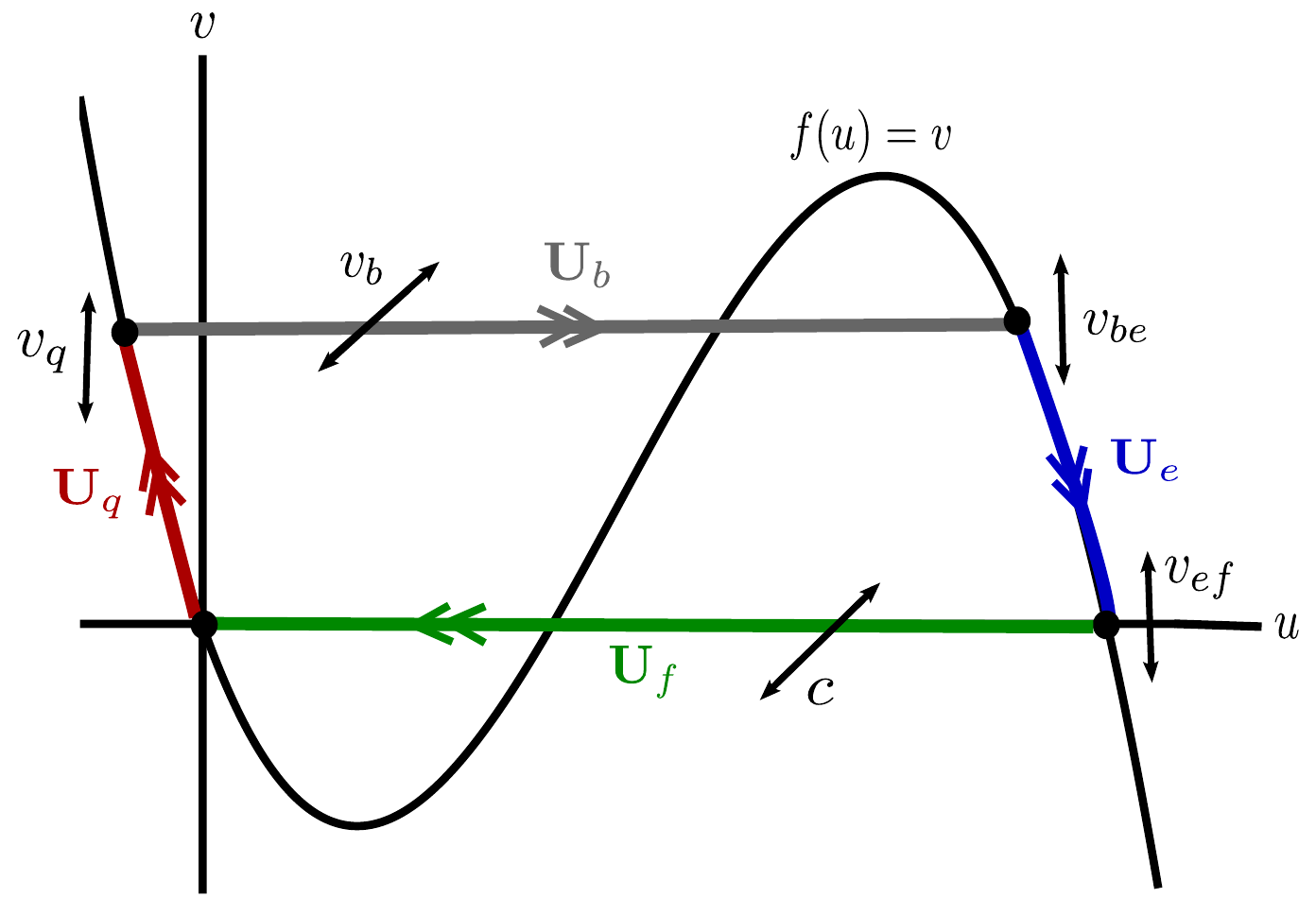}
\caption{Illustration of the assumptions on the nonlinearity $f$, left. To the right, the singular pulse, consisting of the quiescent part $\bU_q$ on the left branch of the slow manifold, the back $\bU_b$ connecting to the excited branch, the excitatory part $\bU_e$, and the front solution $\bU_f$. The five parameters $(v_q,v_b,v_{be},v_{ef},c)$  encode take-off and touch-down points, and ensure invertibility of the linearization at the singular solution.}\label{fig:cubicboth}
\end{figure}

Our first assumption concerns the nonlinearity, which we assume to be of excitable type. 
\begin{Hypothesis}{(H1)}
The nonlinearity $f$ is a $\cC^\infty$-smooth function with $f(0)=f(1)=0$, $f'(0)<0$ and $f'(1)<0$. Moreover, we assume that $\gamma>0$ is small enough so that $f(\gamma v)\neq v$. Lastly, we assume that $f(u)-v$ is of bistable type for $v\in(v_{min},v_{max})$, fixed, that is, it possesses precisely three nondegenerate zeroes. 
\end{Hypothesis}
The assumptions on $f$ are illustrated in Figure \ref{fig:cubicboth}. We denote the left and right zeroes of $f(u)-v$ by $u_q=\varphi_q(v)$ and $u_e=\varphi_e(v)$ and denote by $I_q$ and $I_e$ the ranges of $\varphi_{q}$ and $\varphi_e$. 

Our second assumption concerns the convolution kernel $\K$. For any $\eta \in \R$, we define  the space of exponentially  weighted functions on the real line equipped with its usual norm \[
L^1_\eta:=\left\{\bu:\R\rightarrow\R~|~ \int_\R e^{\eta|\xi|}|\bu(\xi)|d\xi<\infty \right\}.\]   We also write $\delta(\xi)$ for the dirac distribution with $\int\delta=1$. 
\begin{Hypothesis}{(H2)}
We suppose that the kernel $\K$ can be written as a sum $\K_{cont}+\K_{disc}$ with the following properties:
\begin{itemize}
\item There exists $\eta_0>0$ such that $\K_{cont} \in L^1_{\eta_0}$;
\item $\K_{disc}=\sum_{j\in\Z} a_j\delta(\xi-\xi_j)$, and $\sum_j|a_j|e^{\eta_0|\xi_j|}<\infty$;
\item the Fourier transform $\widehat\K(i\ell)$ of $\K$ satisfies $\widehat{\K}(0)=1$ and $\widehat{\K}(i\ell)-1<0$ for $\ell\neq 0$. 
\end{itemize}
\end{Hypothesis}
The first two assumptions, on regularity and on localization, mimic the assumptions in \cite{faye-scheel:13}, where Fredholm properties of nonlocal operators were established. The assumption $\widehat{\K}(0)=\int K=1$ is merely a normalization condition and can be achieved by scaling and redefining $f$. The last assumption can be slightly relaxed to 
\[
\widehat{\K}(i\ell)-1-f'(u)<0, \mbox{ for all } \ell\neq 0, \ u\in[\varphi_q(v_*),0]\cup [\varphi_e(v_*),1],
\]
where $v_*$ is defined in Hypothesis (H3), below.  Our assumptions do cover typical exponential or Gaussian kernels, as well as infinite-range pointwise interactions. A few comments on the last assumption are in order. Exponential localization guarantees that 
\bqs
\widehat{\K}(z)=\int_\R \K(x)e^{-zx}dx
\eqs
is analytic in a strip $|\Re(z)|< \eta_*$. Values of the characteristic function 
\bqs
\Delta_{j,v,c}(z)=zc-1+\widehat{\K}(z)+f'\left(u\right),
\eqs
determine the spectrum of the linearization at a constant state $u$. Our assumption then guarantees that constant states with $f'(u)<0$ do not possess zero spectrum,  also in spaces with exponential weights $|\eta|<\eta_*$ sufficiently small.


The last  assumption refers to the $u$-system with $\bv\equiv const$. Consider therefore 
\bqq
\label{eq:FBNL}
-c_* \frac{d}{d\xi}\bu(\xi)=-\bu(\xi)+\K\ast\bu\left(\xi\right)+f(\bu(\xi))-v_0,
\eqq
and the corresponding linearized operator
\bqq
\label{eq:FBNLL}
\cL(\bu_*)\bu(\xi)=c_* \frac{d}{d\xi}\bu(\xi)-\bu(\xi)+\K\ast\bu\left(\xi\right)+f'(\bu_*(\xi))\bu(\xi).
\eqq

\begin{Hypothesis}{(H3)}
We assume that there exists non-degenerate front and back solutions with equal speed. More precisely, there exists $c_*,v_*>0$ such that \eqref{eq:FBNL} possesses a \textbf{front} solution $\bu_f$ and a \textbf{back} solution $\bu_b$ with equal speed $c=c_*$, and $v$-values $v=0$  and $0<v=v_*<v_{max}$, respectively, that satisfy the limits
\begin{align*}
\underset{\xi\rightarrow-\infty}{\lim}\bu_f(\xi)=1,&\quad  \underset{\xi\rightarrow+\infty}{\lim}\bu_f(\xi)=0,\\
\underset{\xi\rightarrow-\infty}{\lim}\bu_b(\xi)=\varphi_q(v_*),&\quad  \underset{\xi\rightarrow+\infty}{\lim}\bu_b(\xi)=\varphi_e(v_*).
\end{align*}
Moreover, the operators $\cL(\bu_f)$ and $\cL(\bu_b)$ each possess an algebraically simple eigenvalue $\lambda=0$. 
\end{Hypothesis}
We remark that both linearized operators are automatically Fredholm of index zero \cite{faye-scheel:13}, so that the algebraic multiplicity of the eigenvalue $\lambda=0$ is finite. Since the derivatives of front and back profile contribute to the kernel, multiplicity is at least one. 

While hypotheses (H1) and (H2) are direct assumptions on nonlinearity and kernel, (H3) is an indirect assumption on both. For positive and even kernels, existence and stability can be established using comparison principles and monotonicity arguments; see for instance \cite{bates-etal:97,chen:97,bates-chen:06} for the specific case where $f(u)=u(1-u)(u-a)$, with $0<a<\frac{1}{2}$.  We also mention the early work of Ermentrout \& McLeod \cite{ermentrout-mcleod:93} who proved the existence of traveling front solutions for the neural field system \eqref{eq:nfe} with no adaptation. In a slightly different direction, De Masi \textit{et al.} proved existence and stability results for traveling fronts in nonlocal equations arising in Ising systems with Glauber dynamics and Kac potentials \cite{demasi-etal:95}. In all these cases, fronts are in fact monotone, a property that is however not needed in our construction. 

On the other hand, the set of hypotheses (H1)-(H3) form open conditions on nonlinearity and kernel: non-degenerate fronts can readily seen to persist under small perturbations, using for instance a variation of the methods presented in our proof.

\subsection{Main result -- summary}\label{s:22}

We can now state our main result.

\begin{thm}\label{thm:existence}
Consider the nonlocal FitzHugh-Nagumo equation \eqref{eq:nl} and suppose that Hypotheses (H1)-(H3) are satisfied; then for every sufficiently small $\epsilon>0$, there exist functions $\bu_\epsilon,\bv_\epsilon\in \mathcal{C}^1(\R,\R)$ and a wave speed $c(\epsilon)>0$ that depends smoothly on $\epsilon>0$ with $c(0)=c_*$, such that
\bqq
\label{eq:nlpulse}
(\bu(x,t),\bv(x,t))=(\bu_\epsilon(x-c(\epsilon)t),\bv_\epsilon(x-c(\epsilon)t))
\eqq
is a traveling wave solution of \eqref{eq:nl} that satisfies the limits
\bqq
\label{eq:limitpulse}
 \underset{\xi\rightarrow\pm\infty}{\lim}\left(\bu_\epsilon(\xi),\bv_\epsilon(\xi)\right)=(0,0).
\eqq
\end{thm}

Together with the discussion after Hypothesis (H3), we can state the following somewhat more explicit result.
\begin{cor}
The nonlocal FitzHugh-Nagumo equation, $f(u)=u(1-u)(u-a)$, $0<a<\frac{1}{2}$, $\gamma$ sufficiently small, $\K,\K'\in L^1_{\eta_0}$, $\K$ even, positive, with $\int\K=1$,  possesses a traveling pulse solution.
\end{cor}
Our approach is self-contained, roughly replacing subtle results on exponential dichotomies \cite{hss,mvl,pss} with crude Fredholm theory. 
Given the basic simplicity, we believe that our approach should cover a variety of different solution types and different media. For instance, one can readily see how to prove the existence of periodic wave trains in excitable or oscillatory regimes, or front solutions in bistable regimes. In analogy to the case of discrete media \cite{hupkes-sandstede:09}, we expect different phenomena when $c_*=0$, that is, for $a\sim 1/2$ in the cubic case, or for the slow pulse \cite{bates-non-smooth,krszsa}. Since the convolution operator does not regularize, compactly supported and discontinuous solutions can occur. 

\subsection{Sketch of the proof}\label{subsec:sketch}

Our proof of Theorem \ref{thm:existence} can be roughly divided into four main parts that can be outlined as follows. 

\paragraph{Step 1: Slow manifolds.}In a first step, we shall construct invariant slow manifolds for nonlocal differential equations of the form \eqref{eq:TPbis} for $0<\epsilon\ll1$ and $c>0$. Proving the persistence of invariant slow manifolds in the context of singularly perturbed ODEs was originally shown using graph transform \cite{fenichel:79}. Later, an alternative proof based on variation of constant formulas and exponential dichotomies  for differential equations with slowly varying coefficients was given \cite{sakamoto:90}. This latter approach was extended to ill-posed, forward-backward equations in  \cite{ssdefect,hupkes-sandstede:09}. Our approach completely renounces the concept of a phase space while picking up the main ingredients from the dynamical systems proofs: we modify nonlinearities outside a fixed neighborhood, construct an approximate trial solution, linearize at this ``almost solution'', and find a linear convolution type operator with slowly varying coefficients. We invert this operator by constructing suitable local approximate inverses and conclude the proof by setting up a Newton iteration scheme. We will see that the solution on the slow manifold satisfies a scalar ordinary differential equation, with leading order given by an expression equivalent to the one formally derived in \cite{pinto-ermentrout:01}.

\paragraph{Step 2: The singular solution.} We construct a singular solution using front and back solutions from Hypothesis (H3), together with pieces of slow manifolds from Step 1. We glue those solutions using appropriately positioned partitions of unity. Using partitions of unity instead of the matching procedure in cross-sections to the flow, common in dynamical systems approaches, is a second key difference of our approach. It allows us to avoid the notion of a phase space.  Schematically, the solution is formed by gluing together a quiescent part $\bU_q$ on the left branch of the slow manifold to a back solution $\bU_b$, then to an excitatory part $\bU_e$ on the right branch of the slow manifold, then to a front solution $\bU_f$ as shown in Figure \ref{fig:cubicboth}. On each solution piece, we allow for a correction $\bW$. See also Figure \ref{fig:ansatz} for a detailed picture.

\paragraph{Step 3: Linearizing and counting parameters.} In order to allow for weak interaction between the different corrections to solutions, we use function spaces with appropriately centered exponential weights. The weights, at the same time, encode the facts that solution pieces lie in either strong stable or unstable manifolds, or, in a more subtle way, the Exchange Lemma that tracks inclination of manifolds transverse to stable foliation forward with a flow \cite{jones-etal:96,brunovsky,jonesreview}. Our setup can be viewed as a version of \cite{brunovsky}, without phase space, in the simplest setting of a one-dimensional slow manifold. 

Linearizing at the different solution pieces, we find Fredholm operators with negative index. Roughly speaking, uniform exponential localization of perturbations does not allow corrections in the slow direction. In addition, the linearizations at back and front contribute one-dimensional cokernels, each. In the dynamical systems proofs, matching in cross-sections is accomplished by exploiting 
\begin{itemize}
\item free variables in stable and unstable manifolds;
\item auxiliary parameters, in our case $c$;
\item variations of touchdown and takeoff points on the slow manifolds.
\end{itemize}
We mimic precisely this idea, pairing the negative index Fredholm operators with suitable additional parameters, so that parameter derivatives span cokernels. A more detailed description is encoded in Figure \ref{fig:cubicboth}. We associate to the quiescent part $\bU_q$ the takeoff parameter $v_q\approx v_*$, which encodes the base point of the stable foliation that contains the back.  We associate 
to the excitatory part $\bU_e$ touchdown and takeoff parameters $v_{be}\approx v_*$ and $v_{ef}\approx 0$ that will compensate for the mismatched between the back and front parts. Finally we assign to the back $\bU_b$ the separate touchdown parameter $v_b\approx v_*$ and to the front $\bU_f$ the wave speed $c\approx c_*$. These two parameters effectively compensate for cokernels of front and back linearizations.

\paragraph{Step 4: Errors and fixed point argument.}

Our last step will be to use a fixed point argument to solve an equation of the form
\bqs
\F_\epsilon(\bW,\left(v_q,v_b,v_{be},v_{ef},c \right))=0,
\eqs
that is obtained by substituting our Ansatz directly into the system \eqref{eq:TPbis}. More precisely, we will show that 
\begin{itemize}
\item[(i)] $\| \F_\epsilon(\mathbf{0},\left(v_*,v_*,v_*,0,c_* \right))\| \rightarrow 0$ as $\epsilon \rightarrow 0$ in a suitable norm;
\item[(ii)] $D_{(\bW,\lambda)} \F_\epsilon(\mathbf{0},\left(v_*,v_*,v_*,0,c_* \right))$ is invertible with bounded inverse uniformly in $0<\epsilon \ll 1$;
\item[(iii)] $\F_\epsilon$ possesses a unique zero on suitable Banach spaces using a Newton iteration argument.
\end{itemize}
Here, (ii) follows from Step 3 and (iii) is a simple fixed point iteration. Errors (i) are controlled due to the careful choice of Ansatz and a sequence of commutator estimates between convolution kernels and linear or nonlinear operators.

\section{Persistence of slow manifolds}\label{slow}
In this section, we proof existence of solutions near the quiescent and the excited branch of $f(u)=v$,
\[
\M_q:=\left\{\left(\varphi_q(v),v \right) \right\},\qquad \M_e:=\left\{\left(\varphi_e(v),v \right)\right\}
.\]
We follow the ideas used in the construction of slow manifolds in dynamical systems and use a cut-off function to modify the slow flow outside a neighborhood that is relevant for our construction. We emphasize however that, due to the infinite-range coupling, the concept of solution defined locally in time is not applicable. In other words, the fact that we are modifying the equation outside of a neighborood will create error terms for all $\xi$. 

\begin{figure}[t]
\centering
\includegraphics[width=0.5\textwidth]{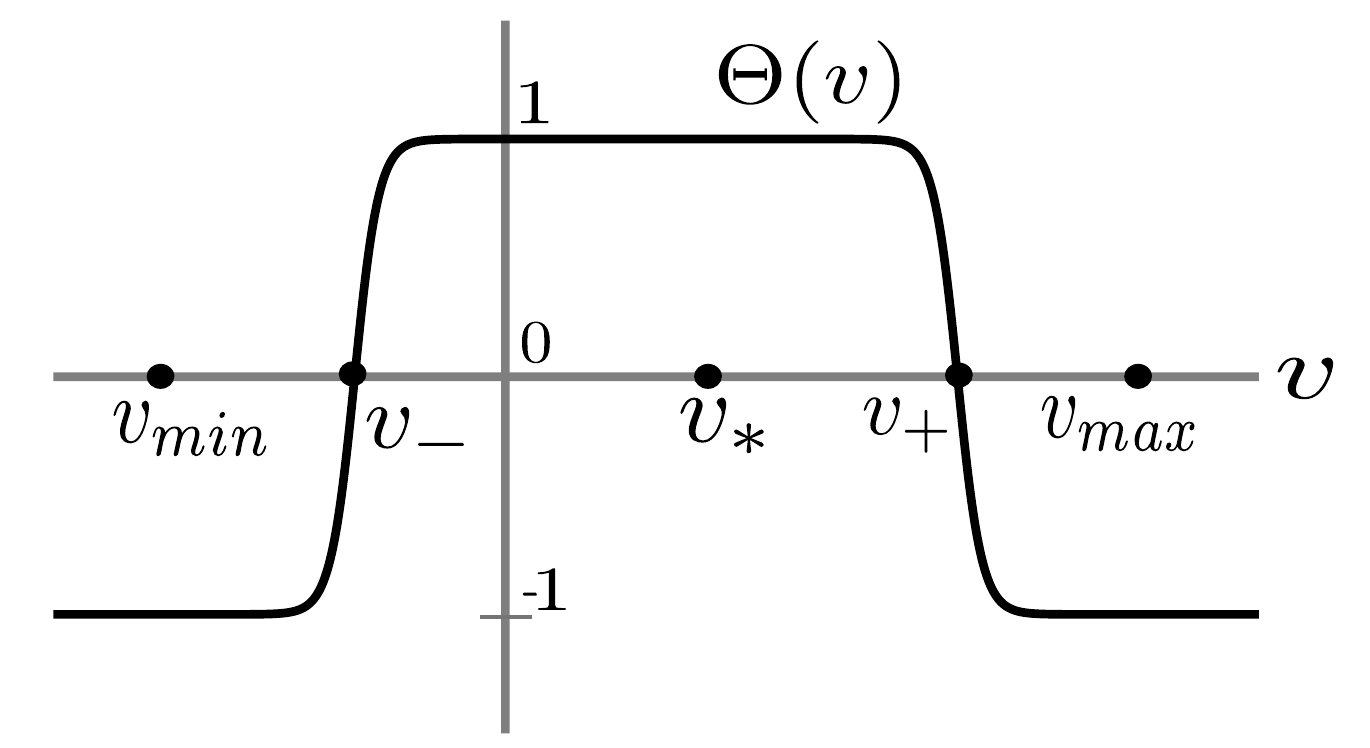}
\caption{The definition of the cut-off function $\Theta(v)$.}
\label{fig:cutoff}
\end{figure}

We use a simple modification of \eqref{eq:TP}, multiplying the right-hand side of the $v$-equation by a cut-off function $\Theta(v)$ as shown in Figure \ref{fig:cutoff}. The modified equation now reads
\begin{subequations}
\label{eq:TPmod}
\begin{align}
-c \frac{d}{d\xi}\bu(\xi)&=-\bu(\xi)+\K \ast \bu(\xi)+f(\bu(\xi))- \bv(\xi), \\
-c \frac{d}{d\xi}\bv(\xi)&=\epsilon(\bu(\xi)-\gamma\bv(\xi))\Theta(\bv(\xi)).
\end{align}
\end{subequations}
Formally, this introduces two equilibria on the slow manifold, with the effect that the solution on the slow manifold is expected to be a simple heteroclinic orbit. In order to exhibit the slow flow, we rescale space by introducing $\zeta = \epsilon \xi$ so that \eqref{eq:TPmod} becomes
\begin{subequations}
\label{eq:slow}
\begin{align}
-\epsilon c \frac{d}{d\zeta} \bu(\zeta) &=-\bu(\zeta)+\K_\epsilon \ast \bu(\zeta) +f(\bu(\zeta))- \bv(\zeta), \\
-c \frac{d}{d\zeta} \bv(\zeta) &=(\bu(\zeta)-\gamma \bv(\zeta))\Theta(\bv(\zeta)),
\end{align}
\end{subequations}
where we have defined the rescaled kernel as $\K_\epsilon (\zeta) := \epsilon^{-1}\K(\epsilon^{-1} \zeta)$. At $\epsilon=0$, the slow system is given by
\begin{subequations}
\label{eq:slow0}
\begin{align}
0 &=f(\bu(\zeta))- \bv(\zeta), \\
-c \frac{d}{d\zeta} \bv(\zeta) &=(\bu(\zeta)-\gamma \bv(\zeta))\Theta(\bv(\zeta)),
\end{align}
\end{subequations}

since formally, $\K_\epsilon \rightarrow\delta$, the Dirac distribution. Now, for each $c>0$, there exists a heteroclinic solution $(\varphi_q(\bv_{h,q}),\bv_{h,q})$ to \eqref{eq:slow0} on the quiescent slow manifold $\M_q$, connecting  the rest state $(0,0)$ to $(\varphi_q(v_+),v_+)$ for which the profile $\bv_{h,q}\in\cC^\infty(\R,\R)$ satisfies
\bqq
\label{eq:heteroclinicL}
-c \frac{d}{d\zeta} \bv =(\varphi_q(\bv)-\gamma\bv)\Theta(\bv), \quad (\varphi_q(\bv),\bv)\in\M_q
\eqq
with limits
\bqq
\label{eq:limitheteroclinicL}
\underset{\zeta\rightarrow-\infty}{\lim}\bv_{h,q}(\zeta)=0\text{ and } \underset{\zeta\rightarrow+\infty}{\lim}\bv_{h,q}(\zeta)=v_+.
\eqq
We normalize the solution so that $\bv_{h,q}(0)=v_*$.  Furthermore, for each $c>0$, there also exists a heteroclinic solution $(\varphi_e(\bv_{h,e}),\bv_{h,e})\in\M_e$ connecting the rest state $(\varphi_e(v_+),v_+)$ to $(\varphi_e(v_-),v_-)$ on the excitatory slow manifold $\M_e$ for which the profile $\bv_{h,e}\in \cC^\infty(\R,\R)$ satisfies
\bqq
\label{eq:heteroclinicR}
-c \frac{d}{d\zeta} \bv =(\varphi_e(\bv)-\gamma\bv)\Theta(\bv), \quad (\varphi_e(\bv),\bv)\in\M_e
\eqq
with limits
\bqq
\label{eq:limitheteroclinicR}
\underset{\zeta\rightarrow-\infty}{\lim}\bv_{h,e}(\zeta)=v_+\text{ and } \underset{\zeta\rightarrow+\infty}{\lim}\bv_{h,e}(\zeta)=v_-.
\eqq
We normalize this solution so that $\bv_{h,e}(0)=0$. 

\begin{figure}[t]
\centering
\includegraphics[width=0.6\textwidth]{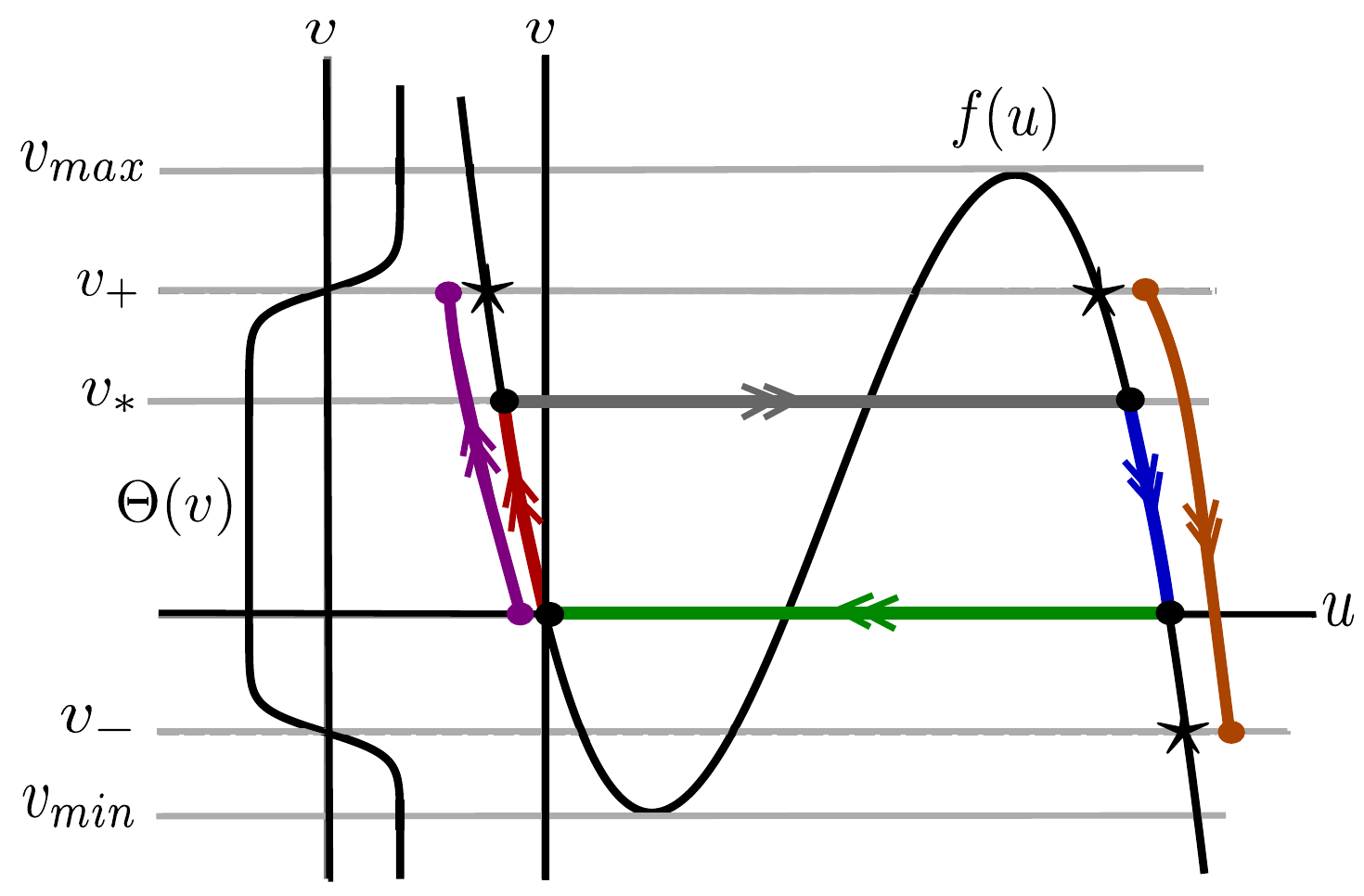}
\caption{Heteroclinics from Proposition \ref{prop:slowsolL} (left,purple) and from Proposition \ref{prop:slowsolR} (right,orange). The upper limit $v_+$ is induced by the cut-off $\Theta$ which is superimposed on the $v$-axis.}
\label{fig:CubicCutoff}
\end{figure}

Our goal in this section is to show that these two heteroclinic solutions persist for $0<\epsilon \ll 1$ using a fixed point argument  for nonlocal differential evolution equations with slowly varying coefficients. We give formal statements of the main result; a schematic picture of these heteroclinics relative to the singular pulse is shown in Figure \ref{fig:CubicCutoff}.

\begin{prop}[Quiescent slow manifold]\label{prop:slowsolL}
For every sufficiently small $\epsilon>0$ and any $c>0$, there exist functions $\bu_q,\bv_q\in\cC^\infty(\R,\R)$ such that
\bqq
\label{eq:slowsolL}
\left(\bu_q(\epsilon \xi),\bv_q(\epsilon \xi) \right)
\eqq
is a heteroclinic solution of \eqref{eq:TPmod} that satisfies the limits
\bqq
\label{eq:limitsolL}
\underset{\zeta\rightarrow-\infty}{\lim}\left(\bu_q(\zeta),\bv_q(\zeta)\right)=(0,0)\text{ and } \underset{\zeta\rightarrow+\infty}{\lim}\left(\bu_q(\zeta),\bv_q(\zeta)\right)=(\varphi_q(v_+),v_+).
\eqq
Up to translation, this solution is locally unique and depends smoothly on $\epsilon$ and $c$. 
\end{prop}

\begin{prop}[Excitatory slow manifold]\label{prop:slowsolR}
For every sufficiently small $\epsilon>0$ and  any $c>0$, there exist functions $\bu_e,\bv_e\in\cC^\infty(\R,\R)$ such that
\bqq
\label{eq:slowsolR}
\left(\bu_e(\epsilon \xi),\bv_e(\epsilon \xi) \right)
\eqq
is a heteroclinic solution of \eqref{eq:TPmod} that satisfies the limits
\bqq
\label{eq:limitsolR}
\underset{\zeta\rightarrow-\infty}{\lim}\left(\bu_e(\zeta),\bv_e(\zeta)\right)=(\varphi_e(u_+),u_+)\text{ and } \underset{\zeta\rightarrow+\infty}{\lim}\left(\bu_e(\zeta),\bv_e(\zeta)\right)=(\varphi_e(v_-),v_-).
\eqq
Up to translation, this solution is locally unique and depends smoothly on $\epsilon$ and $c$. 
\end{prop}

The proofs of these two propositions  will occupy the rest of this section. We remark that this construction of slow manifolds in nonlocal equations is somewhat general but comes with some caveats. First, the construction is simple here, since the slow manifold is one-dimensional and hence consists of a single trajectory, only. As a consequence, smoothness of slow manifolds is trivial, here. Second, the solutions are \emph{not} solutions for the original system, without the modifier $\Theta$, since the equation has infinite-range interaction in time. In other words, the modified piece of the trajectory influences the solution even where the solution takes values in the unmodified range. We will however exploit later that the error terms stemming from this modification are exponentially small due to the exponential localization of the kernel. 

We also note that monotonicity of $\bv_q$ (and similarly $\bv_e$) implies that $\bv_q$ solves a simple first-order differential equation, the ``reduced equation'' on the slow manifold. Again, this equation depends, even locally, on the modifier $\Theta$. From our construction, below, one can easily see that the leading-order vector field in $\epsilon$  is just the one given in \eqref{eq:heteroclinicR}.

\subsection{Set-up of the problem}

The strategy for the proof of  Propositions \ref{prop:slowsolL} and \ref{prop:slowsolR} is as follows. First, we introduce the map
\bqq
\label{eq:mapFeps}
\F^{\epsilon}:\begin{array}{ccl}
(\bu,\bv)&\longmapsto&  \left(c\epsilon \dfrac{d}{d\zeta}\bu-\bu+\K_\epsilon \ast \bu+f(\bu)- \bv, c \dfrac{d}{d\zeta}\bv+(\bu-\gamma \bv)\Theta(\bv)
 \right).
 \end{array}
\eqq
We can immediately confirm that any solution $(\bu,\bv)$ of $\F^{\epsilon}(\bu,\bv)=0$ is, by definition of the map $\F^{\epsilon}$, a solution of system \eqref{eq:slow}. From the above analysis, a natural extension to $\epsilon=0$ is
\bqq
\label{eq:mapFeps0}
\F^0:\begin{array}{ccl}
(\bu,\bv)&\longmapsto&  \left(f(\bu)- \bv, c \dfrac{d}{d\zeta}\bv+(\bu-\gamma \bv)\Theta(\bv)
 \right).
 \end{array}
\eqq
For sufficiently small $\epsilon>0$, $(\varphi_e(\bv_{h,e}),\bv_{h,e})$ should thus be an approximate solution to $\F^{\epsilon}(\bu,\bv)=0$, when $\bv_{h,e}$ is obtained from solving the second component of $\F^{\epsilon}$ with $\bu=\varphi_e(\bv_{h,e})$, \eqref{eq:heteroclinicR}. The following proposition quantifies the corresponding error.
\begin{prop}\label{prop:EstFeps} As $\epsilon \rightarrow 0$, the following estimate holds
\bqq
\label{eq:EstFeps}
\left\|\F^\epsilon~(\varphi_j(\bv_{h,j}),\bv_{h,j})\right\|_{L^2\times L^2}=\cO(\epsilon),
\eqq
for $j=q,e$.
\end{prop}
Suppose for a moment that we are able to prove the following result.
\begin{prop}\label{prop:Invertibility}
Let $D\F^{\epsilon}(\varphi_j(\bv_{h,j}),\bv_{h,j})$ be the linearization of $\F^{\epsilon}$ at the heteroclinic solution $(\bu,\bv)=(\varphi_j(\bv_{h,j}),\bv_{h,j})$, $j=q,e$, and  denote by $\X$ the Banach space $\X:=\left\{\bu\in H^1 ~|~ \bu(0)=0 \right\}$. Then, there exists $\epsilon_0$ and $C>0$ so that for all $0<\epsilon<\epsilon_0$ we have
\begin{itemize}
\item[(i)] $D\F^{\epsilon}(\varphi_j(\bv_{h,j}),\bv_{h,j}):H^1\times \X\to L^2 \times L^2$ is invertible;
\item[(ii)] $\left\| D\F^{\epsilon}(\varphi_j(\bv_{h,j}),\bv_{h,j})^{-1}\right\|\leq C $, uniformly in $\epsilon$;
\end{itemize}
for $j=q,e$.
\end{prop}
We can then set-up a Newton-type iteration scheme
\bqq
\label{eq:iteration}
(\bu_{n+1},\bv_{n+1})=\cS^\epsilon(\bu_n,\bv_n):=(\bu_n,\bv_n)-\left(D\F^{\epsilon}(\varphi_j(\bv_{h,j}),\bv_{h,j}) \right)^{-1}\F^{\epsilon}\left(\varphi(\bv_{h,j})+\bu_n,\bv_{h,j}+\bv_n \right),
\eqq
for $j=q,e$, with starting point $(\bu_0,\bv_0)=(\mathbf{0},\mathbf{0})$.  With the previous observations, we find that the map $\cS^\epsilon: H^1\times \X \longrightarrow H^1\times \X$  possesses the following properties:
\begin{itemize}
\item there exists $C_0>0$, such that $\left\|\cS^\epsilon(0,0)\right\|_{H^1\times \X}\leq C_0\epsilon$ as $\epsilon \rightarrow 0$;
\item $\cS^\epsilon$ is a $\cC^\infty$-map;
\item $D\cS^\epsilon(\mathbf{0},\mathbf{0})=0$;
\item there exist $\delta>0$ and $C_1>0$ such that for all $(\bu,\bv)\in\B_\delta$, the ball of radius $\delta$ centered at $(\mathbf{0},\mathbf{0})$ in $H^1\times \X$, $\left\| D\cS^\epsilon(\bu,\bv) \right\| \leq C_1 \delta$.
\end{itemize}

Now, suppose inductively that $(\bu_k,\bv_k)\in\B_\delta$ for all $1\leq k \leq n$, then
\bqs
\left\| (\bu_{n+1},\bv_{n+1})-(\bu_n,\bv_n) \right\|_{H^1 \times \X} \leq C_1 \delta \left\| (\bu_{n},\bv_{n})-(\bu_{n-1},\bv_{n-1}) \right\|_{H^1 \times \X},
\eqs
so that
\bqs
\left\| (\bu_{n+1},\bv_{n+1}) \right\|_{H^1 \times \X} \leq \frac{C_0}{1-C_1 \delta}\epsilon.
\eqs
For small enough $\epsilon$, we then have $\dfrac{C_0}{1-C_1 \delta}\epsilon<\delta$ and $(\bu_{n+1},\bv_{n+1})\in\B_\delta$, so that the map $\cS^\epsilon$ is a contraction.  Banach's fixed point theorem then gives a fixed point $(\bu^\epsilon,\bv^\epsilon)=\cS^\epsilon(\bu^\epsilon,\bv^\epsilon)$.

As a conclusion, for every sufficiently small $\epsilon>0$ and for each $c>0$, we have constructed functions $\bu_j$ and $\bv_j$ that can be written as
\begin{align*}
\bu_j&=\varphi(\bv_{h,j})+\bu_j^\epsilon\\
\bv_j&=\bv_{h,j}+\bv_j^\epsilon,
\end{align*}
such that $\left(\bu_j(\epsilon\xi), \bv_j(\epsilon \xi)\right)$ is a heteroclinic solution of \eqref{eq:TPmod} with $j=q,e$.

It now remains to prove Propositions \ref{prop:EstFeps} and \ref{prop:Invertibility}. In order to simplify our notation, we will  write 
\bqs
(\bu_h,\bv_h)=(\varphi(\bv_{h,j}),\bv_{h,j}),
\eqs
not distinguishing between the cases $j=q$ and $j=e$ since proofs in both cases are completely equivalent.

\subsection{Proof of Proposition \ref{prop:EstFeps}}
A direct computation shows that
\bqs
\F^\epsilon(\bu_h,\bv_h)(\zeta)=\left(\epsilon c \frac{d}{d\zeta} \bu_h(\zeta)
-\bu_h(\zeta)+\K_\epsilon \ast \bu_h(\zeta),0\right),
\eqs
for all $\zeta\in\R$. In the following, we use  $A\lesssim B$ whenever $A<C B$, with $C$ independent of $\epsilon$. The key ingredient to the proof is a comparison of the rescaled convolution with a dirac delta.
\begin{prop}\label{prop:EstIntF}
For any $\bw \in H^1$, $\| -\bw+\K_\epsilon \ast \bw  \|_{L^2}\lesssim \epsilon \|\bw\|_{H^1} $. 
\end{prop}
\begin{Proof}\
For all $\zeta\in\R$, we have
\bqs
-\bw(\zeta)+\K_\epsilon \ast \bw (\zeta) = \int_\R \K(y)\left( \bw(\zeta-\epsilon y)-\bw(\zeta)\right)dy=-\epsilon \int_\R y \K(y) \int_0^1 \bw'(\zeta-\epsilon y s) ds dy.
\eqs
Then
\bqs
\left( \int_\R y \K(y) \int_0^1 \bw'(\zeta-\epsilon y s) ds dy \right)^2\leq \left( \int_\R y^2 |\K(y)| dy\right)\left(\int_\R|\K(y)|\int_0^1 \bw'(\zeta-\epsilon y s)^2 ds dy \right).
\eqs
Here, we have used the fact that $\K$ is exponentially localized so that its second moment always exists. This readily yields
\begin{align*}
\|-\bw+ \K_\epsilon \ast \bw \|_{L^2}&\leq \epsilon \left( \int_\R y^2 |\K(y)| dy\right)^{\frac{1}{2}} \|\bw\|_{H^1}.
\end{align*}
\end{Proof}

Next, define the affine spaces
$\widetilde{L}^2_j=\Phi_j+L^2_j$, $j=q,e$, with distance given by the $L^2$ norm, where 
\bqs
\Phi_e(\zeta)=\varphi_e(u_+)\chi_-(\zeta)+\varphi_e(v_-)\chi_+(\zeta)\text{ and } 
\Phi_q(\zeta)=\varphi_q(v_+)\chi_+(\zeta),
\eqs
where  $\chi_\pm=(1\pm\tanh(\zeta))/2$. In an analogous fashion, we define $\widetilde{H}^1_j$.
Now the map 
\bqs
\bu \mapsto   \epsilon c \frac{d}{d\zeta} \bu-\bu+\K_\epsilon \ast \bu
\eqs
is well defined from $\widetilde{H}^1_j$ to $L^2$. Combining the fact that $\bu_h \in \widetilde{H}^1$ and the above proposition, we obtain
\bqs
\left\| -\bu_h+\K_\epsilon \ast \bu_h+\epsilon c \frac{d}{d\zeta} \bu_h\right\|_{L^2}\lesssim \epsilon.
\eqs
This implies that $\left\|\F^\epsilon(\bu_h,\bv_h)\right\|_{L^2\times L^2}=\cO(\epsilon)$ and thus completes the proof of Proposition  \ref{prop:EstFeps}.

\subsection{Proof of Proposition \ref{prop:Invertibility}}\label{subsec:DFinv}

We recall that we obtained $D\F^\epsilon(\bu_h,\bv_h)$ by linearizing equation \eqref{eq:slow} around the heteroclinic solution $(\bu_h,\bv_h)$ found for $\epsilon=0$. A convenient way to represent  $D\F^\epsilon(\bu_h,\bv_h)$ is through its matrix form
\bqq
\label{eq:DFeps}
D\F^\epsilon(\bu_h,\bv_h)=\left( 
\begin{matrix}
\cL_{\epsilon,\zeta} & \cL_{u,v} \\
\cL_{v,u} & \cL_{v,v}
\end{matrix}
\right),
\eqq
where the linear operators $\cL_{\epsilon,\zeta}$, $\cL_{u,v}$, $\cL_{v,u}$ and $\cL_{v,v}$ are defined as follows
\begin{subequations}
\begin{align}
\label{eq:Leps}
\cL_{\epsilon,\zeta}&: \left\{\begin{array}{ccl}
H^1&\longrightarrow& L^2\\
 \bu&\longmapsto& -\bu+\K_\epsilon \ast \bu + \epsilon c \frac{d}{d\zeta}\bu +f'\left(\bu_h(\zeta)\right)\bu,
 \end{array}\right.\\
\label{eq:Luv}
\cL_{u,v}&:  \left\{\begin{array}{ccl}
H^1&\longrightarrow& L^2\\
 \bu&\longmapsto& - \bu,
\end{array}\right.\\
\label{eq:Lvu}
\cL_{v,u}&:  \left\{\begin{array}{ccl}
H^1&\longrightarrow& L^2\\
 \bu&\longmapsto&  \Theta\left(\bv_h(\zeta)\right) \bu,
\end{array}\right.\\
\label{eq:Lvv}
\cL_{v,v}&:  \left\{\begin{array}{ccl}
H^1&\longrightarrow& L^2\\
 \bu&\longmapsto& c \frac{d}{d\zeta}\bu- \gamma \left(\Theta\left(\bv_h(\zeta)\right)\bu+ \bv_h(\zeta)\Theta'\left(\bv_h(\zeta)\right)\bu\right).
\end{array}\right.
\end{align}
\end{subequations}

To prove Proposition \ref{prop:Invertibility}, we will solve the linear system
\bqq
\label{eq:linsyst}
\left( 
\begin{matrix}
\cL_{\epsilon,\zeta} & \cL_{u,v} \\
\cL_{v,u} & \cL_{v,v}
\end{matrix}
\right)
\left(\begin{matrix} \bu\\ \bv \end{matrix} \right)=\left(\begin{matrix} \bh\\ \bg \end{matrix} \right),
\eqq
for all $\bh,\bg\in L^2$ and $\bu,\bv\in H^1$. Mimicking the dynamical systems approach of first diagonalizing the frozen system, at $\epsilon=0$, we change variables
\bqq
\label{eq:linsystmod}
\left( 
\begin{matrix}
\cL_{\epsilon,\zeta} & \widetilde{\cL}_{u,v} \\
\cL_{v,u} & \cL_{h,\zeta}
\end{matrix}
\right)
\left(\begin{matrix} \widetilde{\bu}\\ \bv \end{matrix} \right)=\left(\begin{matrix} \bh\\ \bg \end{matrix} \right),
\qquad \mbox{with }\widetilde{\bu}=\bu-\frac{1}{f'\left(\bu_h(\zeta)\right)}\bv,
\eqq
and
\begin{subequations}
\begin{align}
\label{eq:Luvt}
\widetilde{\cL}_{u,v}&:  \left\{\begin{array}{ccl}
H^1&\longrightarrow& L^2\\
 \bu&\longmapsto& -\dfrac{1}{f'\left(\bu_h(\zeta)\right)} \bu+\K_\epsilon \ast \left(\dfrac{1}{f'\left(\bu_h(\zeta)\right)} \bu \right)+ \epsilon c \dfrac{d}{d\zeta}\left( \dfrac{1}{f'\left(\bu_h(\zeta)\right)} \bu\right),
 \end{array}\right.\\
\label{eq:Lhzeta}
\cL_{h,\zeta}&:  \left\{\begin{array}{ccl}
H^1&\longrightarrow& L^2\\
 \bu&\longmapsto& \cL_{v,v}\bu+ \dfrac{\Theta\left(\bv_h(\zeta)\right)}{f'\left(\bu_h(\zeta)\right)} \bu.
\end{array}\right.
\end{align}
\end{subequations}

Using Proposition \ref{prop:EstIntF} and the fact that $\zeta \mapsto \left(f'\left(\bu_h(\zeta)\right)\right)^{-1}$ is a bounded function for all $\zeta\in\R$, we directly obtain that $\| \widetilde{\cL}_{u,v} \|_{H^1\rightarrow L^2}=\cO(\epsilon)$, so that it is sufficient to solve 
\bqq
\label{eq:triangular}
\left( 
\begin{matrix}
\cL_{\epsilon,\zeta} & 0 \\
\cL_{v,u} & \cL_{h,\zeta}
\end{matrix}
\right)
\left(\begin{matrix} \widetilde{\bu}\\ \bv \end{matrix} \right)=\left(\begin{matrix} \bh\\ \bg \end{matrix} \right),
\eqq
with $\epsilon$-uniform bounds. This in turn follows from  obvious $\epsilon$-uniform bounds on $\cL_{v,u}$ and $\epsilon$-uniform invertibility of  $\cL_{\epsilon,\zeta}: H^1\to L^2$ and $\cL_{h,\zeta}:\X\to L^2$.

\subsubsection{Invertibility of $\cL_{\epsilon,\zeta}$}

First, let $\bh\in L^2$ and consider the frozen system
\bqq
\label{eq:frozen}
\cL_{\epsilon,\zeta_0}\bu=\bh
\eqq
with $\zeta_0\in\R$ fixed. The solution is obtained by convolution with the Green's function  $\G_\epsilon(~\cdot~;\zeta_0):\R\rightarrow\C$ which we obtain as follows. We  define $\Delta_{\zeta_0}(i\ell)$ as
\bqq
\label{eq:Delta}
\Delta_{\zeta_0}(i\ell)=-1+\widehat{\K}(i\ell) + f'\left(\bu_h(\zeta_0)\right)+i \ell c, \quad \eta\in\R
\eqq
where we have set 
\bqs
\widehat{\K}(i\ell)=\int_\R \K(\zeta)e^{-i\ell\zeta}d\zeta. 
\eqs

We observe that $\Delta_{\zeta_0}(i\ell)=\cO(|\ell|)$ as $\ell\rightarrow \pm \infty$ and, since $f'\left(u_h(\zeta_0)\right)<0$, $\Delta_{\zeta_0}(i\ell)^{-1}$ is well defined for all $\ell\in\R$, so the function $\ell\mapsto\Delta_{\zeta_0}(i\ell)^{-1}$ belongs to $L^2$.  We may therefore construct its inverse Fourier transform,
\bqs
\G(\zeta;\zeta_0):=\frac{1}{2\pi}\int_\R e^{i\zeta \ell} \Delta_{\zeta_0}(i\ell)^{-1} d\ell\in L^2.
\eqs
Lastly, the Green's function $\G_\epsilon$ is now given through
\bqq
\label{eq:Green}
\G_\epsilon(\zeta;\zeta_0)=\epsilon^{-1}\G\left(\epsilon^{-1}\zeta;\zeta_0\right), \quad \forall \zeta\in\R.
\eqq

\begin{prop}
The operator $\cL_{\epsilon,\zeta_0}:H^1\to L^2$ with $\epsilon$ small, $\zeta_0$ fixed, is an isomorphism, with inverse given by convolution with $\G_\epsilon(\zeta;\zeta_0)$ from \eqref{eq:Green},
\bqq
\label{eq:solfrozen}
\left(\cL_{\epsilon,\zeta_0}^{-1}\bh\right)(\zeta)= \left(\G_\epsilon \ast \bh\right)(\zeta) = \int_\R \G_\epsilon(\zeta-\tilde \zeta;\zeta_0)\bh(\tilde\zeta)d\tilde\zeta.
\eqq
\end{prop}

\begin{Proof}

Interpreting $\G_\epsilon$ as a tempered distribution, we consider the distribution
\bqs
F(\zeta)=\cL_{\epsilon,\zeta_0}\G_\epsilon(\zeta;\zeta_0).
\eqs
We can evaluate the Fourier transform $\widehat F$ of $F$ and find
\bqs
\widehat F(\ell) = \left(-1+\widehat{\K}(i\epsilon \ell) + f'\left(u_h(\zeta_0)\right)+i \epsilon \ell c \right) \widehat\G_\epsilon(\ell;\zeta_0) = \Delta_{\zeta_0}(i\epsilon\ell)\widehat\G(\epsilon \ell;\zeta_0)=\Delta_{\zeta_0}(i\epsilon\ell)\Delta_{\zeta_0}(i\epsilon\ell)^{-1}=1.
\eqs
Thus $F=\delta$ where $\delta$ denotes the Dirac delta distribution. Since $\Delta_{\zeta_0}$ is analytic in a strip, one can readily show that $\G$ and $\G_\epsilon$ are exponentially localized, hence belong to $L^1$. 

We can now define $\bu$ via  convolution $\bu=\G_\epsilon \ast \bh$ and Young's inequality gives 
\bqs
\|\bu\|_{L^2}\leq \| \G_\epsilon \|_{L^1(\R)} \|\bh\|_{L^2}.
\eqs
One readily verifies that  $\bu$ satisfies  \eqref{eq:frozen} in the sense of distributions and we conclude that $\bu\in H^1$. It then follows that $\cL_{\epsilon,\zeta_0}:H^1\longrightarrow L^2$ is onto. It remains to show that $\cL_{\epsilon,\zeta_0}$ is one-to-one. Suppose that $\cL_{\epsilon,\zeta_0}\cdot \bu=0$ for $\bu\in H^1$. Then using Fourier transform on both sides we obtain
\bqs
\left(-1+\widehat{\K}(i\epsilon \ell) + f'\left(\bu_h(\zeta_0)\right)+i \epsilon \ell c \right)\widehat \bu(\ell)=0
\eqs
for all $\ell\in\R$. Then $\widehat \bu(\ell)=0$ and hence $\bu$ is the zero function. 
\end{Proof}

We now return to the construction of a right inverse of $\cL_{\epsilon,\zeta}$.  Let $\bh\in L^2$ and consider the unfrozen system
\bqq
\label{eq:unfrozen}
\cL_{\epsilon,\zeta}\bu=-\bu+\K_\epsilon \ast \bu + \epsilon c \frac{d}{d\zeta}\bu +f'\left(\bu_h(\zeta)\right)=\bh, \quad \forall \zeta\in\R.
\eqq 
Exploiting the fact that coefficients are varying slowly, we use the solution of the frozen system \eqref{eq:solfrozen} as an Ansatz for  \eqref{eq:unfrozen} and show smallness of remainder terms. Therefore, define
\bqq
\label{eq:N}
\widetilde{\bu}(\zeta)= \cN \bh(\zeta):= \int_\R \G_\epsilon(\zeta-\tilde\zeta;\zeta)\bh(\tilde\zeta)d\tilde\zeta
\eqq
for all $\zeta\in\R$. 
\begin{lem}\label{lem:nbd}
The operator $\cN:L^2\to H^1$ is bounded, uniformly in $\epsilon$. 
\end{lem}
\begin{Proof}
From its definition \eqref{eq:N}, we obtain, using Holder's inequality, that
\begin{align*}
 \int_\R \left(\cN \bh(\zeta)\right)^2d\zeta&=\int_\R \left(\int_\R \G_\epsilon(\zeta-\tilde\zeta;\zeta)\bh(\tilde\zeta)d\tilde\zeta\right)^2d\zeta\\
 &\leq \int_\R \left(\int_\R \underset{y\in\R}{\sup~} \left|\G_\epsilon(\zeta-\tilde\zeta;y)\right|\bh(\tilde\zeta)d\tilde\zeta\right)^2d\zeta\\
 &\leq \left(\int_\R \underset{y\in\R}{\sup~} \left|\G_\epsilon(\zeta;y)\right|d\zeta\right)^2\int_\R \bh(\zeta)^2d\zeta.
\end{align*}
The claim now follows from 
\bqs
\int_\R \underset{y\in\R}{\sup~} \left|\G_\epsilon(\zeta;y)\right|d\zeta=\int_\R \underset{y\in\R}{\sup~} \left|\G(\zeta;y)\right|d\zeta=M.
\eqs
\end{Proof}
We next substitute the Ansatz $\bu=\widetilde{\bu}+\bu_1$ into \eqref{eq:unfrozen}, with $\bu_1\in H^1$. We find that $\bu_1$ satisfies 
\bqq
\label{eq:u1}
\cL_{\epsilon,\zeta}\bu_1=\left(\I-\cL_{\epsilon,\zeta}\cN \right) \bh(\zeta)=:\cR_{\epsilon,\zeta}\bh(\zeta).
\eqq

\begin{prop}
For all $\bh\in L^2$, we have
\bqs
\|\cR_{\epsilon,\zeta}\bh\|_{L^2}\lesssim \epsilon \|\bh\|_{L^2}.
\eqs
\end{prop}

\begin{Proof}
First, observe that the differentiability of $\zeta_0\mapsto \G_\epsilon(\zeta;\zeta_0)$ follows directly from  its definition \eqref{eq:Green} and the differentiability of $\zeta \mapsto f'\left(\bu_h(\zeta)\right)$. We then denote $\partial_2 \G_\epsilon(\cdot;\cdot)$ for the partial derivate with respect to the second component. Next, since $\Delta_{\zeta_0}(\eta)=\cO(|\eta|)$ as $\eta\rightarrow \pm \infty$, we have $\partial_{\zeta_0}\left(\Delta_{\zeta_0}(\eta)^{-1}\right)=\cO\left(|\eta|^{-2}\right)$ as $\eta\rightarrow \pm \infty$, and the function $\eta\mapsto \partial_{\zeta_0}\left(\Delta_{\zeta_0}(\eta)^{-1}\right)$ belongs to $L^1$ so that
\bqs
\frac{d}{d\zeta}\int_\R \G_\epsilon(\zeta-\tilde\zeta;\zeta)\bh(\tilde\zeta)d\tilde\zeta=\int_\R \partial_1 \G_\epsilon(\zeta-\tilde\zeta;\zeta)\bh(\tilde\zeta)d\tilde\zeta+\int_\R \partial_2 \G_\epsilon(\zeta-\tilde\zeta;\zeta)\bh(\tilde\zeta)d\tilde\zeta, \quad \forall\zeta\in\R,
\eqs
where $\partial_1 \G_\epsilon(\cdot;\cdot)$ stands for the partial derivate with respect to the first component.
A direct computation shows that
\begin{align*}
\cR_{\epsilon,\zeta}\bh(\zeta)&=\int_\R \K_\epsilon(\zeta-\tilde\zeta)\left(\int_\R\left[\G_\epsilon(\tilde\zeta-\breve\zeta;\zeta)-\G_\epsilon(\tilde\zeta-\breve\zeta;\tilde\zeta) \right]\bh(\breve\zeta)d\breve\zeta\right)d\tilde\zeta-\epsilon c \int_\R \partial_2 \G_\epsilon(\zeta-\tilde\zeta;\zeta)\bh(\tilde\zeta)d\tilde\zeta\\
&=\J_1(\zeta)+\J_2(\zeta).
\end{align*}
For all $(\zeta,\tilde\zeta,\breve\zeta)\in \R^3$, we have
\bqs
\G_\epsilon(\tilde\zeta-\breve\zeta;\zeta)-\G_\epsilon(\tilde\zeta-\breve\zeta;\tilde\zeta) =(\zeta-\tilde\zeta)\int_0^1\partial_2 \G_\epsilon(\tilde\zeta-\breve\zeta;(1-s)\tilde\zeta+s\zeta)ds,
\eqs
so that we can write $\J_1$ as
\begin{align*}
\J_1(\zeta)&= \epsilon \int_\R (\zeta-\tilde\zeta)\K_\epsilon(\zeta-\tilde\zeta)\int_\R\left(\bh(\breve\zeta)\int_0^1\partial_2 \G_\epsilon(\tilde\zeta-\breve\zeta;(1-s)\tilde\zeta+s\zeta)dsd\breve\zeta\right)d\tilde\zeta\\
&= \epsilon \int_\R y\K(y)\int_\R\left(\bh(\breve\zeta)\int_0^1\partial_2 \G_\epsilon(\zeta-\epsilon y-\breve\zeta;\zeta-(1-s)\epsilon y)dsd\breve\zeta\right)dy
\end{align*}
and
\bqs
\| \J_1 \|_{L^2}\leq \epsilon \left(\int_\R|\zeta||\K(\zeta)|d\zeta \right)^{1/2} \left(\int_\R \left| \underset{y\in\R}{\sup~}\partial_2 \G_\epsilon(\zeta;y)\right|d\zeta \right)\| \bh \|_{L^2}.
\eqs
Furthermore, we also have
\bqs
\| \J_2 \|_{L^2}\leq \epsilon c \left(\int_\R \left| \underset{y\in\R}{\sup~}\partial_2 \G_\epsilon(\zeta;y)\right|d\zeta \right)\| \bh \|_{L^2}.
\eqs
which completes the proof.
\end{Proof}

We can construct a solution of \eqref{eq:u1} of the form
\bqs
\bu_1=\widetilde{\bu}_1+\bu_2, \quad \widetilde{\bu}_1=\cN \cR_{\epsilon,\zeta}\bh, \quad \bu_2\in H^1,
\eqs
and $\bu_2$ is solution of 
\bqs
\cL_{\epsilon,\zeta} \bu_2=\left(\cR_{\epsilon,\zeta}- \cL_{\epsilon,\zeta}\cN \cR_{\epsilon,\zeta}\right)\bh=\left(\I-\cL_{\epsilon,\zeta}\cN \right)\cR_{\epsilon,\zeta}\bh=\cR_{\epsilon,\zeta}^2\bh
\eqs
with $\|\cR_{\epsilon,\zeta}^2\|\preceq \epsilon^2$. We inductively construct a sequence of functions 
\bqs
\bu_n=\widetilde{\bu}_n+\bu_{n+1}, \quad \widetilde{\bu}_n=\cN \cR_{\epsilon,\zeta}^nh, \quad \bu_{n+1}\in H^1,
\eqs
where $\bu_{n+1}$ solves
\bqs
\cL_{\epsilon,\zeta} \bu_{n+1}=\cR_{\epsilon,\zeta}^{n+1}\bh, \quad \|\cR_{\epsilon,\zeta}^{n+1}\|\preceq \epsilon^{n+1}.
\eqs

Then, for $\epsilon>0$ small enough, we obtain $\bu$, solution of \eqref{eq:unfrozen}, from the convergent geometric series
\bqq
\label{eq:solunfrozen}
\bu=\cN\left(\sum_{n=0}^{\infty}\cR_{\epsilon,\zeta}^n \right)\bh=\cN\left(\I-\cR_{\epsilon,\zeta} \right)^{-1}\bh, \quad \forall \bh\in L^2.
\eqq
As a consequence, $\cL_{\epsilon,\zeta}:H^1\longrightarrow L^2$ is onto. Next, let $\bu\in H^1$ be such that
\bqs
\cL_{\epsilon,\zeta}\bu=0.
\eqs
Multiplying both sides by $\bu$ and integrating over the real line, we find
\bqs
\int_\R \bu(\zeta)\cL_{\epsilon,\zeta}\bu(\zeta)d\zeta= \int_\R \bu(\zeta)\left(-\bu(\zeta)+\K_\epsilon \ast \bu(\zeta) \right)d\zeta +\int_{\R}f'\left(\bu_h(\zeta)\right)\bu^2(\zeta)d\zeta.
\eqs
Using Parseval's identity on the first term of the above equation, we obtain
\bqs
 \int_\R \bu(\zeta)\left(-\bu(\zeta)+\K_\epsilon \ast \bu(\zeta) \right)d\zeta= 2\pi\int_\R \left(-1+\widehat{\K}(i\epsilon \ell) \right)\widehat{\bu}(\ell)^2d\ell \leq 0.
\eqs
Since $f'\left(\bu_h(\zeta)\right)<0$ for all $\zeta\in\R$, we have $\int_\R \bu(\zeta)\cL_{\epsilon,\zeta}\bu(\zeta)d\zeta<0$ unless $\bu=0$, which proves that $\cL_{\epsilon,\zeta}$ is one-to-one. In conclusion we have proved the following result.

\begin{prop}\label{prop:ConstInverse}
The operator $\cL_{\epsilon,\zeta}:H^1\to L^2$ is an isomorphism  with $\epsilon$-uniformly bounded inverse 
\bqs
\cL_{\epsilon,\zeta}^{-1}=\cN\left(\I-\cR_{\epsilon,\zeta} \right)^{-1},
\eqs
where  $\cN$ and $\cR_{\epsilon,\zeta}$ were defined  in \eqref{eq:N} and \eqref{eq:u1}, respectively.
\end{prop}

\subsubsection{Invertibility of $\cL_{h,\zeta}$}

In this section we show that $\cL_{h,\zeta}$ is invertible from $\X=\left\{\bu\in H^1~|~ \bu(0)=0 \right\}$ to $L^2$. We define the operator $\T:H^1\rightarrow L^2$ as
\bqs
\T \bu=\dfrac{d}{d\zeta}\bu-\A(\zeta)\bu, \quad \A(\zeta)=c^{-1}\left( \gamma\left(\Theta\left(\bv_h(\zeta)\right)+ \bv_h(\zeta)\Theta'\left(\bv_h(\zeta)\right)\right)- \dfrac{\Theta\left(\bv_h(\zeta)\right)}{f'\left(\bu_h(\zeta)\right)} \right),
\eqs
with limiting entries $\A_\pm=\underset{\zeta\rightarrow\pm\infty}{\lim}\A(\zeta)$ given by
\bqs
\A_-=c^{-1}\left(\gamma-\dfrac{1}{f'\left(0\right)} \right)>0 \text{ and } \A_+=c^{-1}\gamma v_+\Theta'(v_+)<0,
\eqs
for the quiescent case, and by
\bqs
\A_-=c^{-1}\left(\gamma-\dfrac{1}{f'\left(\varphi_e(v_+)\right)} \right)>0 \text{ and } \A_+=c^{-1}\gamma v_-\Theta'(v_-)<0,
\eqs
for the excitatory case. The signs of $\A_\pm$ imply that $\T$ is Fredholm with Fredholm index $1$. The kernel is at most one-dimensional since solutions to the ODE are unique. Therefore, restricting to  $\X=\left\{\bu\in H^1~|~ \bu(0)=0 \right\}$, $\cL_{h,\zeta}$ yields an  invertible operator from $\X$ to $L^2$.

\subsubsection{Conclusion of the proof}

\begin{Proof}[ of Proposition \ref{prop:Invertibility}.]
Invertibility of  $\cL_{\epsilon,\zeta}:H^1\to L^2$ and $\cL_{h,\zeta}:\X\to L^2$, smallness of $\cL_{u,v}$, and boundedness of   $\cL_{u,v}$ give invertibility of $D\F^\epsilon(\bu_h,\bv_h)$.
We now show that $\left\|D\F^\epsilon(\bu_h,\bv_h)^{-1}\right\|$ is bounded uniformly in $\epsilon$ which will prove the second and last part of the proposition. Using Proposition \ref{prop:ConstInverse}, we have that $\cL_{\epsilon,\zeta}^{-1}=\cN\left(\I-\cR_{\epsilon,\zeta} \right)^{-1}$ and $\left\| \left(\I-\cR_{\epsilon,\zeta} \right)^{-1}\right\|\leq 1$. 
Since $\cN$ is bounded by \eqref{lem:nbd}, we find uniform bounds on $\cL_{\epsilon,\zeta}^{-1}$ as claimed. Differentiability with respect to parameters is a consequence of differentiability of the function $\F^\epsilon$. A simple bootstrap argument gives smoothness in $\zeta$.
\end{Proof}

\subsection{Invertibility in the fast component}\label{subsec:LdeInv}

In this section, we prove a complementary result that will be useful for the forthcoming sections. For $\epsilon>0$ and $c>0$, we will invert
\bqq
\label{eq:Lsl}
\cL_\epsilon:H^1\longrightarrow L^2,\qquad \left[\cL_\epsilon\bw\right](\xi)= -\bw(\xi)+\K\ast \bw(\xi)+c\frac{d}{d\xi}\bw(\xi)+f'\left(\bu_h(\epsilon \xi) \right)\bw(\xi).
\eqq

\begin{lem}\label{lem:LqeInv}
There exists $0<\eta_h<\eta_0$ and $\epsilon_0>0$ such that for all $|\eta|<\eta_h$, $0<\epsilon<\epsilon_0$ and $c>0$, $\cL_\epsilon$ is an isomorphism from $H^1_\eta$ to $L^2_{\eta}$ with $\epsilon$-uniform bounds on the inverse, depending continuously on $c,\eta$. 
\end{lem}

\begin{Proof}
Fix $\xi_0\in\R$ and consider again the frozen operator 
\bqs
\cL_\epsilon^{\xi_0}:H^1\longrightarrow L^2,\qquad \left[\cL_\epsilon^{\xi_0}\bw\right](\xi)= -\bw(\xi)+\K\ast \bw(\xi)+c\frac{d}{d\xi}\bw(\xi)+f'\left(\bu_h(\epsilon \xi_0) \right)\bw(\xi).
\eqs
Following the previous approach,  we define
\bqs
\Delta^{\epsilon\xi_0}(i\ell)=-1+\widehat{\K}(i\ell) + f'\left(\bu_h(\epsilon\xi_0)\right)+i \ell c,\qquad 
\G(\xi;\epsilon \xi_0)=\frac{1}{2\pi}\int_{\R}e^{i\ell\xi}\left[\Delta^{\epsilon\xi_0}(i\ell)\right]^{-1}d\ell.
\eqs
Since $\G$ is exponentially localized, there exists $0<\eta_h<\eta_0$ such that $\cL_\epsilon^{\xi_0}$ is an isomorphism from $H^1_\eta$ to $L^2_{\eta}$ for all $|\eta|<\eta_h$, with inverse given by convolution
\bqs
\left[\left(\cL_\epsilon^{\xi_0}\right)^{-1}\bh\right](\xi)= \left(\G \ast \bh\right)(\xi) = \int_\R \G(\xi-\tilde \xi;\epsilon \xi_0)\bh(\tilde\xi)d\tilde\xi.
\eqs
We introduce the linear operators $\cN$ and $\cR_\epsilon$, defined as
\bqs
\left[\cN \bh\right](\xi):= \int_\R \G(\xi-\tilde \xi;\epsilon\xi)\bh(\tilde\xi)d\tilde\xi, \quad \left[\cR_\epsilon \bh\right](\xi):=\left[\left(\I-\cL_\epsilon\cN\right)\bh\right](\xi)
\eqs
for all $\xi\in\R$. A direct computation shows that for all $\bh\in L^2_\eta$, we have
\bqs
\|\cR_{\epsilon}\bh\|_{L^2_{\eta}}\lesssim \epsilon \|\bh\|_{L^2_{\eta}}.
\eqs
Then, there exists $\epsilon>0$ such that for all $0<\epsilon<\epsilon_0$, $\cN \left(\I-\cR_\epsilon \right)^{-1}$ is well-defined from $L^2_\eta$ to $H^1_\eta$ and 
\bqs
\bu= \cN \left(\I-\cR_\epsilon \right)^{-1} \bh
\eqs
is a solution of $\cL_\epsilon\bu=\bh$. It is straightforward to check that $\cL_\epsilon$ is also one-to-one.
\end{Proof}

\section{Construction of the traveling pulse solution}\label{ansatz}

\subsection{The Ansatz}

\begin{figure}[htp]
\centering
\subfigure[]{
\label{fig:omega}
\includegraphics[width=0.65\textwidth]{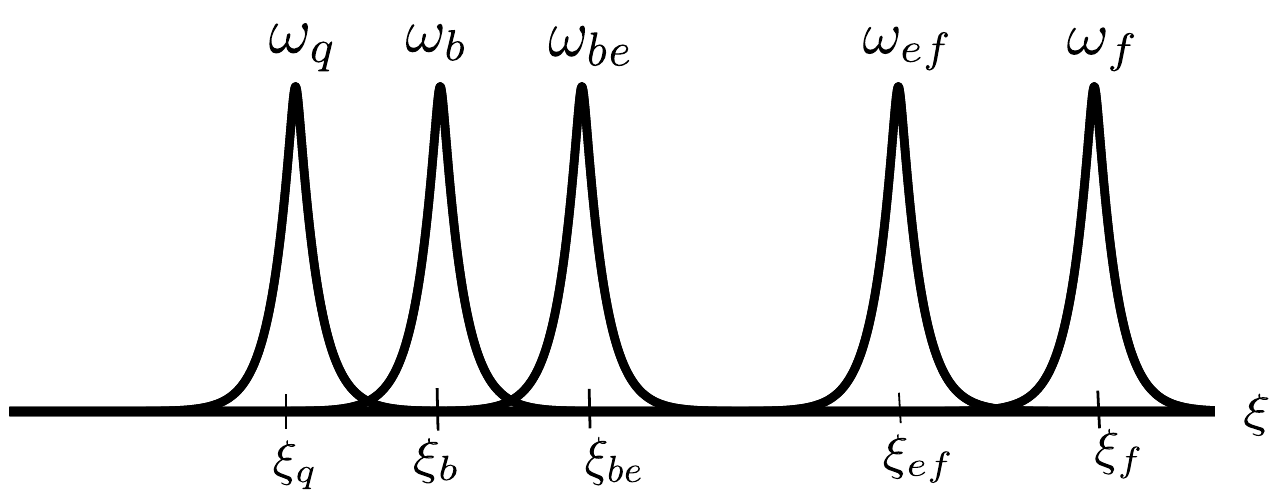}}
\hspace{.3in}
\subfigure[]{
\label{fig:Ucomp}
\includegraphics[width=0.65\textwidth]{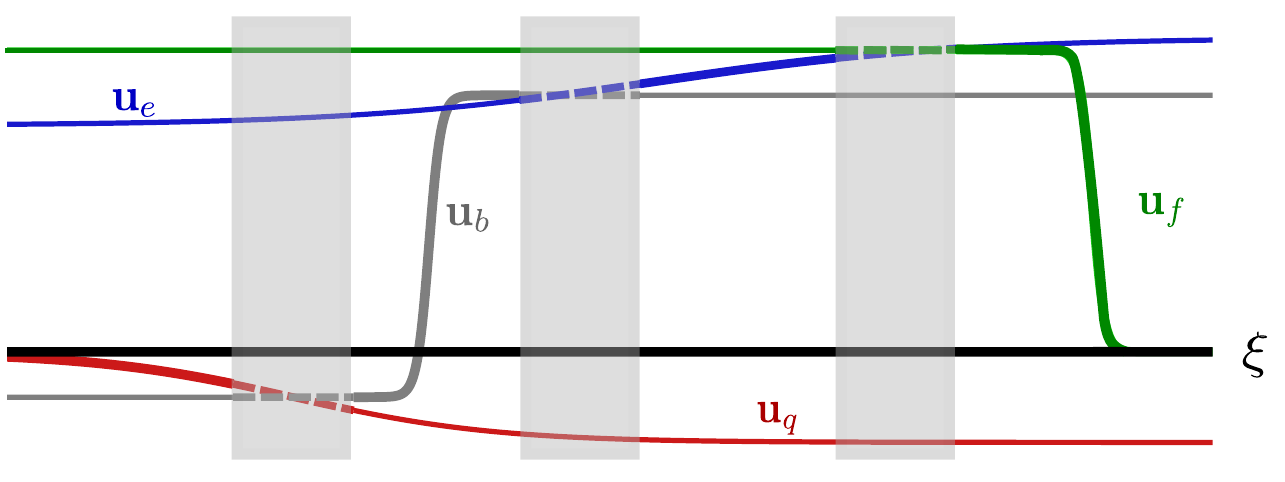}}
\subfigure[]{
\label{fig:Vcomp}
\includegraphics[width=0.65\textwidth]{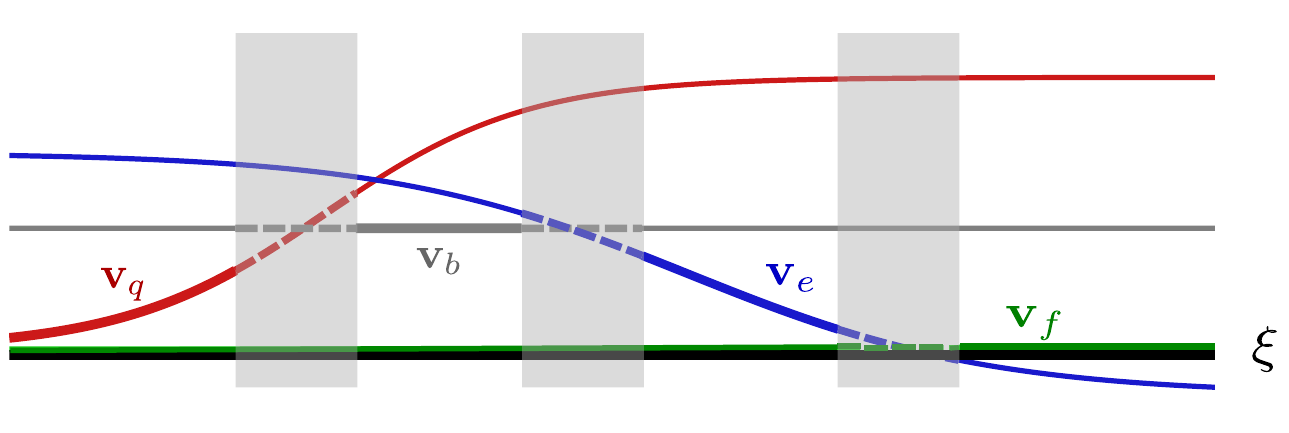}}
\hspace{.3in}
\subfigure[]{
\label{fig:Chi}
\includegraphics[width=0.65\textwidth]{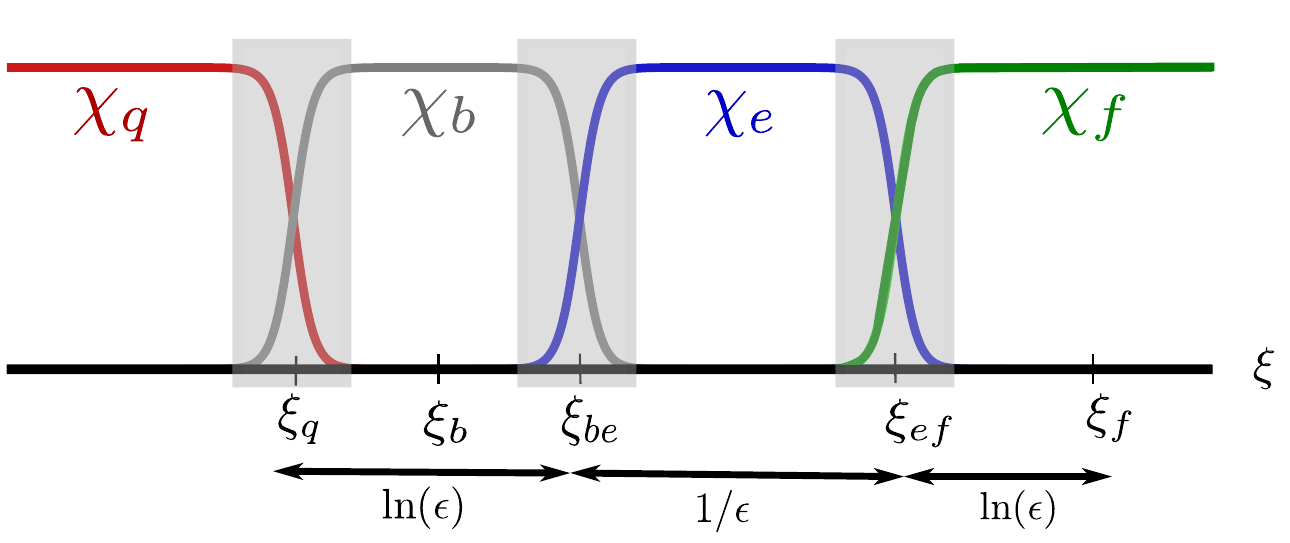}}
\caption{Schematic description of the Ansatz solution \eqref{eq:ansatz}. Envelopes $\omega_j$ for corrections $\bW_j$ as imposed by the weights $\omega_j^{-1}$ and $\bu_j$-components of the different parts of the Ansatz \eqref{eq:ansatz}. Profiles $\chi_j$ of the partition of unity as defined in \eqref{eq:chis}.}
\label{fig:ansatz}
\end{figure}

In the following, we present a decomposition of the solution into the singular pulse and corrections, separated using cut-off functions and exponentially localized weights. A schematic illustration of this procedure is shown in Figure \ref{fig:ansatz}. 

We write $\bU_f=\left(\bu_f,0\right)$ where $\bu_f\in \cC^\infty(\R,\R)$ is the front solution from Hypothesis (H3), solving
\bqq
\label{eq:Front}
-c_* \frac{d}{d\xi}\bu(\xi)=-\bu(\xi)+\int_{\R}\K(\xi-\xi')\bu(\xi')d\xi'+f(\bu(\xi)),
\eqq
with 
\bqs
\underset{\xi\rightarrow-\infty}{\lim}\bu_f(\xi)=1,\  \underset{\xi\rightarrow+\infty}{\lim}\bu_f(\xi)=0, \ \mbox{and }\bu_f(0)=\frac{1}{2} .
\eqs
Similarly, we set $\bU_b=\left(\bu_b,v_b\mathbf{1}\right)$ where $v_b\in (v_*-\delta_b,v_*+\delta_b)\subset (v_{min},v_{max})$ is a free parameter and  $\bu_b\in \cC^\infty(\R,\R)$ is the solution of 
\bqq
\label{eq:Back}
-c_b \frac{d}{d\xi}\bu(\xi)=-\bu(\xi)+\int_{\R}\K(\xi-\xi')\bu(\xi')d\xi'+f(\bu(\xi))-v_b
\eqq
with limits
\bqs
\underset{\xi\rightarrow-\infty}{\lim}\bu_b(\xi)=\varphi_q(v_b),\  \underset{\xi\rightarrow+\infty}{\lim}\bu_b(\xi)=\varphi_e(v_b) \ \mbox{and }\bu_b(0)=\left(\varphi_e(v_b)-\varphi_q(v_b)\right)/2.
\eqs
Again, this solution is obtained from Hypothesis (H3) for $v_b=v_*$. Using the implicit function theorem and simplicity of the zero eigenvalue, we can find the profile $\bu_b\in \mathcal{C}^\infty$ and the wave speed $c$ as smooth functions of $v_b\sim v_*$.

Using Proposition \ref{prop:slowsolL} and \ref{prop:slowsolR}, we define $\bU_{q}=\left(\bu_q,\bv_q\right)$ and $\bU_{e}=\left(\bu_e,\bv_e\right)$ where the heteroclinic solutions $\left(\bu_q(\epsilon \xi),\bv_q(\epsilon \xi) \right)$ and $\left(\bu_e(\epsilon \xi),\bv_e(\epsilon \xi) \right)$ solve
\begin{subequations}
\label{eq:SlowMod}
\begin{align}
-c \frac{d}{d\xi}\bu(\xi)&=-\bu(\xi)+\int_{\R}\K(\xi-\xi')\bu(\xi')d\xi'+f(\bu(\xi))- \bv(\xi) \\
-c \frac{d}{d\xi}\bv(\xi)&=\epsilon(\bu(\xi)-\gamma\bv(\xi))\Theta(\bv(\xi)),
\end{align}
\end{subequations}
with limits
\begin{align*}
\underset{\zeta\rightarrow-\infty}{\lim}\left(\bu_q(\zeta),\bv_q(\zeta)\right)=(0,0)&\text{ and } \underset{\zeta\rightarrow+\infty}{\lim}\left(\bu_q(\zeta),\bv_q(\zeta)\right)=(\varphi_q(v_+),v_+),\\
\underset{\zeta\rightarrow-\infty}{\lim}\left(\bu_e(\zeta),\bv_e(\zeta)\right)=(\varphi_e(v_+),v_+)&\text{ and } \underset{\zeta\rightarrow+\infty}{\lim}\left(\bu_e(\zeta),\bv_e(\zeta)\right)=(\varphi_e(v_-),v_-).
\end{align*}

Let $\delta_q>0$, $\delta_{be}>0$ and $\delta_{ef}>0$ be fixed such that $(v_*-\delta_q,v_*+\delta_q)\subset (v_{min},v_{max})$, $(v_*-\delta_{be},v_*+\delta_{be})\subset (v_{min},v_{max})$ and $(-\delta_{ef},\delta_{ef})\subset (v_{min},v_{max})$. We introduce three parameters $v_q\in (v_*-\delta_q,v_*+\delta_q)$, $v_{be}\in (v_*-\delta_{be},v_*+\delta_{be})$ and $v_{ef}\in(-\delta_{ef},\delta_{ef})$. We normalize the solutions $\bU_q$ and $\bU_e$ so that $\left(\bu_q(0),\bv_q(0)\right)=(\widetilde{\varphi}_q(v_q,\epsilon,c),v_q)$ and $\left(\bu_e(0),\bv_e(0)\right)=(\widetilde{\varphi}_e(v_{ef},\epsilon,c),v_{ef})$. Note that here we exploited monotonicity of the solution in the slow manifold and the implicit function theorem to normalize uniformly in the parameters
\bqs
\bu_j(0)=\widetilde{\varphi}_j(\bv_j(0),\epsilon,c),\quad j=q,e.
\eqs
We also define $T(v_{be},v_{ef})>0$ as the leading order time spent by $(\bu_e,\bv_e)$ on the excitatory slow manifold from $(\widetilde{\varphi}_e(v_{be},\epsilon,c),v_{be})$ to $(\widetilde{\varphi}_e(v_{ef},\epsilon,c),v_{ef})$. Note that $(v_{be},v_{ef})\longmapsto T(v_{be},v_{ef})$ is a continuously differentiable function on $(v_*-\delta_{be},v_*+\delta_{be})\times(-\delta_{ef},\delta_{ef})$.

We introduce a partition of unity through four $\cC^\infty$-functions $\chi_j$, $j\in\left\{q,b,e,f \right\}$, so that:
\bqq
\label{eq:chis}
\chi_q(\xi)+\chi_b(\xi)+\chi_e(\xi)+\chi_f(\xi)=1, \quad \forall \xi\in\R,
\eqq
and
\begin{align*}
\chi_q(\xi)=
\left\{
\begin{array}{ll}
0 & \xi \geq \xi_{q}+1 \\
1 & \xi \leq \xi_{q}-1
\end{array}
\right.,&\quad \chi_b(\xi)=
\left\{
\begin{array}{ll}
0 & \xi \leq \xi_{q}-1\text{ and } \xi \geq \xi_{be}+1 \\
1 & \xi_{q}+1 \leq \xi \leq \xi_{be}-1
\end{array}
\right.
,\\
\chi_f(\xi)=
\left\{
\begin{array}{ll}
0 & \xi \leq -1 \\
1 & \xi \geq 1
\end{array}
\right.
,& \quad
\chi_e(\xi)=
\left\{
\begin{array}{ll}
0 & \xi \leq \xi_{be}-1\text{ and } \xi \geq 1 \\
1 & \xi_{be}+1 \leq \xi \leq -1
\end{array}
\right..
\end{align*}
The constants $\xi_{q}$ and $\xi_{be}$ are defined as
\begin{align*}
\xi_{be}=-\frac{T(v_{be},v_{ef})}{\epsilon},\qquad 
\xi_{q}=\xi_{be}+2\eta_b\ln(\epsilon),
\end{align*}
where  $\eta_b>0$ will be fixed later. We will rely on the exponentially weighted spaces $H^1_\eta$, $L^2_\eta$ with
\begin{align*}
L^2_\eta=\left\{\bu:\R\rightarrow \R~|~ \left\| e^{\eta|\xi|}\bu(\xi) \right\|_{L^2}<+\infty \right\},\qquad 
H^1_\eta=\left\{\bu \in L^2_\eta ~|~\partial_\xi \bu \in L^2_\eta \right\},
\end{align*}
where again $\eta>0$ will be determined later. To find a pulse solution, we start with the following Ansatz
\begin{align}
\label{eq:ansatz}
\bU_a(\xi)=(\bu_a(\xi),\bv_a(\xi))=&~\bU_q\left(\epsilon\left(\xi-\xi_{q}\right)\right)\chi_q(\xi)+\bU_b(\xi-\xi_{b})\chi_b(\xi)+\bU_e(\epsilon\xi)\chi_e(\xi)+\bU_f(\xi-\xi_f)\chi_f(\xi)\nonumber\\
&+\bW_{q}\left(\xi-\xi_{q}\right)+\bW_b(\xi-\xi_{b})+\bW_{e}(\xi)+\bW_f(\xi-\xi_f),
\end{align}
where $\xi_b=\xi_{be}+\eta_b\ln(\epsilon)$, $\xi_f:=-\eta_f \ln(\epsilon)>0$, and 
$
\bW_j=(\bw_j^u,\bw_j^v), \quad j\in J_\bw:=\left\{q,b,e,f \right\}. 
$

The $\bu_a$ and $\bv_a$ components of $\bU_a$ are thus given by
\begin{subequations}
\label{eq:uv}
\begin{align}
\bu_a(\xi)&=\bu_q\left(\epsilon\left(\xi-\xi_{q}\right)\right)\chi_q(\xi)+\bu_b(\xi-\xi_{b})\chi_b(\xi)+\bu_e(\epsilon\xi)\chi_e(\xi)+\bu_f(\xi-\xi_f)\chi_f(\xi)\\
&\quad~+\bw_{q}^u\left(\xi-\xi_{q}\right)+\bw_b^u(\xi-\xi_{b})+\bw_{e}^u(\xi)+\bw_f^u(\xi-\xi_f),\nonumber\\
\bv_a(\xi)&=\bv_q\left(\epsilon\left(\xi-\xi_{q}\right)\right)\chi_q(\xi)+v_b\chi_b(\xi)+\bv_e(\epsilon\xi)\chi_e(\xi)+\bw_{q}^v\left(\xi-\xi_{q}\right)+\bw_b^v(\xi-\xi_{b})\\
&\quad~+\bw_{e}^v(\xi)+\bw_f^v(\xi-\xi_f).\nonumber
\end{align}
\end{subequations}

\begin{rmk}
We retain five free parameters $(c,v_q,v_b,v_{be},v_{ef})$;  $(\delta_q,\delta_b,\delta_{be},\delta_{ef},\eta_b,\eta_f)$ will be fixed later on.
\end{rmk}

\subsection{Deriving equations for the corrections $\bw_j$}

We substitute the expressions of $\bu_a$ and $\bv_a$ into \eqref{eq:TP} and obtain  equations for the corrections $\bw_j$. In the following, we will first make these equations explicit and then split the equations into a weakly coupled system of equations for the corrections $\bw_j$. 

The first component of our traveling-wave system gives the equation
\begin{align}
\label{eq:uExpNew}
0=&~\sum_{j\in J_\bw }\left(c\frac{d}{d\xi}\bw_j^u\left(\xi-\xi_j\right)-\bw_j^u\left(\xi-\xi_j\right)+\K\ast\bw_j^u(\xi-\xi_j)-\bw_j^v\left(\xi-\xi_j\right)\right)\nonumber\\
&+ \int_\R \K (\xi-\tilde \xi)\bu_q(\epsilon(\tilde\xi-\xi_{q}))\chi_q(\tilde\xi)d\tilde\xi -\chi_q(\xi)\int_\R \K (\xi-\tilde \xi)\bu_q(\epsilon(\tilde\xi-\xi_{q}))d\tilde\xi \nonumber\\
&+\int_\R \K (\xi-\tilde \xi)\bu_b(\tilde\xi-\xi_{b})\chi_b(\tilde\xi)d\tilde\xi -\chi_b(\xi)\int_\R \K (\xi-\tilde \xi)\bu_b(\tilde\xi-\xi_{b})d\tilde\xi \nonumber\\
&+\int_\R \K (\xi-\tilde \xi)\bu_e(\epsilon\tilde\xi)\chi_e(\tilde\xi)d\tilde\xi -\chi_e(\xi)\int_\R \K (\xi-\tilde \xi)\bu_e(\epsilon\tilde\xi)d\tilde\xi \nonumber\\
&+\int_\R \K (\xi-\tilde \xi)\bu_f(\tilde\xi-\xi_f)\chi_f(\tilde\xi)d\tilde\xi -\chi_f(\xi)\int_\R \K (\xi-\tilde \xi)\bu_f(\tilde\xi-\xi_f)d\tilde\xi \nonumber\\
&+c\bu_q\left(\epsilon\left(\xi-\xi_{q}\right)\right)\chi_q'(\xi)+c\bu_b(\xi-\xi_{b})\chi_b'(\xi)+c\bu_e(\epsilon\xi)\chi_e'(\xi)+c\bu_f(\xi-\xi_f)\chi_f'(\xi)\nonumber\\
&+(c-c_*)\bu_f'(\xi-\xi_f)\chi_f(\xi)+(c-c_b)\bu_b'(\xi-\xi_b)\chi_b(\xi)+f\left( \bu_a(\xi)\right)\nonumber\\
&-f(\bu_q\left(\epsilon\left(\xi-\xi_{q}\right)\right))\chi_q(\xi)-f(\bu_b(\xi-\xi_{b}))\chi_b(\xi)-f(\bu_e(\epsilon\xi))\chi_e(\xi)-f(\bu_f(\xi-\xi_f))\chi_f(\xi);
\end{align}
and the second component yields
\begin{align}
\label{eq:vExpNew}
0=&~\sum_{j\in J_\bw }\left( c\frac{d}{d\xi}\bw_j^u\left(\xi-\xi_j\right)+\epsilon \left(\bw_j^u\left(\xi-\xi_j\right)-\gamma \bw_j^v\left(\xi-\xi_j\right)\right)\right)+\epsilon\left(\bu_b(\xi-\xi_b)-\gamma v_b \right)\chi_b(\xi)\nonumber\\
&+\epsilon\bu_f(\xi-\xi_f)\chi_f(\xi)+c\bv_q(\epsilon(\xi-\xi_{q}))\chi_q'(\xi)+cv_b\chi'_b(\xi)+c\bv_e(\epsilon\xi)\chi'_e(\xi)\nonumber\\
&+\epsilon\chi_q(\xi)\left(\bu_q(\epsilon(\xi-\xi_{q}))-\gamma \bv_q(\epsilon (\xi-\xi_{q}))\right)\left(1-\Theta(\bv_q(\epsilon (\xi-\xi_{q})))\right)\nonumber\\
&+\epsilon\chi_e(\xi)\left(\bu_e(\epsilon\xi)-\gamma \bv_e(\epsilon \xi)\right)\left(1-\Theta(\bv_e(\epsilon \xi))\right).
\end{align}

We now split this system into five systems, one for each $\bW_j=(\bw_j^u,\bw_j^v)$, $j\in J_\bw$. The right-hand sides of \eqref{eq:uExpNew} and \eqref{eq:vExpNew} will be the sum of the right-hand sides of those equations for the $\bW_j$, below, so that solving the equations for the $\bW_j$ will automatically give us a solution to \eqref{eq:uExpNew} and \eqref{eq:vExpNew}. 

We first introduce some notation in order to facilitate the presentation of these systems. We Taylor expand the nonlinear term $f(\bu_a(\xi))$ in \eqref{eq:uExpNew} at 
\bqs
\bu_0(\xi):=\bu_q\left(\epsilon\left(\xi-\xi_{q}\right)\right)\chi_q(\xi)+\bu_b(\xi-\xi_{b})\chi_b(\xi)+\bu_e(\epsilon\xi)\chi_e(\xi)+\bu_f(\xi-\xi_f)\chi_f(\xi),
\eqs
and get
\bqs
f\left( \bu_a(\xi)\right)=f\left( \bu_0(\xi)\right)+f'\left( \bu_0(\xi)\right)\sum_{j\in J_\bw}\bw_j^u(\xi-\xi_j)+\sum_{j \leq k \in \J_\bw} \bQ_{j,k}(\xi) \bw_j^u(\xi-\xi_j)\bw_k^u(\xi-\xi_k),\eqs
where $\bQ_{j,k}:=\bQ_{j,k}\left(\bu_0,\bw_{q}^u,\bw_{b}^u,\bw_{e}^u,\bw_{f}^u\right)$.  There exist constants $C_{j,k}$, independent of $\epsilon$, such that
\bqs
\left\| \bQ_{j,k}\right\|_{L^\infty(\R)}\leq C_{j,k} \text{ as } \left\|(\bw_{q}^u,\bw_{b}^u,\bw_{e}^u,\bw_{f}^u)\right\|\rightarrow 0.
\eqs

If we define the function $\bF$ as:
\begin{align}
\bF(\xi)=&~ f(\bu_0(\xi))-f(\bu_q\left(\epsilon\left(\xi-\xi_{q}\right)\right))\chi_q(\xi)-f(\bu_b(\xi-\xi_{b}))\chi_b(\xi)-f(\bu_e(\epsilon\xi))\chi_e(\xi)\nonumber\\
&-f(\bu_f(\xi-\xi_f))\chi_f(\xi),\label{eq:ComF}
\end{align}
for $\xi\in\R$, then a direct computation shows that we have
\bqq
\label{eq:ComDefF}
\bF(\xi)=\bF_q(\xi)\mathds{1}_{\xi_q}(\xi)+\bF_{be}(\xi)\mathds{1}_{\xi_{be}}(\xi)+\bF_{ef}(\xi)\mathds{1}_{\xi_{ef}}(\xi),
\eqq
where
\begin{subequations}
\label{eq:ComFj}
\begin{align}
\bF_q(\xi)&= f\left((\bu_q\left(\epsilon\left(\xi-\xi_{q}\right)\right)\chi_q(\xi)+\bu_b(\xi-\xi_{b})\chi_b(\xi)\right)-f(\bu_q\left(\epsilon\left(\xi-\xi_{q}\right)\right))\chi_q(\xi)-f(\bu_b(\xi-\xi_{b}))\chi_b(\xi),\\
\bF_{be}(\xi)&=f\left(\bu_b(\xi-\xi_{b})\chi_b(\xi)+\bu_e(\epsilon\xi)\chi_e(\xi) \right)-f(\bu_b(\xi-\xi_{b}))\chi_b(\xi)-f(\bu_e(\epsilon\xi))\chi_e(\xi),\\
\bF_{ef}(\xi)&=f\left( \bu_e(\epsilon\xi)\chi_e(\xi)+\bu_f(\xi-\xi_f)\chi_f(\xi)\right)-f(\bu_e(\epsilon\xi))\chi_e(\xi)-f(\bu_f(\xi-\xi_f))\chi_f(\xi),
\end{align}
\end{subequations}
and $\mathds{1}_{\xi_j}$ stands for the indicator function of the interval $[\xi_j-1,\xi_j+1]$. We also define 
\begin{subequations}
\label{eq:ComC}
\begin{align}
\bC(\xi)=&~ \K \ast \left[\bu_0(\xi) \right]-\chi_q(\xi)\K \ast \left[\bu_q(\epsilon(\xi-\xi_q)) \right] -\chi_b(\xi) \K \ast \left[\bu_b(\xi-\xi_b) \right] -\chi_e(\xi) \K \ast  \left[\bu_e(\epsilon\xi) \right]\nonumber\\
& -\chi_f(\xi)\K\ast \left[\bu_f(\xi-\xi_f) \right],
\end{align}
\end{subequations}
for all $\xi\in\R$. We denote $\chi_{bq}:=\chi_b+\chi_q$ and $\chi_{ef}:=\chi_e+\chi_f$ so that we have
$
\chi_{bq}(\xi)+\chi_{ef}(\xi)=1, \quad \forall \xi \in \R.
$
Finally, a direct computation shows that we have
\bqq
\label{eq:ComDefC}
\bC(\xi)=\bC_q(\xi)+\bC_{be}(\xi)+\bC_{ef}(\xi),
\eqq
where
\begin{subequations}
\label{eq:ComCj}
\begin{align}
\bC_q(\xi)&= \K\ast \left[ \left(u_q(\epsilon(\xi-\xi_q))-u_b(\xi-\xi_b) \right)\chi_q(\xi) \right]-\chi_q(\xi)\K\ast \left[ u_q(\epsilon(\xi-\xi_q))-u_b(\xi-\xi_b) \right],\\
\bC_{be}(\xi)&=\K\ast \left[ \left(u_b(\xi-\xi_b)-u_e(\epsilon\xi) \right)\chi_{ef}(\xi) \right]-\chi_{ef}(\xi)\K\ast \left[ u_b(\xi-\xi_b)-u_e(\epsilon\xi) \right] ,\\
\bC_{ef}(\xi)&=  \K\ast \left[ \left(u_f(\xi-\xi_f)-u_e(\epsilon\xi) \right)\chi_f(\xi) \right]-\chi_f(\xi)\K\ast \left[ u_f(\xi-\xi_f)-u_e(\epsilon\xi) \right].
\end{align}
\end{subequations}

We are now ready to present the explicit form of the equations for the $\bW_j$. 

\paragraph{\textbf{Equations for the quiescent part:}}
\begin{subequations}
\label{eq:Wq}
\begin{align}
0=&~c\frac{d}{d\xi}\bw_q^u(\xi)-\bw_q^u(\xi)+\K\ast \bw_q^u(\xi)+f'\left(\bu_q(\epsilon\xi) \right)\bw_q^u(\xi)- \bw_q^v(\xi)+\sum_{j\in J_\bw}\bL_{q,j}(\xi+\xi_q)\bw_j^u(\xi-\xi_j+\xi_q)\nonumber\\
&+c\bu_q\left(\epsilon\left(\xi\right)\right)\chi_q'(\xi+\xi_q)+c\bu_b(\xi-\xi_{b}+\xi_q)\chi_b'(\xi+\xi_q)\mathds{1}_{\xi_b}(\xi+\xi_q)\nonumber\\
&+\bF_q(\xi+\xi_q)\mathds{1}_{\xi_q}(\xi+\xi_q)+\bC_q(\xi+\xi_q) +\bQ_{q,q}(\xi+\xi_q)\left(\bw_q^u(\xi)\right)^2\nonumber\\
&+\sum_{j\neq k \in J_\bw}\bQ_{j,k}(\xi+\xi_q)\bw_j^u(\xi-\xi_j+\xi_q)\bw_k^u(\xi-\xi_k+\xi_q)\chi_q(\xi+\xi_q),\label{eq:Wqu}\\
0=&~c\frac{d}{d\xi}\bw_q^v(\xi)+ \epsilon \left(\bw_q^u(\xi)-\gamma \bw_q^v(\xi) \right)+c\bv_q(\epsilon(\xi))\chi_q'(\xi+\xi_q)+cv_b\chi'_b(\xi+\xi_q)\mathds{1}_{\xi_b}(\xi+\xi_q)\nonumber\\
&+\epsilon\chi_q(\xi+\xi_q)\left(\bu_q(\epsilon \xi)-\gamma\bv_q(\epsilon \xi)\right)\left(1-\Theta(\bv_q(\epsilon \xi))\right)\nonumber\\
&+\epsilon \left(\int_\R \left(\bw_b^u(\xi)-\gamma\bw_b^v(\xi) + (\bu_b(\xi)-\gamma v_b)\chi_b(\xi+\xi_b)\right)d\xi\right)\psi\left(\xi-\frac{\xi_b-\xi_q}{2}\right).\label{eq:Wqv}
\end{align}
\end{subequations}

\paragraph{\textbf{Equations for the back part:}}
\begin{subequations}
\label{eq:Wb}
\begin{align}
0=&~c\frac{d}{d\xi}\bw_b^u(\xi)-\bw_b^u(\xi)+\K\ast \bw_b^u(\xi)+f'\left(\bu_b(\xi) \right)\bw_b^u(\xi)- \bw_b^v(\xi)\nonumber\\
&+\sum_{j\in J_\bw}\bL_{b,j}(\xi+\xi_b)\bw_j^u(\xi-\xi_j+\xi_b)+(c-c_b)\bu_b'(\xi)\chi_b(\xi+\xi_b)+\bQ_{b,b}(\xi+\xi_b)\left(\bw_b^u(\xi)\right)^2\nonumber\\
&+\sum_{j\neq k \in J_\bw}\bQ_{j,k}(\xi+\xi_b)\bw_j^u(\xi-\xi_j+\xi_b)\bw_k^u(\xi-\xi_k+\xi_b)\chi_b(\xi+\xi_b),\label{eq:Wbu}\\
0=&~c\frac{d}{d\xi}\bw_b^v(\xi)+ \epsilon \left(\bw_b^u(\xi)-\gamma \bw_b^v(\xi) \right)+\epsilon \left(\bu_b(\xi)-\gamma v_b\right)\chi_b(\xi+\xi_b)\nonumber\\
&-\epsilon \left(\int_\R \left(\bw_b^u(\xi)-\gamma\bw_b^v(\xi) + \left(\bu_b(\xi)-\gamma v_b\right)\chi_b(\xi+\xi_b)\right)d\xi\right)\psi\left(\xi+\frac{\xi_b-\xi_q}{2}\right).\label{eq:Wbv}
\end{align}
\end{subequations}

\paragraph{\textbf{Equations for the excitatory part:}}

\begin{subequations}
\label{eq:We}
\begin{align}
0=&~c\frac{d}{d\xi}\bw_e^u(\xi)-\bw_e^u(\xi)+\K\ast \bw_e^u(\xi)+f'\left(\bu_e(\epsilon\xi) \right)\bw_e^u(\xi)- \bw_e^v(\xi)+\sum_{j\in J_\bw}\bL_{e,j}(\xi)\bw_j^u(\xi-\xi_j)\nonumber\\
&+c\bu_b(\xi-\xi_{b})\chi_b'(\xi)\mathds{1}_{\xi_{be}}(\xi)+c\bu_e(\epsilon\xi)\chi_e'(\xi)+c\bu_f(\xi-\xi_f)\chi_f'(\xi) \nonumber\\
&+\bF_{be}(\xi)\mathds{1}_{\xi_{be}}(\xi)+\bF_{ef}(\xi)\mathds{1}_{\xi_{ef}}(\xi)+\bC_{be}(\xi)+\bC_{ef}(\xi)\nonumber\\
&+\bQ_{e,e}(\xi)\left(\bw_e^u(\xi)\right)^2+\sum_{j\neq k \in J_\bw}\bQ_{j,k}(\xi)\bw_j^u(\xi-\xi_j)\bw_k^u(\xi-\xi_k)\chi_e(\xi),\label{eq:Weu}\\
0=&~c\frac{d}{d\xi}\bw_e^v(\xi)+ \epsilon \left(\bw_e^u(\xi)-\gamma \bw_e^v(\xi) \right)
+c\bv_e(\epsilon\xi)\chi_e'(\xi)+\epsilon\chi_e(\xi)\left(\bu_e(\epsilon \xi)-\gamma\bv_e(\epsilon \xi)\right)\left(1-\Theta(\bv_e(\epsilon \xi))\right)\nonumber\\
&+\epsilon \left(\int_\R \left(\bw_f^u(\xi)-\gamma \bw_f^v(\xi) + \bu_f(\xi)\chi_f(\xi+\xi_f)\right)d\xi\right)\psi(\xi-\xi_f/2). \label{eq:Wev}
\end{align}
\end{subequations}

\paragraph{\textbf{Equations for the front part:}}
\begin{subequations}
\label{eq:Wf}
\begin{align}
0=&~c\frac{d}{d\xi}\bw_f^u(\xi)-\bw_f^u(\xi)+\K\ast \bw_f^u(\xi)+f'\left(\bu_f(\xi) \right)\bw_f^u(\xi)- \bw_f^v(\xi)\nonumber\\
&+\sum_{j\in J_\bw}\bL_{f,j}(\xi+\xi_f)\bw_j^u(\xi-\xi_j+\xi_f)+(c-c_*)\bu_f'(\xi)\chi_f(\xi+\xi_f)+\bQ_{f,f}(\xi+\xi_f)\left(\bw_f^u(\xi)\right)^2\nonumber\\
&+\sum_{j\neq k \in J_\bw}\bQ_{j,k}(\xi+\xi_f)\bw_j^u(\xi-\xi_j+\xi_f)\bw_k^u(\xi-\xi_k+\xi_f)\chi_f(\xi+\xi_f),\label{eq:Wfu} \\ 
0=&~c\frac{d}{d\xi}\bw_f^v(\xi)+ \epsilon \left(\bw_f^u(\xi)-\gamma \bw_f^v(\xi) \right)+\epsilon \bu_f(\xi)\chi_f(\xi+\xi_f)\nonumber\\
&-\epsilon \left(\int_\R \left(\bw_f^u(\xi)-\gamma\bw_f^v(\xi) + \bu_f(\xi)\chi_f(\xi+\xi_f)\right)d\xi\right)\psi(\xi+\xi_f/2). \label{eq:Wfv}
\end{align}
\end{subequations}

The linear terms $\bL_{k,j}(\xi)$ that appear in systems \eqref{eq:Wq}, \eqref{eq:Wb}, \eqref{eq:We} and \eqref{eq:Wf} are defined as follows. For all $j\in J_\bw$, the diagonal terms are equal:
\begin{align}
\label{eq:Ldiag}
\bL_d(\xi):=\bL_{j,j}(\xi)=&~f'(\bu_0(\xi))-f'(\bu_q\left(\epsilon\left(\xi-\xi_{q}\right)\right))\chi_q(\xi)-f'(\bu_b(\xi-\xi_{b}))\chi_b(\xi)-f'(\bu_e(\epsilon\xi))\chi_e(\xi)\nonumber\\
&-f'(\bu_f(\xi-\xi_f))\chi_f(\xi).
\end{align}
We also have for the quiescent part:
\begin{subequations}
\label{eq:Lqj}
\begin{align}
\bL_{q,b}(\xi)=&~\chi_q(\xi)\left(f'(\bu_q(\epsilon(\xi-\xi_q))-f'(\bu_b(\xi-\xi_b)) \right),\\\
\bL_{q,e}(\xi)=&~\chi_q(\xi)\left(f'(\bu_q(\epsilon(\xi-\xi_q))-f'(\bu_e(\epsilon\xi)) \right),\\
\bL_{q,f}(\xi)=&~\chi_q(\xi)\left(f'(\bu_q(\epsilon(\xi-\xi_q))-f'(\bu_f(\xi-\xi_f)) \right),
\end{align}
\end{subequations}
for the back part:
\begin{subequations}
\label{eq:Lbj}
\begin{align}
\bL_{b,q}(\xi)=&~\chi_b(\xi)\left(f'(\bu_b(\xi-\xi_b))-f'(\bu_q(\epsilon(\xi-\xi_q))) \right),\\\
\bL_{b,e}(\xi)=&~\chi_b(\xi)\left(f'(\bu_b(\xi-\xi_b))-f'(\bu_e(\epsilon\xi)) \right),\\
\bL_{b,f}(\xi)=&~\chi_b(\xi)\left(f'(\bu_b(\xi-\xi_b))-f'(\bu_f(\xi-\xi_f)) \right),
\end{align}
\end{subequations}
for the excitatory part:
\begin{subequations}
\label{eq:Lej}
\begin{align}
\bL_{e,q}(\xi)=&~\chi_e(\xi)\left(f'(\bu_e(\epsilon\xi))-f'(\bu_q(\epsilon(\xi-\xi_q))) \right),\\\
\bL_{e,b}(\xi)=&~\chi_e(\xi)\left(f'(\bu_e(\epsilon\xi))-f'(\bu_b(\xi-\xi_b)) \right),\\
\bL_{e,f}(\xi)=&~\chi_e(\xi)\left(f'(\bu_e(\epsilon\xi))-f'(\bu_f(\xi-\xi_f)) \right),
\end{align}
\end{subequations}
and for the front part
\begin{subequations}
\label{eq:Lfj}
\begin{align}
\bL_{f,q}(\xi)=&~\chi_f(\xi)\left(f'(\bu_f(\xi-\xi_f))-f'(\bu_q(\epsilon(\xi-\xi_q))) \right),\\\
\bL_{f,b}(\xi)=&~\chi_f(\xi)\left(f'(\bu_f(\xi-\xi_f))-f'(\bu_b(\xi-\xi_b)) \right),\\
\bL_{f,e}(\xi)=&~\chi_f(\xi)\left(f'(\bu_f(\xi-\xi_f))-f'(\bu_e(\epsilon\xi)) \right).
\end{align}
\end{subequations}

The function $\psi:\R\rightarrow \R$, that appears in equations \eqref{eq:Wqv}, \eqref{eq:Wbv}, \eqref{eq:Wev} and \eqref{eq:Wfv}, is chosen to be $\cC^\infty$, exponentially localized around $\xi=0$ with compact support, and such that
\bqs
\int_\R \psi(\xi)d\xi =1 \text{ and } \left\| \psi \right\| _{L^2_\eta} <\infty,\quad \forall\eta>0.
\eqs
It effectively shifts mass between different components $\bw_j$. In particular, our choice of $\psi$ guarantees that \eqref{eq:Wqv} is satisfied upon integration. Anticipating some of the later analysis, we remark that the operator $\frac{d}{d\xi}$ which appears in \eqref{eq:Wqv} and \eqref{eq:Wbv} possesses a cokernel spanned by the constant functions. In the original systems, compensating for this cokernel requires one additional parameter. Splitting the equations into different components artificially inflates this cokernel, and we compensate for this fact by artificially transferring mass between the different parts of the system.

From the above definition of $\psi$, we directly have the estimates:
\begin{align}
\left\| \psi\left(\cdot \pm \frac{\xi_b-\xi_q}{2}\right) \right\| _{L^2_\eta} &\lesssim \epsilon^{-\frac{\eta\eta_b}{2} } \left\| \psi \right\| _{L^2_\eta},\label{eq:EstPsib}\\
\left\| \psi\left(\cdot \pm \frac{\xi_f}{2}\right) \right\| _{L^2_\eta} &\lesssim \epsilon^{- \frac{\eta\eta_f}{2} } \left\| \psi \right\| _{L^2_\eta},\label{eq:EstPsif}
\end{align}
as $\epsilon \rightarrow 0$. If we suppose that $\bw_j^u$ and $\bw_j^v$ belong to $L^2_\eta$ for $j\in\left\{ q,f \right\}$, then we have that 
\bqs
\int_\R \bw_j^u(\xi)-\gamma \bw_j^v(\xi)d\xi <\infty, \quad j\in\left\{ q,f \right\}.
\eqs
On the other hand, we have that
\begin{align*}
\left|\int_\R \left(\bu_b(\xi)-\gamma v_b\right)\chi_b(\xi+\xi_b)d\xi \right|&= \cO(|\ln \epsilon|),\\
\left|\int_\R \bu_f(\xi)\chi_f(\xi+\xi_f)d\xi \right|&= \cO(|\ln \epsilon|),
\end{align*}
for $\epsilon \rightarrow 0$.
As a consequence we obtain
\begin{subequations}
\begin{align}
\left\| \epsilon \left(\int_\R \left(\bw_b^u(\xi)-\gamma \bw_b^v(\xi) +\left(\bu_b(\xi)-\gamma v_b\right)\chi_b(\xi+\xi_b) \right)d\xi\right) \psi\left(\cdot \pm  \frac{\xi_b-\xi_q}{2}\right) \right\| _{L^2_\eta}&= \cO\left( \epsilon^{1-\frac{ \eta \eta_b}{2} }|\ln \epsilon | \right),\label{eq:EstB}\\
\left\| \epsilon \left(\int_\R \left(\bw_f^u(\xi)-\bw_f^v(\xi) + \bu_f(\xi)\chi_f(\xi+\xi_f)\right)d\xi\right) \psi\left(\cdot \pm \frac{\xi_f}{2}\right) \right\| _{L^2_\eta}&= \cO\left( \epsilon^{1-\frac{ \eta\eta_f}{2} }|\ln \epsilon | \right),\label{eq:EstF}
\end{align}
\end{subequations}
and, provided that $1- \eta \eta_b/2>0$ and $1- \eta \eta_f/2$, the above quantities are small as $\epsilon \rightarrow 0$. From now on, we assume that $\eta < \min\left(\frac{2}{\eta_b},\frac{2}{\eta_f} \right)$. 

The correction of the excitatory part $\bW_e$ is commonly constructed using the Exchange Lemma in a dynamical systems based approach \cite{jonesreview}. Those corrections are exponentially localized close to touchdown and takeoff points. Rather than encoding this localization at two diverging points $0$ and $\xi_{be}$, with a varying family of weights, we prefer to again split $\bW_e$ into  $\bW_{be}$ and $\bW_{ef}$, 
\bqs
\bW_e(\xi)=\bW_{be}(\xi-\xi_{be})+\bW_{ef}(\xi)
\eqs 
with $\bW_{be}=(\bw_{be}^u,\bw_{be}^v)\in H^1_\eta\times H^1_\eta$ and $\bW_{ef}=(\bw_{ef}^u,\bw_{ef}^v)\in H^1_\eta\times H^1_\eta$. Again, we separate the system  \eqref{eq:We} into two parts, as follows.

\paragraph{\textbf{Equations for the back/excitatory part:}}

\begin{subequations}
\label{eq:Wbe}
\begin{align}
0=&~c\frac{d}{d\xi}\bw_{be}^u(\xi)-\bw_{be}^u(\xi)+\K\ast \bw_{be}^u(\xi)+f'\left(\bu_e(\epsilon\xi+\epsilon \xi_{be}) \right)\bw_{be}^u(\xi)- \bw_{be}^v(\xi)\nonumber\\
&+\sum_{j\in \widetilde{J}_\bw}\bL_{be,j}(\xi+\xi_{be})\bw_j^u(\xi-\xi_j+\xi_{be})+c\bu_b(\xi-\xi_{b}+\xi_{be})\chi_b'(\xi+\xi_{be})\mathds{1}_{\xi_{be}}(\xi+\xi_{be})\nonumber\\
&+c\bu_e(\epsilon\xi+\epsilon\xi_{be})\chi_e'(\xi+\xi_{be})\mathds{1}_{\xi_{be}}(\xi+\xi_{be})+\bF_{be}(\xi+\xi_{be})\mathds{1}_{\xi_{be}}(\xi+\xi_{be})+\bC_{be}(\xi+\xi_{be})\nonumber\\
&+\bQ_{be,be}(\xi)\left(\bw_{be}^u(\xi)\right)^2+\sum_{j\neq k \in \widetilde{J}_\bw}\bQ_{j,k}(\xi)\bw_j^u(\xi-\xi_j+\xi_{be})\bw_k^u(\xi-\xi_k+\xi_{be})\chi_e(\xi+\xi_{be}),\label{eq:Wbeu}\\
0=&~c\frac{d}{d\xi}\bw_{be}^v(\xi)+ \epsilon \left(\bw_{be}^u(\xi)-\gamma \bw_{be}^v(\xi) \right)+cv_b\chi'_b(\xi+\xi_{be})\mathds{1}_{\xi_{be}}(\xi+\xi_{be})\nonumber\\
&+c\bv_e(\epsilon\xi+\epsilon\xi_{be})\chi_e'(\xi+\xi_{be})\mathds{1}_{\xi_{be}}(\xi+\xi_{be})\nonumber\\
&+\epsilon\chi_e(\xi+\xi_{be})\left(\bu_e(\epsilon \xi+\epsilon\xi_{be})-\gamma\bv_e(\epsilon \xi+\epsilon\xi_{be})\right)\left(1-\Theta(\bv_e(\epsilon \xi+\epsilon\xi_{be}))\right)\mathds{1}_{\xi_{be}}(\xi+\xi_{be}). \label{eq:Wbev}
\end{align}
\end{subequations}

\paragraph{\textbf{Equations for the excitatory/front part:}}

\begin{subequations}
\label{eq:Wef}
\begin{align}
0=&~c\frac{d}{d\xi}\bw_{ef}^u(\xi)-\bw_{ef}^u(\xi)+\K\ast \bw_{ef}^u(\xi)+f'\left(\bu_e(\epsilon\xi) \right)\bw_{ef}^u(\xi)- \bw_{ef}^v(\xi)+\sum_{j\in \widetilde{J}_\bw}\bL_{ef,j}(\xi)\bw_j^u(\xi-\xi_j)\nonumber\\
&+c\bu_e(\epsilon\xi)\chi_e'(\xi)\mathds{1}_{\xi_{ef}}(\xi)+c\bu_f(\xi-\xi_f)\chi_f'(\xi)+\bF_{ef}(\xi)\mathds{1}_{\xi_{ef}}(\xi)+\bC_{ef}(\xi) \nonumber\\
&+\bQ_{ef,ef}(\xi)\left(\bw_{ef}^u(\xi)\right)^2+\sum_{j\neq k \in \widetilde{J}_\bw}\bQ_{j,k}(\xi)\bw_j^u(\xi-\xi_j)\bw_k^u(\xi-\xi_k)\chi_e(\xi),\label{eq:Wefu}\\
0=&~c\frac{d}{d\xi}\bw_{ef}^v(\xi)+ \epsilon \left(\bw_{ef}^u(\xi)-\gamma \bw_{ef}^v(\xi) \right)+c\bv_e(\epsilon\xi)\chi_e'(\xi)\mathds{1}_{\xi_{ef}}(\xi)\nonumber\\
&+\epsilon\chi_e(\xi)\left(\bu_e(\epsilon \xi)-\gamma\bv_e(\epsilon \xi)\right)\left(1-\Theta(\bv_e(\epsilon \xi))\right)\mathds{1}_{\xi_{ef}}(\xi) \nonumber\\
&+\epsilon \left(\int_\R \left(\bw_f^u(\xi)-\gamma \bw_f^v(\xi) + \bu_f(\xi)\chi_f(\xi+\xi_f)\right)d\xi\right)\psi(\xi-\xi_f/2). \label{eq:Wefv}
\end{align}
\end{subequations}

\subsection{Formulation of the problem}

To conclude the setup, we rewrite systems \eqref{eq:Wq}, \eqref{eq:Wb}, \eqref{eq:Wbe}, \eqref{eq:Wef}, and \eqref{eq:Wf} in the more general and compact form:
\bqq
\label{eq:Compact}
0=\mathcal{L}_\epsilon \left(\bW,\lambda-\lambda_*\right) +\mathcal{R}_\epsilon + \mathcal{N}_\epsilon(\bW,\lambda-\lambda_*):=\F_\epsilon(\bW,\lambda),
\eqq
where we have set
\begin{align*}
\bW&:=\left(\left(\bw_q^u,\bw_q^v \right),\left(\bw_b^u,\bw_b^v \right),\left( \bw_{be}^u,\bw_{be}^v\right),\left(\bw_{ef}^u,\bw_{ef}^v \right),\left( \bw_f^u,\bw_f^v\right) \right)\in \X:=\left(H^1_\eta \times H^1_\eta \right)^5,\\
\lambda&:=\left(c,v_q,v_b,v_{be},v_{ef} \right)\in \V,\\
\lambda_*&:=\left(c_*,v_*,v_*,v_*,0 \right).
\end{align*}
Here, $\cL_\epsilon$ represents all the linear terms, $\mathcal{R}_\epsilon$ collects all the error terms and $\cN_\epsilon$ all the nonlinear terms. We define the nonlinear map $\F_\epsilon$  as follows
\bqq
\label{eq:CompactF}
\F_\epsilon: \begin{array}{ccl}
\X\times\V &\longrightarrow& \Y \\
 \left(\bW,\lambda \right)&\longmapsto& \F_\epsilon(\bW,\lambda)
 \end{array}
\eqq
where $\X:=\left(H^1_\eta \times H^1_\eta \right)^5$, $\Y:=\left(L^2_\eta \times L^2_\eta \right)^5$ and $\V:=(c_*-\delta_c,c_*+\delta_c)\times (v_*-\delta_b,v_*+\delta_b)\times (v_*-\delta_b,v_*+\delta_b)\times (v_*-\delta_{be},v_*+\delta_{be})\times (-\delta_{ef},\delta_{ef})$ is a neighborhood of $\lambda_*=(c_*,v_*,v_*,v_*,0)$ in $\R^5$. In the following sections, our strategy will be to show that 
\begin{enumerate}
\item the map $\F_\epsilon$ is well-defined from $\X\times\V$ to $\Y$ and is $\cC^\infty$;
\item $\mathcal{R}_\epsilon=\F_\epsilon(\mathbf{0},\lambda_*)\longrightarrow 0$ as $\epsilon \longrightarrow 0$;
\item $\cL_\epsilon=D\F_\epsilon(\mathbf{0},\lambda_*)$ can be decomposed in two parts:
\bqq
\label{eq:DecompositonL}
\mathcal{L}_\epsilon=\mathcal{L}_\epsilon^i+\mathcal{L}_\epsilon^p,
\eqq
where $\mathcal{L}_\epsilon^i$ is invertible with bounded inverse on suitable Banach spaces and $\mathcal{L}_\epsilon^p$ is an $\epsilon$-perturbation: $\mathcal{L}_\epsilon^p\longrightarrow0$ as $\epsilon \longrightarrow 0$.
\end{enumerate}
Then, to conclude the proof of Theorem \ref{thm:existence}, we will use a fixed point iteration argument on the map $\F_\epsilon$ which will give the existence of $\left(\bW(\epsilon),\lambda(\epsilon) \right)$, solution of \eqref{eq:Compact}, in a neighborhood of $(\mathbf{0},\lambda_*)$ for small values of $\epsilon>0$. 

\section{Proof of Theorem \ref{thm:existence}}\label{proof}

\subsection{Estimates for the error terms $\mathcal{R}_\epsilon$}\label{subsec:R}

In this section, we provide estimates for the error terms $\mathcal{R}_\epsilon$ as stated in the  following result.
\begin{prop}\label{prop:EstR}
With $\cR_\epsilon$ defined in equations \eqref{eq:Compact} and  \eqref{eq:DefRq}, \eqref{eq:DefRb}, \eqref{eq:DefRe}, \eqref{eq:DefRf} below, we have 
\bqq
\label{eq:LimitR}
\underset{\epsilon \rightarrow 0}{\lim}\;\left\| \mathcal{R}_\epsilon \right\|_\Y=0,
\eqq
where $\Y=\left(L^2_\eta \times L^2_\eta \right)^5$.
\end{prop}
We first make the various error terms that enter $\cR_\epsilon$ explicit: we find, setting  $\left(\bW,\lambda\right)=(\mathbf{0},\lambda_*)$, in systems \eqref{eq:Wq}, \eqref{eq:Wb}, \eqref{eq:Wbe}, \eqref{eq:Wef}, and \eqref{eq:Wf},   for the quiescent part:
\begin{subequations}
\label{eq:DefRq}
\begin{align}
\cR_\epsilon^{u,q}(\xi)=&\; c_*\bu_q\left(\epsilon\left(\xi\right)\right)\chi_q'(\xi+\xi_q)+c_*\bu_b(\xi-\xi_{b}+\xi_q)\chi_b'(\xi+\xi_q)\mathds{1}_{\xi_b}(\xi+\xi_q) +\bF_q(\xi+\xi_q)\mathds{1}_{\xi_q}(\xi+\xi_q)\nonumber\\
&+\bC_q(\xi+\xi_q),\label{eq:Ruq} \\
\cR_\epsilon^{v,q}(\xi)=&\;c_*\bv_q(\epsilon(\xi))\chi_q'(\xi+\xi_q)+c_*v_*\chi'_b(\xi+\xi_q)\mathds{1}_{\xi_b}(\xi+\xi_q)\nonumber\\
&+\epsilon\chi_q(\xi+\xi_q)\left(\bu_q(\epsilon \xi)-\gamma\bv_q(\epsilon \xi)\right)\left(1-\Theta(\bv_q(\epsilon \xi))\right)\nonumber\\
&+\epsilon \left(\int_\R (\bu_b(\xi)-\gamma v_*)\chi_b(\xi+\xi_b)d\xi\right)\psi\left(\xi-\frac{\xi_b-\xi_q}{2}\right);\label{eq:Rvq}
\end{align}
\end{subequations}
for the back part:
\begin{subequations}
\label{eq:DefRb}
\begin{align}
\cR_\epsilon^{u,b}(\xi)=&\; 0,\label{eq:Rub}\\
\cR_\epsilon^{v,b}(\xi)=&\; \epsilon \left(\bu_b(\xi)-\gamma v_*\right)\chi_b(\xi+\xi_b)-\epsilon \left(\int_\R  \left(\bu_b(\xi)-\gamma v_*\right)\chi_b(\xi+\xi_b)d\xi\right)\psi\left(\xi+\frac{\xi_b-\xi_q}{2}\right);\label{eq:Rvb}
\end{align}
\end{subequations}
for the excitatory parts:
\begin{subequations}
\label{eq:DefRe}
\begin{align}
\cR_\epsilon^{u,be}(\xi)=&\; c_*\bu_b(\xi-\xi_{b}+\xi_{be})\chi_b'(\xi+\xi_{be})\mathds{1}_{\xi_{be}}(\xi+\xi_{be})+c_*\bu_e(\epsilon\xi+\epsilon\xi_{be})\chi_e'(\xi+\xi_{be})\mathds{1}_{\xi_{be}}(\xi+\xi_{be})\nonumber\\
&+\bF_{be}(\xi+\xi_{be})\mathds{1}_{\xi_{be}}(\xi+\xi_{be})+\bC_{be}(\xi+\xi_{be}), \label{eq:Rube}\\
\cR_\epsilon^{v,be}(\xi)=&\; c_*v_b\chi'_b(\xi+\xi_{be})\mathds{1}_{\xi_{be}}(\xi+\xi_{be})+c_*\bv_e(\epsilon\xi+\epsilon\xi_{be})\chi_e'(\xi+\xi_{be})\mathds{1}_{\xi_{be}}(\xi+\xi_{be}),\nonumber\\
&+\epsilon\chi_e(\xi+\xi_{be})\left(\bu_e(\epsilon \xi+\epsilon\xi_{be})-\gamma\bv_e(\epsilon \xi+\epsilon\xi_{be})\right)\left(1-\Theta(\bv_e(\epsilon \xi+\epsilon\xi_{be}))\right)\mathds{1}_{\xi_{be}}(\xi+\xi_{be})\label{eq:Rvbe}\\
\cR_\epsilon^{u,ef}(\xi)=&\; c_*\bu_e(\epsilon\xi)\chi_e'(\xi)\mathds{1}_{\xi_{ef}}(\xi)+c_*\bu_f(\xi-\xi_f)\chi_f'(\xi)+\bF_{ef}(\xi)\mathds{1}_{\xi_{ef}}(\xi)+\bC_{ef}(\xi),\label{eq:Ruef}\\
\cR_\epsilon^{v,ef}(\xi)=&\; c_*\bv_e(\epsilon\xi)\chi_e'(\xi)\mathds{1}_{\xi_{ef}}(\xi)+\epsilon\chi_e(\xi)\left(\bu_e(\epsilon \xi)-\gamma\bv_e(\epsilon \xi)\right)\left(1-\Theta(\bv_e(\epsilon \xi))\right)\mathds{1}_{ef}(\xi)\nonumber\\
&+\epsilon \left(\int_\R  \bu_f(\xi)\chi_f(\xi+\xi_f)d\xi\right)\psi\left(\xi-\frac{\xi_f}{2}\right);\label{eq:Rvef}
\end{align}
\end{subequations}
and for the front part:
\begin{subequations}
\label{eq:DefRf}
\begin{align}
\cR_\epsilon^{u,f}(\xi)=&\; 0,\label{eq:Ruf}\\
\cR_\epsilon^{v,f}(\xi)=&\; \epsilon \bu_f(\xi)\chi_f(\xi+\xi_f)-\epsilon \left(\int_\R  \bu_f(\xi)\chi_f(\xi+\xi_f)d\xi\right)\psi\left(\xi+\frac{\xi_f}{2}\right).\label{eq:Rvf}
\end{align}
\end{subequations}

We see that in the definition of $\cR_\epsilon$, through equations \eqref{eq:DefRq}, \eqref{eq:DefRb}, \eqref{eq:DefRe}, and \eqref{eq:DefRf},  terms of the same "nature" appear several times. Indeed, we can find terms that only involve derivatives of the partition functions, $\chi'_j$ for $j\in J_\bw$. We can also find commutator terms of two types. The first type, denoted $\bF_q$, $\bF_{be}$, and $\bF_{ef}$, represents the commutators between  nonlinearity $f$ and partition of unity defined in equations \eqref{eq:ComF}, \eqref{eq:ComDefF}, and \eqref{eq:ComFj}. The second type, denoted $\bC_q$, $\bC_{be}$, and $\bC_{ef}$,  involves the commutators between  convolution and partition of unity as defined in equations \eqref{eq:ComC}, \eqref{eq:ComDefC}, and \eqref{eq:ComCj}. There are also error terms that stem from the fact that, in our Ansatz, the back $\bU_b$ and the front $\bU_f$ are only solutions at $\epsilon=0$, and from the fact that $\bU_q$ and $\bU_e$ are solutions of the modified equation \eqref{eq:SlowMod}, only. Finally, there are terms 
that involve the function $\psi$ and we have already seen in equations \eqref{eq:EstB} and \eqref{eq:EstF} that
\begin{align*}
\left\| \epsilon \left(\int_\R \left(\bu_b(\xi)-\gamma v_b\right)\chi_b(\xi+\xi_b) d\xi\right) \psi\left(\cdot \pm  \frac{\xi_b-\xi_q}{2}\right) \right\| _{L^2_\eta}&= \cO\left( \epsilon^{1-\frac{ \eta \eta_b}{2} }|\ln \epsilon | \right),\\
\left\| \epsilon \left(\int_\R\bu_f(\xi)\chi_f(\xi+\xi_f)d\xi\right) \psi\left(\cdot \pm \frac{\xi_f}{2}\right) \right\| _{L^2_\eta}&= \cO\left( \epsilon^{1-\frac{ \eta\eta_f}{2} }|\ln \epsilon | \right),
\end{align*}
as $\epsilon \rightarrow0$. 

We have divided the proof of Proposition \ref{prop:EstR} into three Propositions \ref{prop:EstChi}, \ref{prop:EstCom}, and \ref{prop:EstErr}, where we respectively provide estimates for the $\chi_j'$-terms, the commutator terms and the error terms.

\begin{prop}[Estimates for $\chi_j'$-terms]\label{prop:EstChi} The following estimates hold as $\epsilon \rightarrow0$:
\begin{itemize}
\item $\left\| c_*\bu_q\left(\epsilon~\cdot\right)\chi_q'(\cdot+\xi_q)+c_*\bu_b(\cdot-\xi_{b}+\xi_q)\chi_b'(\cdot+\xi_q)\mathds{1}_{\xi_b}(\cdot+\xi_q) \right\|_{L^2_\eta}=\cO(\epsilon^{\min(1,\eta_*\eta_b)})$;
\item $\left\| c_*\bv_q(\epsilon~\cdot)\chi_q'(\cdot+\xi_q)+c_*v_*\chi'_b(\cdot+\xi_q)\mathds{1}_{\xi_b}(\cdot+\xi_q) \right\|_{L^2_\eta}=\cO(\epsilon)$; 
\item $\left\| c_*\bu_b(\cdot-\xi_{b}+\xi_{be})\chi_b'(\cdot+\xi_{be})\mathds{1}_{\xi_{be}}(\cdot+\xi_{be})+c_*\bu_e(\epsilon\cdot+\epsilon\xi_{be})\chi_e'(\cdot+\xi_{be})\mathds{1}_{\xi_{be}}(\cdot+\xi_{be}) \right\|_{L^2_\eta}=\cO(\epsilon^{\min(1,\eta_*\eta_b)})$; 
\item $\left\| c_*v_b\chi'_b(\cdot+\xi_{be})\mathds{1}_{\xi_{be}}(\cdot+\xi_{be})+c_*\bv_e(\epsilon\cdot+\epsilon\xi_{be})\chi_e'(\cdot+\xi_{be})\mathds{1}_{\xi_{be}}(\cdot+\xi_{be})\right\|_{L^2_\eta}=\cO(\epsilon)$; 
\item $\left\| c_*\bu_e(\epsilon~\cdot)\chi_e'\mathds{1}_{\xi_{ef}}+c_*\bu_f(\cdot-\xi_f)\chi_f'\right\|_{L^2_\eta}=\cO(\epsilon^{\min(1,\eta_*\eta_f)})$; 
\item $\left\| c_*\bv_e(\epsilon~\cdot)\chi_e'\mathds{1}_{\xi_{ef}} \right\|_{L^2_\eta}=\cO(\epsilon)$.
\end{itemize}
\end{prop}
\begin{Proof}
We only prove the last two estimates as the others can easily be deduced following similar types of argument. As for all $|\xi|\geq 1$, $\chi'_e(\xi)\mathds{1}_{\xi_{ef}}(x)=\chi'_f(\xi)=0$, we have that $\chi_e'\mathds{1}_{\xi_{ef}} \in L^2_\eta$ and $\chi_f'\in L^2_\eta$ for any $\eta>0$. For all $|\xi| \leq 1$, the following asymptotic estimates hold
\begin{align*}
\bu_f(\xi-\xi_f)-1&=\cO\left( \epsilon^{\eta_*\eta_f}\right), \\
\bu_e(\epsilon \xi)-1&= \cO(\epsilon ),\\
\bv_e(\epsilon \xi)&=\cO(\epsilon),
\end{align*}
uniformly in $\xi$ as $\epsilon\rightarrow0$. Noticing that 
\bqs
\chi_e'(\xi)\mathds{1}_{\xi_{ef}}(\xi)=-\chi'_f(\xi), \quad \forall \xi \in \R,
\eqs
we then obtain
\bqs
\left\| c_*\bu_e(\epsilon~\cdot)\chi_e'\mathds{1}_{\xi_{ef}}+c_*\bu_f(\cdot-\xi_f)\chi_f'\right\|_{L^2_\eta}\lesssim \epsilon \left\|\chi_e'\mathds{1}_{\xi_{ef}}\right\|_{L^2_\eta}+\epsilon^{\eta_*\eta_f}\left\|\chi_f'\right\|_{L^2_\eta},
\eqs
and
\bqs
\left\| c_*\bv_e(\epsilon~\cdot)\chi_e'\mathds{1}_{\xi_{ef}} \right\|_{L^2_\eta} \lesssim \epsilon \left\|\chi_e'\mathds{1}_{\xi_{ef}}\right\|_{L^2_\eta}.
\eqs
This gives the desired estimates.
\end{Proof}

\begin{prop}[Estimates for commutator terms]\label{prop:EstCom} The following estimates hold as $\epsilon \rightarrow0$:
\begin{itemize}
\item $\left\| \bF_q(\cdot+\xi_q)\mathds{1}_{\xi_q}(\cdot+\xi_q) \right\|_{L^2_\eta}=\cO(\epsilon^{2\min(1,\eta_*\eta_b)})$;
\item $\left\| \bF_{be}(\cdot+\xi_{be})\mathds{1}_{\xi_{be}}(\cdot+\xi_{be}) \right\|_{L^2_\eta}=\cO(\epsilon^{2\min(1,\eta_*\eta_b)})$;
\item $\left\| \bF_{ef} \mathds{1}_{\xi_{ef}} \right\|_{L^2_\eta}=\cO(\epsilon^{2\min(1,\eta_*\eta_f)})$;
\item $\left\| \bC_q(\cdot+\xi_q) \right\|_{L^2_\eta}=\cO(\epsilon^{\min(1,\eta_*\eta_b,\eta_0\eta_b)})$;
\item $\left\| \bC_{be}(\cdot+\xi_{be}) \right\|_{L^2_\eta}=\cO(\epsilon^{\min(1,\eta_*\eta_b,\eta_0\eta_b)})$;
\item $\left\| \bC_{ef} \right\|_{L^2_\eta}=\cO(\epsilon^{\min(1,\eta_*\eta_f,\eta_0\eta_f)})$.
\end{itemize}
\end{prop}
\begin{Proof} 
We prove the third and last estimates, the others being easily deduced from them. First, we recall the definition of $\bF_{ef}$:
\bqs
\bF_{ef}(\xi)=f\left( \bu_e(\epsilon\xi)\chi_e(\xi)+\bu_f(\xi-\xi_f)\chi_f(\xi)\right)-f(\bu_e(\epsilon\xi))\chi_e(\xi)-f(\bu_f(\xi-\xi_f))\chi_f(\xi),
\eqs
for all $\xi\in \R$. For all $|\xi|\leq 1$, we write $\bu_{ef}(\xi)=(\bu_e(\epsilon\xi)-1)\chi_e(\xi)+(\bu_f(\xi-\xi_f)-1)\chi_f(\xi)$ and we have
\begin{align*}
f\left(\bu_e(\epsilon\xi)\chi_e(\xi)+\bu_f(\xi-\xi_f)\chi_f(\xi)\right)&=f\left(1+\left(\bu_e(\epsilon\xi)-1\right)\chi_e(\xi)+\left(\bu_f(\xi-\xi_f)-1\right)\chi_f(\xi)\right)\\
&=f'(1)\bu_{ef}(\xi) + \left(\bu_{ef}(\xi) \right)^2\int_0^1f''(1+\tau \bu_{ef}(\xi))(1-\tau)d\tau,\\
f(\bu_e(\epsilon\xi))&=f'(1)(\bu_e(\epsilon\xi)-1)+(\bu_e(\epsilon\xi)-1)^2\int_0^1 f''(1+\tau (\bu_e(\epsilon\xi)-1))(1-\tau)d\tau  ,\\
f(\bu_f(\xi-\xi_f))&=f'(1)(\bu_f(\xi-\xi_f)-1)\\
&\quad +(\bu_f(\xi-\xi_f)-1)^2\int_0^1 f''(1+\tau (\bu_f(\xi-\xi_f)-1))(1-\tau)d\tau,\\
\end{align*}
as $\epsilon\rightarrow 0$. We see that we only get corrections at quadratic order and thus
\bqs
\left\| \bF_{ef} \mathds{1}_{\xi_{ef}} \right\|_{L^2_\eta} \lesssim \left(\epsilon^2+\epsilon^{2\eta_*\eta_f}\right)\left\| \mathds{1}_{\xi_{ef}} \right\|_{L^2_\eta}.
\eqs

For the last estimate, we note that, by assumption, there exists $\eta_0>0$ such that $\left\| \K \right\|_{L^1_{\eta_0}}<\infty$. For all $(\xi,\zeta)\in\R^2$, the following estimates holds
\bqs
\bG(\xi,\zeta):=e^{-\eta_0|\xi-\zeta|}\left| \bu_e(\epsilon\zeta)-\bu_f(\zeta-\xi_f) \right|\left| \chi_f(\zeta)-\chi_f(\xi) \right| \lesssim e^{-\eta_0|\xi|} \epsilon^{\min(1,\eta_*\eta_f,\eta_0\eta_f)}.
\eqs
This can be seen by evaluating $\bG$ in different regions of the plane:
\begin{itemize}
\item  for $\xi \geq 1 $ and $\zeta \geq 1$ we have $\bG(\xi,\zeta)=0$;
\item  for $\xi \leq -1 $ and $\zeta \leq -1$ we have $\bG(\xi,\zeta)=0$;
\item for $\xi \geq 1$ and $\zeta \leq 1$ we have
\bqs
\bG(\xi,\zeta) \leq e^{-\eta_0\xi}~\underset{\zeta \leq 1}{\sup}\left(e^{\eta_0\zeta} \left| \bu_e(\epsilon\zeta)-\bu_f(\zeta-\xi_f) \right|\right) \lesssim e^{-\eta_0\xi} \left( \epsilon +\epsilon^{\eta_*\eta_f} \right);
\eqs 
\item for $\xi \leq -1$ and $\zeta \geq -1$ we have
\bqs
\bG(\xi,\zeta) \leq e^{\eta_0\xi}~\underset{\zeta \geq -1}{\sup}\left(e^{-\eta_0\zeta} \left| \bu_e(\epsilon\zeta)-\bu_f(\zeta-\xi_f) \right|\right) \lesssim e^{\eta_0\xi} \epsilon^{\min(1,\eta_*\eta_f,\eta_0\eta_f)};
\eqs
\item for $|\xi|\leq 1$ and $|\zeta|\leq 1$, we have
\bqs
\bG(\xi,\zeta) \lesssim e^{-\eta_0|\xi|} \left( \epsilon +\epsilon^{\eta_*\eta_f} \right);
\eqs
\item for $|\xi|\leq 1$ and $\zeta\leq -1$, we have
\bqs
\bG(\xi,\zeta) \lesssim e^{-\eta_0|\xi|} \left( \epsilon +\epsilon^{\eta_*\eta_f} \right);
\eqs
\item for $|\xi|\leq 1$ and $\zeta\geq 1$, we have
\bqs
\bG(\xi,\zeta) \lesssim e^{-\eta_0|\xi|} \epsilon^{\min(1,\eta_*\eta_f,\eta_0\eta_f)}.
\eqs
\end{itemize}
Finally, if $\eta<\eta_0$, we obtain
\bqs
\left\|  \bC_{ef} \right\|_{ L^2_\eta} \lesssim \epsilon^{\min(1,\eta_*\eta_f,\eta_0\eta_f)} \left\| \K \right\|_{L^1_{\eta_0}(\R)}.
\eqs

\end{Proof}

\begin{prop}[Estimates for the Ansatz error terms]\label{prop:EstErr} The following estimates hold as $\epsilon \rightarrow0$:
\begin{itemize}
\item $\left\| \epsilon\chi_q(\cdot+\xi_q)\left(\bu_q(\epsilon ~\cdot)-\gamma\bv_q(\epsilon ~\cdot)\right)\left(1-\Theta(\bv_q(\epsilon~\cdot))\right) \right\|_{L^2_\eta}=\cO(\epsilon)$;
\item $\left\| \epsilon \left(\bu_b-\gamma v_*\right)\chi_b(\cdot+\xi_b) \right\|_{L^2_\eta}=\cO\left( \epsilon^{1- \eta \eta_b}\right)$;
\item $\left\| \epsilon\chi_e(\cdot+\xi_{be})\left(\bu_e(\epsilon \cdot+\epsilon\xi_{be})-\gamma\bv_e(\epsilon \cdot+\epsilon\xi_{be})\right)\left(1-\Theta(\bv_e(\epsilon \cdot+\epsilon\xi_{be}))\right)\mathds{1}_{\xi_{be}}(\cdot+\xi_{be})\right\|_{L^2_\eta}=\cO(\epsilon)$;
\item $\left\| \epsilon\chi_e\left(\bu_e(\epsilon ~\cdot)-\gamma\bv_e(\epsilon ~\cdot)\right)\left(1-\Theta(\bv_e(\epsilon ~\cdot))\right)\mathds{1}_{ef} \right\|_{L^2_\eta}=\cO(\epsilon)$;
\item $\left\| \epsilon \bu_f\chi_f(\cdot+\xi_f) \right\|_{L^2_\eta}=\cO\left( \epsilon^{1- \eta \eta_f }\right) $.
\end{itemize}
\end{prop}

\begin{Proof} 
Once again, we only prove the last two estimates. First, we see that $\chi_e \left(1-\Theta(\bv_e(\epsilon ~\cdot))\right)\mathds{1}_{ef}$ has a compact support in $[-1,1]$ so that $ \chi_e \left(1-\Theta(\bv_e(\epsilon ~\cdot))\right)\mathds{1}_{ef} \in L^2_\eta$ for all $\eta>0$. And we directly obtain
\bqs
\left\| \epsilon\chi_e\left(\bu_e(\epsilon ~\cdot)-\gamma\bv_e(\epsilon ~\cdot)\right)\left(1-\Theta(\bv_e(\epsilon ~\cdot))\right)\mathds{1}_{ef} \right\|_{L^2_\eta} \lesssim \epsilon  \left\| \chi_e \left(1-\Theta(\bv_e(\epsilon ~\cdot))\right)\mathds{1}_{ef}\right\|_{L^2_\eta}.
\eqs
Second, we use the definition of the $L^2_\eta$ norm of $\bu_f \chi_f(\cdot+\xi_f)$ and the property of $\chi_f$ to obtain
\begin{align*}
\int_\R e^{2\eta |\xi|}\left(\bu_f(\xi) \chi_f(\xi+\xi_f)\right)^2d\xi&= \int_{\R} e^{2\eta |\xi-\xi_f|}\left(\bu_f(\xi-\xi_f) \chi_f(\xi)\right)^2d\xi \\
&=\int_{-1}^1 e^{2\eta |\xi-\xi_f|}\left(\bu_f(\xi-\xi_f) \chi_f(\xi)\right)^2d\xi  + \int_1^\infty e^{2\eta |\xi-\xi_f|}\left(\bu_f(\xi-\xi_f)\right)^2d\xi  \\
&=\int_{-1}^1 e^{2\eta |\xi-\xi_f|}\left(\bu_f(\xi-\xi_f) \chi_f(\xi)\right)^2d\xi  + \int_{1+\xi_f}^\infty e^{2\eta |\xi|}\left(\bu_f(\xi)\right)^2d\xi  \\
&=\I_1(\epsilon)+\I_2(\epsilon).
\end{align*}
A straightforward computation gives
\bqs
\I_1(\epsilon)=\cO(\epsilon^{-2\eta\eta_f}) \text{ and } \I_2(\epsilon)=\cO(\epsilon^{2(\eta_*-\eta)\eta_f}) \text{ as } \epsilon\rightarrow 0,
\eqs
which completes the proof as $\mu < \mu_*$ by definition.
\end{Proof}

\begin{Proof}[ of Proposition \ref{prop:EstR}.] We can now combine the estimates of Propositions \ref{prop:EstChi}, \ref{prop:EstCom} and \ref{prop:EstErr} to obtain the limit \eqref{eq:LimitR} which concludes the proof.
\end{Proof}

\subsection{Study of the linear part $\cL_\epsilon^i$}\label{subsec:Li}

In this section, we shall prove that the linear operator $\cL_\epsilon^i$ is invertible with bounded inverse on a suitable Banach space. We define the linear operator $\cL_\epsilon^i$ as follows
\bqs
\cL_\epsilon^i: \begin{array}{ccl}
\X \times \R^5  &\longrightarrow& \Y,
 \end{array}
\eqs
where $\cL_\epsilon^i$ can be written in matrix form as
\bqq
\label{eq:LiEps}
\cL_\epsilon^i=\left( \A_\bW^i(\epsilon) | \A_\lambda ^i(\epsilon) \right) 
\eqq

\bqs
\A_\bW^i(\epsilon)=\left( 
\begin{matrix}
\cL(\bu_q)  & -1 & 0  & 0 & 0  & 0 & 0 & 0 & 0 & 0\\
0  & c_*\frac{d}{d\xi} & 0 & 0 & 0  & 0 & 0  & 0 & 0  & 0\\
0  & 0 & \cL(\bu_b)  & -1 & 0  & 0 & 0  & 0 & 0  & 0\\
0  & 0 & 0  & c_*\frac{d}{d\xi} & 0  & 0 & 0  & 0 & 0  & 0\\
0  & 0 & 0  & 0 & \cL(\tau_{\xi_{be}}\cdot\bu_e)  & -1 & 0 & 0 & 0 & 0\\
0  & 0 & 0  & 0 & 0  & c_*\frac{d}{d\xi} & 0  & 0 & 0  & 0\\
0  & 0 & 0  & 0 & 0  & 0 & \cL(\bu_e)  & -1 & 0  & 0\\
0  & 0 & 0  & 0 & 0  & 0 & 0  & c_*\frac{d}{d\xi} & 0  & 0\\
0  & 0 & 0  & 0 & 0  & 0 & 0  & 0 & \cL(\bu_f)  & -1\\
0  & 0 & 0  & 0 & 0  & 0 & 0  & 0 & 0  & c_*\frac{d}{d\xi}
\end{matrix}
\right),
\eqs
where we have defined
\begin{subequations}
\begin{align}
\label{eq:Lquiescent}
\cL(\bu_q)&: \left\{\begin{array}{ccl}
H^1_\eta &\longrightarrow& L^2_\eta \\
 \bw&\longmapsto& -\bw+\K \ast \bw +  c_* \frac{d}{d\xi}\bw +f'\left(\bu_q(\epsilon\xi)\right)\bw,
 \end{array}\right. \\
\label{eq:Lback}
\cL(\bu_b)&: \left\{\begin{array}{ccl}
H^1_\eta &\longrightarrow& L^2_\eta \\
 \bw&\longmapsto& -\bw+\K \ast \bw +  c_* \frac{d}{d\xi}\bw +f'\left(\bu_b(\xi)\right)\bw,
 \end{array}\right.\\
\label{eq:Lexcitatory}
\cL(\bu_e)&: \left\{\begin{array}{ccl}
H^1_\eta &\longrightarrow& L^2_\eta \\
 \bw&\longmapsto& -\bw+\K \ast \bw +  c_* \frac{d}{d\xi}\bw +f'\left(\bu_e(\epsilon\xi)\right)\bw,
 \end{array}\right.\\
\label{eq:Lfront}
\cL(\bu_f)&: \left\{\begin{array}{ccl}
H^1_\eta &\longrightarrow& L^2_\eta \\
 \bw&\longmapsto& -\bw+\K \ast \bw +  c_* \frac{d}{d\xi}\bw +f'\left(\bu_f(\xi)\right)\bw,
 \end{array}\right.
\end{align}
\end{subequations}
Finally, the matrix operator $\A_\lambda^i(\epsilon)$ has the following form
\bqq
\label{eq:AiLambda}
\A_\lambda^i(\epsilon)= \left( \bA^{i,1}_\lambda(\epsilon) ~ \bA_\lambda^{i,2}(\epsilon) ~ \bA_\lambda^{i,3}(\epsilon) ~ \bA_\lambda^{i,4}(\epsilon) ~ \bA_\lambda^{i,5}(\epsilon) \right)
\eqq
where the columns $\bA^{i,j}_\lambda(\epsilon)$ are defined as
\bqs
 \bA^{i,1}_\lambda(\epsilon)=\left( 
\begin{matrix}
0\\
0\\
\bu_b'\chi_b(\cdot+\xi_b)\\
0\\
0\\
0\\
0\\
0\\
\bu_f'\chi_f(\cdot+\xi_f)\\
0
\end{matrix}
\right),\quad
 \bA^{i,2}_\lambda(\epsilon)=\left( 
\begin{matrix}
c_*\partial_{v_q}\left( \bu_q(\epsilon\cdot) \right)\chi'_q(\cdot+\xi_q)\\
c_*\partial_{v_q}\left( \bv_q(\epsilon\cdot) \right)\chi'_q(\cdot+\xi_q) \\
0\\
0\\
0\\
0\\
0\\
0\\
0\\
0
\end{matrix}
\right),
\eqs
\bqs
 \bA^{i,3}_\lambda(\epsilon)=\left( 
\begin{matrix}
c_* \partial_{v_b}\left(\tau_{\xi_{q}}\cdot \left[\bu_b(\cdot-\xi_b))\chi'_b\mathds{1}_{\xi_{b}}\right]\right)\\
c_*\tau_{\xi_b}\cdot\left[ \chi'_b\mathds{1}_{\xi_q} \right] \\
-c'_b\bu_b'\chi_b(\cdot+\xi_b\\
0\\
c_* \partial_{v_b}\left(\tau_{\xi_{be}}\cdot \left[\bu_b(\cdot-\xi_b))\chi'_b\mathds{1}_{\xi_{be}}\right]\right)\\
c_* \tau_{\xi_{be}}\cdot \left[\chi'_b\mathds{1}_{\xi_{be}}\right]\\
0\\
0\\
0\\
0
\end{matrix}
\right),\quad
 \bA^{i,4}_\lambda(\epsilon)=\left( 
\begin{matrix}
0\\
0\\
0\\
0\\
c_* \partial_{v_{be}}\left(\tau_{\xi_{be}}\cdot \left[\bu_e(\epsilon\cdot))\chi'_e\mathds{1}_{\xi_{be}}\right]\right)\\
c_* \partial_{v_{be}}\left(\tau_{\xi_{be}}\cdot \left[\bv_e(\epsilon\cdot))\chi'_e\mathds{1}_{\xi_{be}}\right]\right)\\
c_* \partial_{v_{be}}(\bu_e(\epsilon\cdot))\chi'_e\mathds{1}_{\xi_{be}}\\
c_* \partial_{v_{be}}(\bv_e(\epsilon\cdot))\chi'_e\mathds{1}_{\xi_{be}}\\
0\\
0
\end{matrix}
\right),
\eqs
and
\bqs
 \bA^{i,5}_\lambda(\epsilon)=\left( 
\begin{matrix}
0\\
0\\
0\\
0\\
c_* \partial_{v_{ef}}\left(\tau_{\xi_{be}}\cdot \left[\bu_e(\epsilon\cdot))\chi'_e\mathds{1}_{\xi_{ef}}\right]\right)\\
c_* \partial_{v_{ef}}\left(\tau_{\xi_{be}}\cdot \left[\bv_e(\epsilon\cdot))\chi'_e\mathds{1}_{\xi_{ef}}\right]\right)\\
c_* \partial_{v_{ef}}(\bu_e(\epsilon\cdot))\chi'_e\mathds{1}_{\xi_{ef}}\\
c_* \partial_{v_{ef}}(\bv_e(\epsilon\cdot))\chi'_e\mathds{1}_{\xi_{ef}}\\
0\\
0
\end{matrix}
\right),
\eqs
where $c'_b=\partial_{v_b}c_b$. In order to clearly see that the operator $\cL_\epsilon^i$ is invertible, we will rewrite it in a different basis so that the new operator expressed in that basis has a triangular form, and each entry on the diagonal is invertible. More precisely, by permuting the columns and the rows, we have that
\bqq
\label{eq:matrixLiEpsEq}
\cL_\epsilon^i \sim \left(
\begin{matrix}
\cL_f^i & \mathbf{0}_{2,3} & \mathbf{0}_{2,3} & \mathbf{0}_{2,6} \\
\cL_{b,f}^i & \cL_b^i & \mathbf{0}_{2,3} & \mathbf{0}_{2,6} \\
\mathbf{0}_{2,3} & \cL_{q,b}^i & \cL_q^i & \mathbf{0}_{2,6} \\
\mathbf{0}_{4,3} & \cL_{e,b}^i & \mathbf{0}_{4,3} & \cL_e^i
\end{matrix}
 \right),
\eqq
where the diagonal operators appearing in the above matrix \eqref{eq:matrixLiEpsEq} are
\begin{subequations}
\begin{align}
\cL_f^i &=\left(
\begin{matrix}
\cL(\bu_f)  & -1 & \bu_f'\chi_f(\cdot+\xi_f)\\
 0  & c_*\frac{d}{d\xi} & 0
\end{matrix}
\right), \quad \longleftrightarrow \quad \underline{\bW}_f = \left(\bw_f^u,\bw_f^v ,c\right), \label{eq:LiFront}\\
\cL_b^i &=\left(
\begin{matrix}
\cL(\bu_b)  & -1 & -c'_b\bu_b'\chi_b(\cdot+\xi_b)\\
 0  & c_*\frac{d}{d\xi} & 0
\end{matrix}
\right), \quad \longleftrightarrow \quad \underline{\bW}_b = \left(\bw_b^u,\bw_b^v ,v_b\right), \label{eq:LiBack}\\
\cL_q^i &=\left(
\begin{matrix}
\cL(\bu_q)  & -1 & c_*\partial_{v_q}\left( \bu_q(\epsilon\cdot) \right)\chi'_q(\cdot+\xi_q) \\
 0  & c_*\frac{d}{d\xi} & c_*\partial_{v_q}\left( \bv_q(\epsilon\cdot) \right)\chi'_q(\cdot+\xi_q)
\end{matrix}
\right), \quad \longleftrightarrow \quad \underline{\bW}_q = \left(\bw_q^u,\bw_q^v ,v_q\right), \label{eq:LiQuiescent}\\
\cL_e^i &=\left(
\begin{matrix}
\cL(\tau_{\xi_{be}}\bu_e)  & -1 & 0 & 0 & c_* \partial_{v_{be}}\left(\tau_{\xi_{be}}\cdot \left[\bu_e(\epsilon\cdot))\chi'_e\mathds{1}_{\xi_{be}}\right]\right) & c_* \partial_{v_{ef}}\left(\tau_{\xi_{be}}\cdot \left[\bu_e(\epsilon\cdot))\chi'_e\mathds{1}_{\xi_{ef}}\right]\right) \\
 0  & c_*\frac{d}{d\xi} & 0 & 0 & c_* \partial_{v_{be}}\left(\tau_{\xi_{be}}\cdot \left[\bv_e(\epsilon\cdot))\chi'_e\mathds{1}_{\xi_{be}}\right]\right) & c_* \partial_{v_{ef}}\left(\tau_{\xi_{be}}\cdot \left[\bv_e(\epsilon\cdot))\chi'_e\mathds{1}_{\xi_{ef}}\right]\right) \\
 0 & 0 & \cL(\bu_e) & -1 & c_* \partial_{v_{be}}(\bu_e(\epsilon\cdot))\chi'_e\mathds{1}_{\xi_{be}} & c_* \partial_{v_{ef}}(\bu_e(\epsilon\cdot))\chi'_e\mathds{1}_{\xi_{ef}}\\
 0 & 0 &  0  & c_*\frac{d}{d\xi} & c_* \partial_{v_{be}}(\bv_e(\epsilon\cdot))\chi'_e\mathds{1}_{\xi_{be}} & c_* \partial_{v_{ef}}(\bv_e(\epsilon\cdot))\chi'_e\mathds{1}_{\xi_{ef}}
 \end{matrix}
\right), \label{eq:Liexcitatory}
\end{align}
\end{subequations}
where the last matrix operator is expressed in the coordinates  $\underline{\bW}_e = \left(\bw_{be}^u,\bw_{be}^v,\bw_{ef}^u,\bw_{ef}^v ,v_{be},v_{ef}\right)$. The remaining three off-diagonal operators are thus
\bqs
\cL_{b,f}^i =\left(
\begin{matrix}
0 & 0 & \bu_b'\chi_b(\cdot+\xi_b)\\
0 & 0 & 0
\end{matrix}
\right), \quad \cL_{q,b}^i =\left(
\begin{matrix}
0 & 0 & c_* \partial_{v_b}\left(\tau_{\xi_{q}}\cdot \left[\bu_b(\cdot-\xi_b))\chi'_b\mathds{1}_{\xi_{b}}\right]\right)\\
0 & 0 & c_*\tau_{\xi_b}\cdot\left[ \chi'_b\mathds{1}_{\xi_q} \right] 
\end{matrix}
\right),
\eqs
and
\bqs
\cL_{e,b}^i =\left(
\begin{matrix}
0 & 0 & c_* \partial_{v_b}\left(\tau_{\xi_{be}}\cdot \left[\bu_b(\cdot-\xi_b))\chi'_b\mathds{1}_{\xi_{be}}\right]\right)\\
0 & 0 & c_* \tau_{\xi_{be}}\cdot \left[\chi'_b\mathds{1}_{\xi_{be}}\right] 
\end{matrix}
\right).
\eqs

We would like to show that each of the operators appearing on the diagonal is bounded invertible on a suitable Banach space. We treat each case separately in the following sections.

\subsubsection{Invertibility of $\cL_f^i$ and $\cL_b^i$}

In this section, we will show that both $\cL_f^i$ and $\cL_b^i$ are bounded invertible. We therefore introduce the following spaces:
\begin{subequations}
\label{eq:BanachSpace}
\begin{align}
\X_f&:= \left(H^1_\eta \cap \left\{ \bu:\R\rightarrow\R~|~ \langle \bu, \bu'_f\rangle = 0 \right\} \right)\times  H^1_\eta \times\R,\\
\X_b&:= \left(H^1_\eta \cap \left\{ \bu:\R\rightarrow\R~|~ \langle \bu, \bu'_b\rangle = 0 \right\} \right)\times  H^1_\eta \times\R,\\
\Y_0&:= L^2_\eta\times \cZ,\\
\cZ&:=L^2_\eta \cap \left\{ \bu:\R\rightarrow\R~|~ \langle \bu, \mathbf{1} \rangle = 0 \right\}.
\end{align}
\end{subequations}
Finally we define the operator $\cL_f^i$ from $\X_f$ to $\Y_0$, with the entries of $\cL_f^i$ being given in equation \eqref{eq:LiFront} and $\cL_b^i$ from $\X_b$ to $\Y_0$, with the entries of $\cL_b^i$ being given in equation \eqref{eq:LiBack}.

\begin{lem}[Invertibility of the front and the back linearization]\label{lem:LiBackFront}The following assumptions hold true:
\begin{itemize}
\item[(i)] The operator $\cL_f^i:\X_f\longrightarrow \Y_0$ is invertible with bounded inverse, uniformly in $\epsilon>0$.
\item[(ii)] The operator $\cL_b^i:\X_b\longrightarrow \Y_0$ is invertible with bounded inverse, uniformly in $\epsilon>0$.
\end{itemize}
\end{lem}
\begin{Proof} 
\begin{itemize}
\item[(i)] We first remark that the operator $\cL(\bu_f)$ is a Fredholm operator from $H^1_\eta$ to $L^2_\eta$ whenever $\eta<\eta_0$, and its Fredholm index  is $0$ \cite{faye-scheel:13}. Its kernel is spanned by $\bu'_f$ and its cokernel is spanned by $\mathbf{e}_f^*\in H^1_\eta$, a solution of the adjoint equation 
\bqs
\cL^*(\bu_f)\,\mathbf{e}_f^*=0,
\eqs
where the adjoint operator is defined as
\bqs
\cL^*(\bu_f): \left\{\begin{array}{ccl}
H^1_\eta &\longrightarrow& L^2_\eta \\
 \bw&\longmapsto& -\bw+\K \ast \bw +  -c_* \frac{d}{d\xi}\bw +f'\left(\bu_f(\xi)\right)\bw.
 \end{array}\right.
\eqs
Second, we note that the operator $\dfrac{d}{d\xi}$ is Fredholm from $H^1_\eta$ to $L^2_\eta$, for all $\eta>0$, and its Fredholm index is $-1$ with cokernel spanned by $\mathbf{1}$ (the constants). Because of our specific choice of the target space $\Y_0$, we see that the Fredholm operator $\dfrac{d}{d\xi}$ is defined from $H^1_\eta$ to $\cZ$ and thus is Fredholm index  0 on $\cZ$. Finally, we notice that 
\bqs
\int _\R \bu_f'(\xi)\chi_f(\xi+\xi_f)\mathbf{e}_f^*(\xi)d\xi \to \int _\R \bu_f'(\xi)\mathbf{e}_f^*(\xi)d\xi \neq0,
\eqs
for $\epsilon\to 0$. Convergence is due to the fact that $\xi_f\to \infty$ as $\epsilon\to 0$. The latter integral is nonzero since the zero eigenvalue is algebraically simple by (H3). In summary, 
\begin{itemize}
\item[$\bullet$] $\cL(\bu_f)$ is Fredholm from $H^1_\eta\cap \left\{ \bu:\R\rightarrow\R~|~ \langle \bu, \bu'_f\rangle = 0 \right\} $ (orthogonal to the kernel $\bu_f'$) to $L^2_\eta$, with index $0$;
\item[$\bullet$] $\dfrac{d}{d\xi}$ is Fredholm from $H^1_\eta$ to $\cZ$, with index $0$;
\item[$\bullet$] $\int _\R \bu_f'(\xi)\chi_f(\xi+\xi_f)\mathbf{e}_f^*(\xi)d\xi \neq0$.
\end{itemize}
Thus $\cL_f^i:\X_f\longrightarrow \Y_0$ is invertible. The fact that its inverse is bounded uniformly in $\epsilon>0$ is straightforward from the explicit form of $\cL_f^i$ in \eqref{eq:LiFront}.
\item[(ii)] The exact same argument applies to the back and we have that:
\begin{itemize}
\item[$\bullet$] $\cL(\bu_b)$ is Fredholm from $H^1_\eta\cap \left\{ \bu:\R\rightarrow\R~|~ \langle \bu, \bu'_b\rangle = 0 \right\} $ (orthogonal to the kernel $\bu_b'$) to $L^2_\eta$, with index $0$;
\item[$\bullet$] $\dfrac{d}{d\xi}$ is Fredholm from $H^1_\eta$ to $\cZ$, with index $0$;
\item[$\bullet$] $-c'_b\int _\R \bu_b'(\xi)\chi_b(\xi+\xi_b)\mathbf{e}_b^*(\xi)d\xi \neq0$, where $\mathbf{e}_b^*$ spans the kernel of the adjoint operator $\cL^*(\bu_b)$. Note that $c_b' \neq 0$ based on our Hypothesis on the back and front solution.
\end{itemize}
We conclude that $\cL_b^i:\X_b\longrightarrow \Y_0$ is invertible. The fact that its inverse is bounded uniformly in $\epsilon>0$ again is immediate from the explicit form of $\cL_b^i$ in \eqref{eq:LiBack}.
\end{itemize}
\end{Proof}

\subsubsection{Invertibility of $\cL_q^i$ and $\cL_e^i$}

This section is devoted to establish the following result.
\begin{lem}[Invertibility of the quiescent and excitatory parts]\label{lem:LiQuiescentExcitatory} The linear operator associated to the quiescent (respectively excitatory) part $\cL_q^i$ (respectively $\cL_e^i$)  is invertible from $H^1_\eta\times H^1_\eta \times \R$ (respectively $\left(H^1_\eta\times H^1_\eta\right)^2 \times \R^2$) to $L^2_\eta \times L^2_\eta$ (respectively $\left(L^2_\eta \times L^2_\eta \right)^2$), with bounded inverse, uniformly in $\epsilon>0$. 
\end{lem}

\begin{Proof} 
The proof of the lemma  relies on Lemma \ref{lem:LqeInv} of Section \ref{subsec:LdeInv}, which ensures the existence of  $0<\eta_h<\eta_0$  such that for all $0<\eta<\eta_h$ the operators  $\cL(\bu_q)$ and $\cL(\bu_e)$ are both  isomorphisms from $H^1_\eta$ to $L^2_{\eta}$ and the norms $\| \cL(\bu_q)^{-1} \|$ and $\| \cL(\bu_e)^{-1} \|$ can be bounded independently of $\eta$ and $\epsilon$. Note that the shifted operator $\cL(\tau_{\xi_{be}}\cdot\bu_e)$ also satisfies the same properties.

\begin{itemize}
\item \underline{Quiescent:} To conclude the proof for the quiescent part, one needs to show that
\bqs
c_*\int_\R \partial_{v_q}\left( \bv_q(\epsilon\xi) \right)\chi'_q(\xi+\xi_q)d\xi\neq 0,
\eqs
as $\epsilon\rightarrow 0$. Indeed, $\chi'_q(\cdot+\xi_q)$ vanishes for all $|\xi|\geq 1$, so that the above integral simplifies to
\bqs
\int_\R \partial_{v_q}\left( \bv_q(\epsilon\xi) \right)\chi'_q(\xi+\xi_q)d\xi=\int_{-1}^1 \partial_{v_q}\left( \bv_q(\epsilon\xi) \right)\chi'_q(\xi+\xi_q)d\xi.
\eqs
For all $|\xi|\leq 1$, $\bv_q(\epsilon\xi)=v_b+\cO(\epsilon)$ as $\epsilon\rightarrow 0$, so that $\partial_{v_q}\left( \bv_q(\epsilon\xi) \right)\sim 1$ as $\epsilon\rightarrow 0$ and then
\bqs
c_*\int_\R \partial_{v_q}\left( \bv_q(\epsilon\xi) \right)\chi'_q(\xi+\xi_q)d\xi \sim -c_*\neq0. 
\eqs
This result combined with the fact the differential operator $\frac{d}{d\xi}$ from $H^1_\eta$ to $L^2_\eta$ is Fredholm, with index $-1$ and cokernel spanned by the constant $\mathbf{1}$, ensures that $\cL_q^i$ is invertible from $H^1_\eta\times H^1_\eta \times \R$ to $L^2_\eta \times L^2_\eta$, with bounded inverse, uniformly in $\epsilon>0$.
\item \underline{Excitatory:} To conclude the proof for the excitatory part, we need to check that the following integrals do not vanish for small $\epsilon$:
\begin{align*}
c_*\int_\R \partial_{v_{be}}(\bv_e(\epsilon\xi+\epsilon\xi_{be}))\chi'_e(\xi+\xi_{be})\mathds{1}_{\xi_{be}}(\xi+\xi_{be})d\xi&\neq0,\\
c_*\int_\R \partial_{v_{ef}}(\bv_e(\epsilon\xi))\chi'_e(\xi)\mathds{1}_{\xi_{ef}}(\xi)d\xi&\neq0.
\end{align*}
Once again, we use the definition of $\chi'_e(\cdot+\xi_{be})\mathds{1}_{\xi_{be}}(\cdot+\xi_{be})$ and $\chi'_e\mathds{1}_{\xi_{ef}}$ combined with the fact that, for all $|\xi|\leq 1$, $\bv_e(\epsilon\xi+\epsilon\xi_{be})=v_{be}+\cO(\epsilon)$ and $\bv_e(\epsilon\xi)=v_{ef}+\cO(\epsilon)$ as $\epsilon\rightarrow 0$. We then obtain
\begin{align*}
c_*\int_\R \partial_{v_{be}}(\bv_e(\epsilon\xi+\epsilon\xi_{be}))\chi'_e(\xi+\xi_{be})\mathds{1}_{\xi_{be}}(\xi+\xi_{be}) d\xi\sim c_* \int_\R\chi'_e(\xi)\mathds{1}_{\xi_{be}}(\xi)d\xi&\neq0,\\
c_*\int_\R \partial_{v_{ef}}(\bv_e(\epsilon\xi))\chi'_e(\xi)\mathds{1}_{\xi_{ef}}(\xi)d\xi\sim c_* \int_\R\chi'_e(\xi)\mathds{1}_{\xi_{ef}}(\xi)d\xi&\neq0.
\end{align*}
This result combined with the fact the differential operator $\frac{d}{d\xi}$ from $H^1_\eta$ to $L^2_\eta$ is Fredholm, with index $-1$ and cokernel spanned by the constant $\mathbf{1}$, ensures that $\cL_e^i$ is invertible from $\left(H^1_\eta\times H^1_\eta\right)^2 \times \R^2$ to $L^2_\eta\times L^2_\eta$, with bounded inverse, uniformly in $\epsilon>0$.
\end{itemize}
\end{Proof}

\subsection{Estimates for the cross-linear terms $\cL_\epsilon^p$}\label{subsec:Lp}

In this section, we provide some estimates for the cross-linear terms $\cL_\epsilon^p$ as $\epsilon \rightarrow 0$. We define the linear operator $\cL_\epsilon^p$ as follows
\bqs
\cL_\epsilon^p: \begin{array}{ccl}
\X \times \R^4  &\longrightarrow& \Y,
 \end{array}
\eqs
or, more explicitly, in matrix form as
\bqq
\label{eq:LpEps}
\cL_\epsilon^p=  \left( \T \circ \A_\bW^p(\epsilon)  \circ \T^{-1} ~ | ~ \A_\lambda^p(\epsilon) \right) ,
\eqq
where the shift operator $\T$ is defined as 
\bqs
\T\bW =  \left(\left(\tau_{\xi_q}\cdot\bw_q^u,\bw_q^v \right),\left(\tau_{\xi_b}\cdot\bw_b^u,\bw_b^v \right),\left(\tau_{\xi_{be}}\cdot \bw_{be}^u,\bw_{be}^v\right),\left(\tau_{\xi_{ef}}\cdot\bw_{ef}^u,\bw_{ef}^v \right),\left( \tau_{\xi_f}\cdot\bw_f^u,\bw_f^v\right) \right),
\eqs
and, for all $\xi\in\R$,
\bqs
\tau_{\xi_j}\cdot\bw_j^u(\xi)=\bw_j^u(\xi+\xi_j).
\eqs
The operators $\A_\bW^p(\epsilon)$ and $\A_\lambda^p(\epsilon) $ are defined in equations \eqref{eq:ApW} and \eqref{eq:ApLambda}, respectively, and are studied separately in the following two sections.

\subsubsection{Estimates for $\A_\bW^p(\epsilon)$}

The matrix of linear operators $\A_\bW^p(\epsilon): \X \rightarrow \Y$ defined in \eqref{eq:LpEps} is given by
\bqq
\label{eq:ApW}
\A_\bW^p(\epsilon)=\left(
\begin{matrix}
\cL_d  & 0 & \cL_{q,b}  & 0 & \cL_{q,be}  & 0 & \cL_{q,ef}  & 0 & \cL_{q,f}  & 0\\
\epsilon  & -\epsilon\gamma & \cL_{q,b}^u  & \cL_{q,b}^v & 0  & 0 & 0  & 0 & 0  & 0\\
\cL_{b,q}  & 0 & \cL_d  & 0 & \cL_{b,be}  & 0 & \cL_{b,ef}  & 0 & \cL_{b,f}  & 0\\
0  & 0 & \epsilon  & -\epsilon\gamma & 0  & 0 & 0  & 0 & 0  & 0\\
\cL_{be,q}  & 0 & \cL_{be,b}  & 0 & \cL_d  & 0 & \cL_{be,ef}  & 0 & \cL_{be,f} & 0\\
0  & 0 & 0  & 0 & \epsilon  & -\epsilon\gamma & 0  & 0 & 0  & 0\\
\cL_{ef,q}  & 0 & \cL_{ef,b}  & 0 & \cL_{ef,be}  & 0 & \cL_d  & 0 & \cL_{ef,f}  & 0\\
0  & 0 & 0  & 0 & 0  & 0 & \epsilon  & -\epsilon\gamma &  \cL_{ef,f}^u  & \cL_{ef,f}^v\\
\cL_{f,q}  & 0 & \cL_{f,b}  & 0 & \cL_{f,be}  & 0 & \cL_{f,ef}  & 0 & \cL_d  & 0\\
0  & 0 & 0  & 0 & 0  & 0 & 0  & 0 & \cL_f^u  & \cL_f^v
\end{matrix}
 \right).
\eqq
The operator $\cL_d$ is defined for all $\bw\in L^2_\eta$ as
\bqs
\cL_d\cdot\bw(\xi)=\bL_d(\xi)\bw(\xi),\quad \forall \xi \in \R,
\eqs
where $\bL_d$ is defined in equation \eqref{eq:Ldiag}. The different multiplication operators $\cL_{j,k}$ with $j,k\in \widetilde{J}_\bw$ are also implicitly defined using equations \eqref{eq:Lqj}, \eqref{eq:Lbj}, \eqref{eq:Lej} and \eqref{eq:Lfj} through
\bqs
\cL_{j,k}\cdot\bw(\xi)=\bL_{j,k}(\xi)\bw(\xi),\quad \forall \xi \in \R,\quad \forall \bw\in L^2_\eta.
\eqs
The remaining operators are defined as follows
\begin{subequations}
\label{eq:LuvJ}
\begin{align}
\cL _{q,b}^u\cdot\bw(\xi)&= \epsilon \left(\int_\R \bw(\xi)d\xi\right)\psi\left(\xi-\frac{\xi_b-\xi_q}{2}\right),\\
\cL _{q,b}^v\cdot\bw(\xi)&= -\epsilon\gamma \left(\int_\R \bw(\xi)d\xi\right)\psi\left(\xi-\frac{\xi_b-\xi_q}{2}\right),\\
\cL _{b}^u\cdot\bw(\xi)&= \epsilon \bw(\xi)-\epsilon \left(\int_\R \bw(\xi)d\xi\right)\psi\left(\xi+\frac{\xi_b-\xi_q}{2}\right),\\
\cL _{b}^v\cdot\bw(\xi)&= -\epsilon\gamma \bw(\xi)+\epsilon\gamma \left(\int_\R \bw(\xi)d\xi\right)\psi\left(\xi+\frac{\xi_b-\xi_q}{2}\right),\\
\cL _{ef,f}^u\cdot\bw(\xi)&= \epsilon \left(\int_\R \bw(\xi)d\xi\right)\psi\left(\xi-\frac{\xi_f}{2}\right),\\
\cL _{ef,f}^v\cdot\bw(\xi)&= -\epsilon\gamma \left(\int_\R \bw(\xi)d\xi\right)\psi\left(\xi-\frac{\xi_f}{2}\right),\\
\cL _{f}^u\cdot\bw(\xi)&= \epsilon \bw(\xi)-\epsilon \left(\int_\R \bw(\xi)d\xi\right)\psi\left(\xi+\frac{\xi_f}{2}\right),\\
\cL _{f}^v\cdot\bw(\xi)&= -\epsilon \gamma \bw(\xi)\epsilon\gamma \left(\int_\R \bw(\xi)d\xi\right)\psi\left(\xi+\frac{\xi_f}{2}\right).
\end{align}
\end{subequations}

We can immediately confirm that the following estimates hold for the above terms:
\begin{itemize}
\item $ \left\|  \cL_{q,b}^{u,v}\cdot\bw \right\| _{L^2_\eta} \lesssim \epsilon^{1-\frac{\eta\eta_b}{2}} \left\| \bw \right\| _{L^2_\eta} $ and $ \left\|  \cL_{b}^{u,v}\cdot\bw \right\| _{L^2_\eta} \lesssim \epsilon^{1-\frac{\eta\eta_b}{2}} \left\| \bw \right\| _{L^2_\eta} $;
\item $ \left\|  \cL_{ef,f}^{u,v}\cdot\bw \right\| _{L^2_\eta} \lesssim \epsilon^{1-\frac{\eta\eta_f}{2}} \left\| \bw \right\| _{L^2_\eta} $ and $ \left\|  \cL_{f}^{u,v}\cdot\bw \right\| _{L^2_\eta} \lesssim \epsilon^{1-\frac{\eta\eta_f}{2}} \left\| \bw \right\| _{L^2_\eta} $.
\end{itemize}
For the last two types of estimates we have used the estimates \eqref{eq:EstPsib} and \eqref{eq:EstPsif}. Our goal for this subsection is to prove the following result.
\begin{prop}\label{prop:ApWEps}
The cross-linear terms represented by $\A_\bW^p(\epsilon)$ are small in the operator norm,
\bqq
\label{eq:LimitApWEps}
\underset{ \epsilon \rightarrow 0}{\lim}~\left\| \A_\bW^p(\epsilon) \right\|=0.
\eqq
\end{prop}

In order to prepare the proof of Proposition \ref{prop:ApWEps}, we first notice that $\bL_d$ can be decomposed as follows
\bqq
\label{eq:LdiagD}
\bL_d(\xi)=\bD_q(\xi)\mathds{1}_{\xi_q}(\xi)+\bD_{be}(\xi)\mathds{1}_{\xi_{be}}(\xi)+\bD_{ef}(\xi)\mathds{1}_{\xi_{ef}}(\xi), \quad \forall\xi\in\R,
\eqq
where
\begin{subequations}
\label{eq:ComDj}
\begin{align}
\bD_q(\xi)&= f'\left((\bu_q\left(\epsilon\left(\xi-\xi_{q}\right)\right)\chi_q(\xi)+\bu_b(\xi-\xi_{b})\chi_b(\xi)\right)-f'(\bu_q\left(\epsilon\left(\xi-\xi_{q}\right)\right))\chi_q(\xi)-f'(\bu_b(\xi-\xi_{b}))\chi_b(\xi),\\
\bD_{be}(\xi)&=f'\left(\bu_b(\xi-\xi_{b})\chi_b(\xi)+\bu_e(\epsilon\xi)\chi_e(\xi) \right)-f'(\bu_b(\xi-\xi_{b}))\chi_b(\xi)-f'(\bu_e(\epsilon\xi))\chi_e(\xi),\\
\bD_{ef}(\xi)&=f'\left( \bu_e(\epsilon\xi)\chi_e(\xi)+\bu_f(\xi-\xi_f)\chi_f(\xi)\right)-f'(\bu_e(\epsilon\xi))\chi_e(\xi)-f'(\bu_f(\xi-\xi_f))\chi_f(\xi).
\end{align}
\end{subequations}

\begin{lem}\label{lem:Ldiag}
For all $j\in\widetilde{J}_\bw$ and for all $\bw \in L^2_\eta$, we have
\bqs
\left\| \bL_d(\cdot+\xi_j)\bw \right\|_{L^2_\eta}\leq C(\epsilon) \left\| \bw \right\|_{L^2_\eta}, \quad C(\epsilon) \longrightarrow 0,
\eqs
as $\epsilon \rightarrow 0$.
\end{lem}

\begin{Proof}
Using the same arguments as in Proposition \ref{prop:EstCom}, we directly have that
\begin{align*}
\left\| \bD_q\mathds{1}_{\xi_q}\right\|_{L^\infty}&=\cO(\epsilon^{2\min(1,\eta_*\eta_b)}),\\
\left\| \bD_{be}\mathds{1}_{\xi_{be}} \right\|_{L^\infty}&=\cO(\epsilon^{2\min(1,\eta_*\eta_b)}),\\
\left\| \bD_{ef} \mathds{1}_{\xi_{ef}} \right\|_{L^\infty}&=\cO(\epsilon^{2\min(1,\eta_*\eta_f)}).
\end{align*}
Then, for all $j\in\widetilde{J}_\bw$
\bqs
C(\epsilon)=\left\| \bL_d(\cdot+\xi_j) \right\|_{L^\infty}=\cO(\epsilon^{2\min(1,\eta_*\eta_b)})+\cO(\epsilon^{2\min(1,\eta_*\eta_f)})\longrightarrow 0 \text{ as } \epsilon \rightarrow 0.
\eqs
\end{Proof}

We are now ready to give the proof of Proposition \ref{prop:ApWEps}.
\begin{Proof}[ of Proposition \ref{prop:ApWEps}.]
We will only give the proof for the last two components of $\cL_\epsilon^p$, the other components being treated in the exact same way. For all $(\bW,\lambda)\in \X\times\R^4$, the first two components of $\cL_\epsilon^p(\bW,\lambda)$ are given by
\begin{align*}
\left[\cL_\epsilon^p(\bW,\lambda)\right]_f^u(\xi)=&~\bL_d(\xi+\xi_f)\bw_f^u(\xi)+\bL_{f,q}(\xi+\xi_f)\bw_q^u(\xi-\xi_q+\xi_f)+\bL_{f,b}(\xi+\xi_f)\bw_{b}^u(\xi-\xi_{b}+\xi_f)\\
&+\bL_{f,be}(\xi+\xi_f)\bw_{be}^u(\xi-\xi_{be}+\xi_f)+\bL_{f,ef}(\xi+\xi_f)\bw_{ef}^u(\xi+\xi_f),\\
\left[\cL_\epsilon^p(\bW,\lambda)\right]_f^v(\xi)=&~\epsilon \left(\int_\R \left(\bw_f^u(\xi)-\gamma\bw_f^v(\xi)\right) d\xi\right)\psi\left(\xi-\frac{\xi_f}{2}\right).
\end{align*}
Using our previous estimates, we directly have
\bqs
\left\| \left[\cL_\epsilon^p(\bW,\lambda)\right]_f^v \right\|_{L^2_\eta} \lesssim \epsilon^{1-\frac{\eta\eta_f}{2} }\left\| \bW \right\|_\X.
\eqs
We next treat each term of the first component separately. The first term $\bL_d(\xi+\xi_f)\bw_f^u(\xi)$ has been considered in Lemma \ref{lem:Ldiag} and we have
\bqs
\left\| \bL_d(\cdot+\xi_f)\bw_f^u \right\|_{L^2_\eta}\lesssim \epsilon^{2\min(1,\eta_*\eta_b,\eta_*,\eta_f)} \left\| \bw _f^u \right\|_{L^2_\eta}.
\eqs
We now show that
\begin{itemize}
\item[(i)] $\left\|\bL_{f,j}\left(\cdot+\xi_f \right)\bw_j^u(\cdot-\xi_j+\xi_f)\right\|_{ L^2_\eta} \lesssim e^{-\eta\frac{T(v_{be},v_{ef})}{\epsilon}} \left\| \bw _j^u\right\|_{ L^2_\eta}$, for all $j\in\left\{ q,b,be \right\}$;
\item[(ii)] $\left\|\bL_{f,ef}\left(\cdot+\xi_f \right)\bw_{ef}^u(\cdot+\xi_f)\right\|_{ L^2_\eta} \lesssim \left(\epsilon^{\eta_f \min\left((\eta_*-\eta),\eta\right)} + \epsilon^{1- \eta\eta_f }|\ln \epsilon| \right)\left\| \bw_{ef}^u \right\|_{ L^2_\eta}$.
\end{itemize}
\underline{Estimate (i):} As $\bw_j^u \in L^2_\eta$, we have that $\xi \mapsto e^{\eta|\xi-\xi_j+\xi_f|}\bw_j^u(\xi-\xi_j+\xi_f) $ belongs to $L^2$ with 
\bqs
\left\| \bw_j^u \right\|_{ L^2_\eta}=\left\| e^{\eta|\cdot-\xi_j+\xi_f|}\bw_j^u(\cdot -\xi_j+\xi_f) \right\|_{ L^2}.
\eqs
Then we have for all $j \in \left\{ q,b,be \right\}$,
\begin{align*}
\left\|\bL_{f,j}\left(\cdot+\xi_f \right)\bw_j^u(\cdot-\xi_j+\xi_f)\right\|_{ L^2_\eta}&=\left\|e^{\eta|\cdot|}\bL_{f,j}\left(\cdot+\xi_f \right)e^{-\eta|\cdot-\xi_j+\xi_f|}e^{\eta|\cdot-\xi_j+\xi_f|}\bw_j^u(\cdot-\xi_j+\xi_f)\right\|_{ L^2}\\
&\leq \left\|e^{\eta|\cdot|}\bL_{f,j}\left(\cdot+\xi_f \right)e^{-\eta|\cdot-\xi_j+\xi_f|}\right\|_{ L^\infty(\R)} \left\| \bw_j^u \right\|_{ L^2_\eta}.
\end{align*}
Using the definition of $\bL_{f,j}\left(\cdot+\xi_f \right)$ and the property of $\chi_f$, we find that
\bqs
\left\|e^{\eta|\cdot|}\bL_{f,j}\left(\cdot+\xi_f \right)e^{-\eta|\cdot-\xi_j+\xi_f|}\right\|_{ L^\infty} \leq C \underset{\xi+\xi_f\geq -1 }{\sup}~ e^{\eta \left(|\xi|- |\xi-\xi_j+\xi_f| \right)}.
\eqs
As $\epsilon\rightarrow 0$, the $\sup$ is obtained when $\xi \sim - \xi_f $ and we have for all $j \in \left\{ q,b,be\right\}$
\bqs
\left\|e^{\eta|\cdot|}\bL_{f,j}\left(\cdot+\xi_f\right)e^{-\eta|\cdot-\xi_j+\xi_f|}\right\|_{ L^\infty} \leq C_0 e^{-\eta \frac{T(v_{be},v_{ef})}{\epsilon}}.
\eqs
\underline{Estimate (ii):} For the second estimate, we recall that
\bqs
\bL_{f,ef}(\xi+\xi_f)=\chi_f(\xi+\xi_f)\left( f'(\bu_f(\xi))-f'(\bu_e(\epsilon\xi+\epsilon\xi_f)) \right)
\eqs
for all $\xi\in\R$. Thus, we need to evaluate
\bqs
\underset{\xi\in\R}{\sup}~ \left| e^{\eta\left(|\xi|-|\xi+\xi_f|\right)}\chi_f(\xi+\xi_f)\left( f'(\bu_f(\xi))-f'(\bu_e(\epsilon\xi+\epsilon\xi_f)) \right) \right|,
\eqs
when $\epsilon \rightarrow 0$. It is not difficult to see that this $\sup$ is realized for values of $\xi$ in $[-\xi_f,0]$. We set $\xi=\xi_1\eta_f \ln\epsilon$ for $\xi_1 \in [0,1]$ and look for
\bqs
\underset{\xi_1\in[0,1]}{\sup}~ \left| \epsilon^{(1-2\xi_1)\eta\eta_f}\left(f'\left(\bu_f(\xi_1\eta_f\ln\epsilon) \right)-f'\left(\bu_e((\xi_1-1)\eta_f \epsilon \ln\epsilon)\right) \right)\right|.
\eqs
We have the asymptotic estimates as $\epsilon \rightarrow 0$ for all $\xi \in [0,1)$
\begin{align*}
\bu_f(\xi_1 \eta_f \ln\epsilon)-1&=\cO\left( \epsilon^{\xi_1\eta_f\eta_*}\right) \\
\bu_e((\xi_1-1)\eta_f \epsilon \ln\epsilon)&-1= \cO(\epsilon \ln\epsilon).
\end{align*}
Then,
\bqs
\underset{\xi_1\in[0,1]}{\sup}~ \left| \epsilon^{(1-2\xi_1)\eta\eta_f}\left( \epsilon^{ \xi_1\eta_f\eta_*} + \epsilon \ln\epsilon \right) \right| \leq \epsilon^{\eta_f \min\left((\eta_*-\eta),\eta\right)} + \epsilon^{1- \eta\eta_f }|\ln \epsilon|.
\eqs
Regrouping all our estimates, we have shown that
\bqs
\left\| \left[\cL_\epsilon^p(\bW,\lambda)\right]_f^u \right\|_{L^2_\eta} \lesssim C(\epsilon) \left\| \bW \right\|_\X,
\eqs
with $C(\epsilon)\rightarrow0$ as $\epsilon\rightarrow0$. This concludes the proof.
\end{Proof}

\subsubsection{Estimates for $\A_\lambda^p(\epsilon)$}
The matrix of linear operators $\A_\lambda^p(\epsilon): \R^4 \rightarrow \Y$ defined in \eqref{eq:LpEps} is given by
\bqq
\label{eq:ApLambda}
\A_\lambda^p(\epsilon)= \left( \bA^{p,1}_\lambda(\epsilon) ~ \bA_\lambda^{p,2}(\epsilon) ~ \bA_\lambda^{p,3}(\epsilon) ~ \bA_\lambda^{p,4}(\epsilon)~ \bA_\lambda^{p,5}(\epsilon) \right)
\eqq
where the columns $\bA^{p,j}_\lambda(\epsilon)$ are defined as
\bqs
\bA^{p,1}_\lambda(\epsilon) =\left(
\begin{matrix}
\bu_q(\epsilon\cdot)\chi'_q(\cdot+\xi_q)+\bu_b(\cdot-\xi_b+\xi_q)\chi'_b(\cdot+\xi_q)\mathds{1}_{\xi_b}(\xi+\xi_q) \\
\bv_q(\epsilon\cdot)\chi'_q(\cdot+\xi_q)+v_*\chi'_b(\cdot+\xi_q)\mathds{1}_{\xi_b}(\xi+\xi_q) \\
0 \\
0 \\
\bu_b(\cdot-\xi_b+\xi_{be})\chi'_b(\cdot+\xi_{be})\mathds{1}_{\xi_{be}}(\cdot+\xi_{be}) +\bu_e(\epsilon\cdot+\epsilon\xi_{be})\chi'_e(\cdot+\xi_{be})\mathds{1}_{\xi_{be}}(\cdot+\xi_{be})\\
v_*\chi'_b(\cdot+\xi_{be})\mathds{1}_{\xi_{be}}(\cdot+\xi_{be}) +\bv_e(\epsilon\cdot+\epsilon\xi_{be})\chi'_e(\cdot+\xi_{be})\mathds{1}_{\xi_{be}}(\cdot+\xi_{be}) \\
\bu_e(\epsilon\cdot)\chi'_e\mathds{1}_{\xi_{ef}}+\bu_f(\cdot-\xi_f)\chi'_f \\
\bv_e(\epsilon\cdot)\chi'_e\mathds{1}_{\xi_{ef}} \\
0 \\
0 \\
\end{matrix}
 \right),
\eqs

\bqs
\bA^{p,2}_\lambda(\epsilon) =\left( 
\begin{matrix}
\partial_{v_q}\left( \tau_{\xi_q}\cdot\left[\bF_q\mathds{1}_{\xi_q}+\bC_q\right] \right)\\
\epsilon \partial_{v_q}\left(\chi_q(\cdot+\xi_q)\left(\bu_q(\epsilon \cdot)-\gamma\bv_q(\epsilon \cdot)\right)\left(1-\Theta(\bv_q(\epsilon \cdot))\right) \right)\\
0\\
0\\
0\\
0\\
0\\
0\\
0\\
0
\end{matrix}
\right),
\eqs
\bqs
\bA^{p,3}_\lambda(\epsilon) =\left( 
\begin{matrix}
 \partial_{v_b}\left(\tau_{\xi_{q}}\cdot \left[\bF_q\mathds{1}_{\xi_q}+\bC_q\right]\right)\\
\epsilon \partial_{v_b}\left(\int_\R (\bu_b(\xi)-\gamma v_b)\chi_b(\xi+\xi_b)d\xi\right)\psi\left(\xi-\frac{\xi_b-\xi_q}{2}\right)\\
0\\
\epsilon(\partial_{v_b}\bu_b-\gamma )\chi_b(\cdot+\xi_b)-\epsilon \partial_{v_b}\left(\int_\R (\bu_b(\xi)-\gamma v_b)\chi_b(\xi+\xi_b)d\xi\right)\psi\left(\xi-\frac{\xi_b-\xi_q}{2}\right)\\
\partial_{v_b}\left(\tau_{\xi_{be}}\cdot \left[\bF_{be}\mathds{1}_{\xi_{be}}+\bC_{be}\right]\right)\\
0\\
0\\
0\\
0\\
0
\end{matrix}
\right),
\eqs
\bqs
\bA^{p,4}_\lambda(\epsilon) =\left( 
\begin{matrix}
0\\
0\\
0\\
0\\
\partial_{v_{be}}\left(\tau_{\xi_{be}}\cdot \left[\bF_{be}\mathds{1}_{\xi_{be}}+\bC_{be}\right]\right)\\
\epsilon \partial_{v_{be}}\left( \tau_{\xi_{be}}\cdot\left[ \chi_e\left(\bu_e(\epsilon\cdot)-\gamma\bv_e(\epsilon \cdot)\right)\left(1-\Theta(\bv_e(\epsilon \cdot))\right)\mathds{1}_{\xi_{be}} \right] \right)\\
\partial_{v_{be}}\left(\bF_{ef}\mathds{1}_{\xi_{ef}}+\bC_{ef}\right)\\
\epsilon \partial_{v_{be}}\left( \chi_e\left(\bu_e(\epsilon\cdot)-\gamma\bv_e(\epsilon \cdot)\right)\left(1-\Theta(\bv_e(\epsilon \cdot))\right)\mathds{1}_{\xi_{ef}}  \right)\\
0\\
0\\
\end{matrix}
\right),
\eqs
and
\bqs
\bA^{p,5}_\lambda(\epsilon) =\left( 
\begin{matrix}
0\\
0\\
0\\
0\\
\partial_{v_{ef}}\left(\tau_{\xi_{be}}\cdot \left[\bF_{be}\mathds{1}_{\xi_{be}}+\bC_{be}\right]\right)\\
\epsilon \partial_{v_{ef}}\left( \tau_{\xi_{be}}\cdot\left[ \chi_e\left(\bu_e(\epsilon\cdot)-\gamma\bv_e(\epsilon \cdot)\right)\left(1-\Theta(\bv_e(\epsilon \cdot))\right)\mathds{1}_{\xi_{be}} \right] \right)\\
\partial_{v_{ef}}\left(\bF_{ef}\mathds{1}_{\xi_{ef}}+\bC_{ef}\right)\\
\epsilon \partial_{v_{ef}}\left( \chi_e\left(\bu_e(\epsilon\cdot)-\gamma\bv_e(\epsilon \cdot)\right)\left(1-\Theta(\bv_e(\epsilon \cdot))\right)\mathds{1}_{\xi_{ef}}  \right)\\
0\\
0\\
\end{matrix}
\right).
\eqs

\begin{prop}\label{prop:ApLambdaEps}
The following limit holds true in operator norm
\bqq
\label{eq:LimitApLambdaEps}
\underset{ \epsilon \rightarrow 0}{\lim}~\left\| \A_\lambda^p(\epsilon) \right\|=0.
\eqq
\end{prop}
\begin{Proof}
The proof follows closely the computations developed in Propositions \ref{prop:EstChi}, \ref{prop:EstCom} and \ref{prop:EstErr}. 
\begin{itemize}
\item The fact that 
\bqs
\underset{ \epsilon \rightarrow 0}{\lim}~\left\| \bA_\lambda^{p,1}(\epsilon) \right\|=0,
\eqs
is a direct consequence of Proposition \ref{prop:EstChi}.
\item For the commutators terms, let us show that
\bqs
\underset{ \epsilon \rightarrow 0}{\lim}~\left\| \partial_{v_{ef}}(\bF_{ef}\mathds{1}_{ef}) \right\|_{L^2_\eta}=0.
\eqs
The proofs for the other $\bF_j\mathds{1}_j$'s and $\bC_j$'s are analogous. From the definition of $\bF_{ef}$, we see that for all $|\xi| \leq 1$,
\bqs
\partial_{v_{ef}}\left(\bF_{ef}(\xi)\right)=\partial_{v_{ef}}\left(\bu_e(\epsilon\xi)\right)\left[f'\left(\bu_e(\epsilon\xi)\chi_e(\xi)+\bu_f(\xi-\xi_f)\chi_f(\xi) \right)-f'(\bu_e(\epsilon\xi)\chi_e(\xi))\right].
\eqs
We can Taylor expand $f'$ at $1$ and  obtain:
\begin{align*}
f'\left(\bu_e(\epsilon\xi)\chi_e(\xi)+\bu_f(\xi-\xi_f)\chi_f(\xi)\right)=&~f'(1)+f''(1) \left[\left(\bu_e(\epsilon\xi)-1\right)\chi_e(\xi)+\left(\bu_f(\xi-\xi_f)-1\right)\chi_f(\xi)\right]\\
&+\mathbf{R}_\epsilon(\xi),
\end{align*}
where
\begin{align*}
\mathbf{R}_\epsilon(\xi)&=\left(\bu_{ef}(\xi)\right)^2\int_0^1f''\left( 1+\tau \bu_{ef}(\xi) \right)(1-\tau)d\tau,\\
\bu_{ef}(\xi)&=\left(\bu_e(\epsilon\xi)-1\right)\chi_e(\xi)+\left(\bu_f(\xi-\xi_f)-1\right)\chi_f(\xi).
\end{align*}
Similarly, we have
\bqs
f'\left(\bu_e(\epsilon\xi)\right)=f'(1)+f''(1) \left(\bu_e(\epsilon\xi)-1\right)+\left(\bu_e(\epsilon\xi)-1\right)^2\int_0^1f''\left( 1+\tau \left(\bu_e(\epsilon\xi)-1 \right) \right)(1-\tau)d\tau.
\eqs
Finally using the asymptotic expansions for $\bu_e$ and $\bu_f$, we obtain the following estimate
\bqs
\left\| \partial_{v_{ef}}(\bF_{ef} \mathds{1}_{\xi_{ef}}) \right\|_{L^2_\eta} \lesssim \left(\epsilon^{\eta_*\eta_f}+\epsilon^2+\epsilon^{2\eta_*\eta_f}\right)\left\| \mathds{1}_{\xi_{ef}} \right\|_{L^2_\eta},
\eqs
which gives the result.
\item A direct computation shows that
\bqs
\left\| \epsilon \partial_{v_b}\left(\int_\R (\bu_b(\xi)-\gamma v_b)\chi_b(\xi+\xi_b)d\xi\right)\psi\left(\cdot-\frac{\xi_b-\xi_q}{2}\right) \right\| _{L^2_\eta}\sim C \epsilon^{1-\frac{\eta_*\eta_b}{2}}|\ln\epsilon|,
\eqs
as $\epsilon\rightarrow 0$, where the constant $C>0$ is given by
\bqs
C=\eta_b \left| \frac{1}{f'(\varphi_q(v_*))}+\frac{1}{f'(\varphi_e(v_*))}-2\gamma \right|.
\eqs
\item All the remaining terms can be analyzed using Proposition \ref{prop:EstErr}.
\end{itemize}
\end{Proof}

\subsection{Conclusion of the proof of Theorem \ref{thm:existence}}\label{subsec:ProofEx}

In this section, we gather all the information collected so far and prove Theorem  \ref{thm:existence}. We recall that we want to prove the existence of $\left(\bW(\epsilon),\lambda(\epsilon) \right)$, for $0<\epsilon<\epsilon_0$,  solution of the equation
\bqs
\F_\epsilon(\bW,\lambda)=0,
\eqs
where $\F_\epsilon$ is defined in equation \eqref{eq:Compact} as the collection of systems \eqref{eq:Wq}, \eqref{eq:Wb}, \eqref{eq:Wbe}, \eqref{eq:Wef} and \eqref{eq:Wf}, with
\begin{align*}
\bW&:=\left(\left(\bw_q^u,\bw_q^v \right),\left(\bw_b^u,\bw_b^v \right),\left( \bw_{be}^u,\bw_{be}^v\right),\left(\bw_{ef}^u,\bw_{ef}^v \right),\left( \bw_f^u,\bw_f^v\right) \right),\\
\lambda&:=\left(c,v_q,v_b,v_{be},v_{ef} \right).
\end{align*}
Based on the analysis of the previous sections, we define two new Banach spaces $\X_*$ and $\Y_*$ through
\begin{align*}
\X_*&:=\left(H^1_\eta \times H^1_\eta \right)\times \left(\X_b \times H^1_\eta \right) \times \left(H^1_\eta \times H^1_\eta \right) \times \left(H^1_\eta \times H^1_\eta \right) \times  \left(\X_f \times H^1_\eta \right),\\
\Y_*&:=\left(L^2_\eta \times L^2_\eta \right)\times \left(L^2_\eta \times \cZ \right) \times \left(L^2_\eta \times L^2_\eta \right) \times \left(L^2_\eta \times L^2_\eta \right) \times  \left(L^2_\eta \times \cZ \right),
\end{align*}
where $\X_b$, $\X_f$ and $\cZ$ have been defined in \eqref{eq:BanachSpace}. Note that the map $\F_\epsilon$ is well-defined from $\X_*\times \V$ into $\Y_*$ where $\V=(c_*-\delta_c,c_*+\delta_c)\times (v_*-\delta_b,v_*+\delta_b)\times (v_*-\delta_b,v_*+\delta_b)\times (v_*-\delta_{be},v_*+\delta_{be})\times (-\delta_{ef},\delta_{ef})$ is a neighborhood of $\lambda_*=(c_*,v_*,v_*,v_*,0)$ in $\R^5$. Note that the zero-mass conditions encoded in the space $\cZ$ is satisfied due to our particular choice of mass distribution via $\psi$ in  \eqref{eq:Wbv} and \eqref{eq:Wfv}. Using the fact that $\F_\epsilon$ is $\cC^\infty$ in its two arguments, we directly see that
\bqs
\F_\epsilon(\bW,\lambda)=\cR_\epsilon+\cL_\epsilon(\bW,\lambda-\lambda_*)+\cN_\epsilon(\bW,\lambda-\lambda_*),
\eqs
where
\begin{align*}
\cR_\epsilon&=\F_\epsilon(\mathbf{0},\lambda_*),\\
\cL_\epsilon(\bW,\lambda-\lambda_*)&=D_\bW\F_\epsilon(\mathbf{0},\lambda_*)\bW+D_\lambda\F_\epsilon(\mathbf{0},\lambda_*)(\lambda-\lambda_*),\\
\cN_\epsilon(\bW,\lambda-\lambda_*)&=\F_\epsilon(\bW,\lambda)-\F_\epsilon(\mathbf{0},\lambda_*)-D_\bW\F_\epsilon(\mathbf{0},\lambda_*)\bW-D_\bW\F_\epsilon(\mathbf{0},\lambda_*)(\lambda-\lambda_*).
\end{align*}
The error term $\cR_\epsilon$ has been defined in equations  \eqref{eq:Compact} and  \eqref{eq:DefRq}, \eqref{eq:DefRb}, \eqref{eq:DefRe}, \eqref{eq:DefRf} and satisfies the limit  
\bqs
\underset{\epsilon \rightarrow 0}{\lim}\;\left\| \mathcal{R}_\epsilon \right\|_{\Y_*}=0,
\eqs
as proved in Proposition \ref{prop:EstR}. The linear part $\cL_\epsilon$ can be decomposed into two parts: an invertible part with bounded inverse on $\X_*$ and a perturbation part that converges to zero as $\epsilon\rightarrow0$ in operator norm. More precisely, we have that
\bqs
\cL_\epsilon=\cL_\epsilon^i+\cL_\epsilon^p,
\eqs 
where $\cL_\epsilon^i$ is defined in equation \eqref{eq:LiEps} and $\cL_\epsilon^p$ in equation \eqref{eq:LpEps}. Lemma \ref{lem:LiBackFront} and  \ref{lem:LiQuiescentExcitatory} combined show that $\cL_\epsilon^i:\X_*\times \R^5\rightarrow \Y_*$ is invertible with inverse bounded independent of $\epsilon>0$. That is, there exists $M>0$, independent of $\epsilon>0$, so that
\bqs
\left\| \left(\cL_\epsilon^i\right)^{-1}\right\|\leq M.
\eqs
Using Proposition \ref{prop:ApWEps} and \ref{prop:ApLambdaEps}, we have the following limit for the perturbation $\cL_\epsilon^p$
\bqs
\underset{\epsilon\rightarrow0}{\lim}~\left\| \cL_\epsilon^p\right\|=0.
\eqs
Then a perturbation argument ensures that, for $\epsilon$ small, $\cL_\epsilon:\X_*\times \R^5\rightarrow \Y_*$ is invertible with inverse bounded independent of $\epsilon>0$. Finally, the nonlinear term $\cN_\epsilon$ is quadratic in $\bW$ and $\lambda-\lambda_*$, so that we have
\bqs
\left\|\cN_\epsilon(\bW,\lambda-\lambda_*)\right\|_{\Y_*}=\cO(\left\| \bW \right\|^2_{\X_*}+\left\| \lambda-\lambda_*\right\|^2_{\R^5}), \text{ as } (\bW,\lambda)\rightarrow (\mathbf{0},\lambda_*).
\eqs
Note that the quadratic terms in $\bW$ appear explicitly in systems \eqref{eq:Wq}, \eqref{eq:Wb}, \eqref{eq:Wbe}, \eqref{eq:Wef} and \eqref{eq:Wf} while the nonlinear terms in $\lambda-\lambda_*$ are only defined implicitly through equation \eqref{eq:Compact}.

We are now ready to use a fixed point iteration argument on the map $\F_\epsilon$ which will give us the existence of $\left(\bW(\epsilon),\lambda(\epsilon) \right)$ solution of equation \eqref{eq:Compact} in a neighborhood of $(\mathbf{0},\lambda_*)$ for small values of $\epsilon>0$. As for the proof of Proposition \ref{prop:slowsolL} and \ref{prop:slowsolR}, we introduce a map $\cS_\epsilon:\X_*\times \U \rightarrow \Y_*$ defined as
\bqs
\cS_\epsilon(\bW,\rho)=(\bW,\rho)-\cL_\epsilon^{-1}\left(\F_\epsilon(\bW,\rho+\lambda_*)\right),
\eqs
where $\U\subset \R^5$ is a neighborhood of $(0,0,0,0,0)$.
Based on the conclusions stated above, the map $\cS_\epsilon$ satisfies the following properties:
\begin{itemize}
\item $\left\|\cS_\epsilon(\mathbf{0},\mathbf{0})\right\|_{\Y_*}\rightarrow 0$ as $\epsilon \rightarrow 0$;
\item $\cS_\epsilon$ is a $\cC^\infty$-map;
\item $D_{(\bW,\rho)}\cS_\epsilon(\mathbf{0},\mathbf{0})=0$;
\item there exist $\delta>0$ and $C_1>0$ such that for all $(\bW,\rho)\in\B_\delta$, the ball of radius $\delta$ centered at $(\mathbf{0},\mathbf{0})$ in $\X_*\times \U$, we have
\bqs
\left\| D_{(\bW,\rho)}\cS_\epsilon(\bW,\rho) \right\| \leq C_1 \delta.
\eqs
\end{itemize}

We can now define an iteration scheme as follows
\bqs
(\bW_{n+1},\rho_{n+1})=S_\epsilon(\bW_n,\rho_n)=(\bW_n,\rho_n)-\cL_\epsilon^{-1}\left(\F_\epsilon(\bW_n,\rho_n+\lambda_*)\right)
, \quad n\geq 0,
\eqs
with initial point $(\bW_0,\rho_0)=(\mathbf{0},\mathbf{0})$.
Suppose, by induction, that $(\bW_k,\rho_k)\in\B_\delta$ for all $1\leq k \leq n$, then
\bqs
\left\| (\bW_{n+1},\rho_{n+1})-(\bW_n,\rho_n) \right\| \leq C_1 \delta \left\| (\bW_{n},\rho_{n})-(\bW_{n-1},\rho_{n-1}) \right\|,
\eqs
so that
\bqs
\left\| (\bW_{n+1},\rho_{n+1}) \right\| \leq \frac{C_0}{1-C_1 \delta}\epsilon.
\eqs
For small enough $\epsilon$, we have $\dfrac{C_0}{1-C_1 \delta}\epsilon<\delta$ and $(\bW_{n+1},\rho_{n+1})\in\B_\delta$ so that we have a contraction. We can then apply the Banach's fixed point theorem to find a solution $(\bW(\epsilon),\rho(\epsilon))=\underset{n\rightarrow\infty}{\lim~}(\bW_n,\rho_n)$ such that $(\bW(\epsilon),\rho(\epsilon))=\cS_\epsilon(\bW(\epsilon),\rho(\epsilon))$. As a conclusion, for every sufficiently small $\epsilon>0$, we have proved the existence of a traveling pulse solution to \eqref{eq:TP}.

\bibliography{plain}

\end{document}